\documentclass{amsart}
\usepackage[utf8]{inputenc}
\usepackage{lmodern}
\usepackage[margin=3cm]{geometry}
\usepackage{extdash}
\usepackage{amsmath,amsthm,amsfonts,amssymb, microtype}
\usepackage{xr-hyper}
\usepackage[colorlinks=true]{hyperref}
\usepackage{verbatim}
\usepackage{pdfsync}
\usepackage{stmaryrd}
\usepackage{marvosym}
\usepackage{colonequals}
\usepackage{mathrsfs}

\usepackage{mathtools}
\usepackage{enumitem}
\usepackage{tikz}
\usetikzlibrary{arrows,calc,positioning}

\usetikzlibrary{backgrounds}

\usepackage{extdash}
\usepackage{amsthm,amsfonts,amssymb,amsmath, microtype}
\usepackage{xr-hyper}
\usepackage[colorlinks=true]{hyperref}
\usepackage{verbatim}
\usepackage{pdfsync}
\usepackage{comment} 

\usetikzlibrary{circuits.logic.IEC}
\usetikzlibrary{cd}
\usepackage{xcolor}
\usetikzlibrary{shapes}
\usetikzlibrary{decorations.markings}

\tikzstyle heightone=[scale=.7,shift={(0,-.3)}]
\tikzstyle heightones=[scale=.8,xscale=.35,shift={(0,.1)}]
\tikzstyle heightoneonehalf=[scale=.9,shift={(0,-.2)}]
\tikzstyle heighttwo=[scale=.9,shift={(0,-.4)}]
\tikzstyle heighttwos=[scale=.5,xscale=.6,shift={(0,-.1)}]
\tikzstyle heightthree=[scale=.6,shift={(0,-.9)}]
\tikzstyle heightthrees=[scale=.4,xscale=.7,shift={(0,-.2)}]

\tikzstyle arrowstyle=[blue,semitransparent,scale=2]

\tikzstyle basiclabel=[draw=none,fill=none,shape=rectangle,inner sep=2pt,scale=.8]
\tikzstyle leftlabel=[basiclabel,anchor=east]
\tikzstyle rightlabel=[basiclabel,anchor=west]
\tikzstyle bottomlabel=[basiclabel,anchor=north]
\tikzstyle toplabel=[basiclabel,anchor=south]

\tikzstyle vertex=[circle,draw,fill=black,inner sep=1pt]
\tikzstyle ciliation=[circle,draw=none,fill=red,inner sep=1pt,semitransparent]
\tikzstyle ciliatednode=[vertex,pin={[pin distance=1mm,pin edge={semitransparent,red},ciliation]#1:{}}]

\tikzstyle vector=[black,thick,rectangle,draw=gray!50!yellow,top color=yellow!30,bottom color=black!10,scale=.8,inner sep=2pt]
\tikzstyle small vector=[vector,scale=.8]
\tikzstyle plain vector=[rectangle,draw=none,fill=white,scale=.7]

\tikzstyle my signal=[black,thick,signal,signal pointer angle=120,draw=blue!50,top color=blue!20,bottom color=black!10,scale=.8,inner sep=2pt]
\tikzstyle matrix=[my signal,signal from=south,signal to=north]
\tikzstyle reverse matrix=[my signal,signal from=north,signal to=south]

\tikzstyle small matrix=[matrix,scale=.7]
\tikzstyle reverse small matrix=[reverse matrix,scale=.7]
\tikzstyle matrix on edge=[small matrix,sloped,rotate=-90]
\tikzstyle reverse matrix on edge=[small matrix,sloped,rotate=90]

\tikzstyle trivalent=[very thick]
\tikzstyle dotdotdot=[decorate,decoration={markings,
    mark=at position .3 with{\node{.};},
    mark=at position .5 with {\node{.};},
    mark=at position .7 with {\node{.};}}]

\tikzstyle wavyup=[out=90,in=-90]
\tikzstyle wavydown=[out=-90,in=90]

\tikzstyle symmetrizer=[rectangle,fill=gray!10,draw=black]
\tikzstyle permutation=[symmetrizer]
\tikzstyle antisymmetrizer=[rectangle,fill=black,draw=black]
\tikzstyle symlabel=[draw=none,fill=none,black,scale=.8]
\tikzstyle asymlabel=[draw=none,fill=none,white,scale=.8]

\newcommand{\PPP}{\mathcal{P}}
\newcommand{\I}{\mathrm{I}}
\newcommand{\Z}{\mathbb{Z}}

\newcommand{\N}{\mathbb{N}}

\newcommand{\obj}{\mathrm{obj}}

\newcommand{\id}{\mathrm{Id}}

\newcommand{\op}{\mathrm{op}}
\newcommand{\h}{\mathrm{h}}
\newcommand{\HHH}{\mathrm{HH}}
\newcommand{\Ext}{\mathrm{Ext}}

\newcommand{\Qch}{\mathrm{Qch}}
\newcommand{\Coh}{\mathrm{coh}}
\newcommand{\modules}{\mathrm{-mod}}
\newcommand{\Rhom}{\mathrm{RHom}}

\newcommand{\spec}{\mathrm{Spec}}

\newcommand{\A}{\mathcal{A}}
\newcommand{\C}{\mathcal{C}}
\newcommand{\OO}{\mathcal{O}}

\newcommand{\m}{\mathrm{m}}

\newcommand{\G}{\mathcal{G}}
\newcommand{\T}{\mathcal{T}}
\newcommand{\F}{\mathcal{F}}

\newcommand{\D}{\mathcal{D}}
\newcommand{\X}{\mathcal{X}}
\newcommand{\Y}{\mathcal{Y}}
\newcommand{\M}{\mathcal{M}}

\renewcommand{\I}{\mathcal{I}}
\newcommand{\cN}{\mathcal{N}}

\newcommand{\HH}{\mathrm{H}}
\newcommand{\PP}{\mathbb{P}}
\newcommand{\B}{\mathbb{B}}

\newcommand{\CC}{\mathrm{C}}

\newcommand{\End}{\mathrm{End}}

\newcommand{\Vect}{\mathrm{Vect}}

\newcommand{\Hom}{\mathrm{Hom}}

\newcommand{\kk}{\Bbbk}

\newcommand{\RcHom}{\mathrm{R}\mathscr{Hom}\text{om}}

\newcommand{\An}{\A_n}
\newcommand{\Am}{\A_m}
\newcommand{\Ainfty}{\A_\infty}

\newcommand{\nocontentsline}[3]{}
\newcommand{\tocless}[2]{\bgroup\let\addcontentsline=\nocontentsline#1{#2}\egroup}

\newtheorem{theorem}{Theorem}[section]
\newtheorem{proposition}[theorem]{Proposition}
\newtheorem{lemma}[theorem]{Lemma}
\newtheorem{corollary}[theorem]{Corollary}
\newtheorem*{theorem*}{Theorem}

\theoremstyle{definition}
\newtheorem{definition}[theorem]{Definition}
\newtheorem{lemma/definition}[theorem]{Definition/Lemma}
\newtheorem{example}[theorem]{Example}

\newtheorem{remark}[theorem]{Remark}

\newtheorem{construction}[theorem]{Construction}

\begin{document}

\title{Twisted Hodge diamonds give rise to non-Fourier-Mukai functors}
\author{Felix K\"ung}
\begin{abstract}

We apply computations of twisted Hodge diamonds to construct an infinite number of non-Fourier-Mukai functors with well behaved target and source spaces.

To accomplish this we first study the characteristic morphism introduced in \cite{Buchweitz2006} in order to control it for tilting bundles. Then we continue by applying twisted Hodge diamonds of hypersurfaces embedded in projective space to compute the Hochschild dimension of these spaces. This allows us to compute the kernel of the embedding into the projective space in Hochschild cohomology. Finally we use the above computations to apply the construction in \cite{Rizzardo2019} of non-Fourier-Mukai functors and verify that the constructed functors indeed cannot be Fourier-Mukai for odd dimensional quadrics.  

Using this approach we prove that there are a large number of Hochschild cohomology classes that can be used for the construction of \cite{Rizzardo2019}. Furthermore, our results allow the application of computer-based calculations to construct candidate functors for arbitrary degree hypersurfaces in arbitrary high dimensions. Verifying that these are not Fourier-Mukai still requires the existence of a tilting bundle. 

In particular we prove that there is at least one non-Fourier-Mukai functor for every odd dimensional smooth quadric.
\end{abstract} 
\maketitle

\section{introduction}

\subsection{Background and Results}

The concept of Fourier-Mukai functors generalizes the idea of a correspondence to the categorical level.

\begin{definition}
A functor $f: \D^b\left(X\right) \to \D^b\left(Y\right)$ between bounded derived categories of schemes is called Fourier-Mukai if there exists an object $M\in \D^b\left(Y\times X\right)$ such that $f\cong \Phi_M:= \mathrm{R}\pi_{Y,*}\left(M\stackrel{L}{\otimes} \mathrm{L}\pi^*_X \left(\_\right)\right)$. In this case $M$ is called the Fourier-Mukai kernel.
\end{definition}

In particular these functors can be understood geometrically as $\Phi_M$ admits a complete charactersiation by $M \in \D^b\left(Y \times X\right)$. 

It also turns out that most functorial constructions done in algebraic geometry are Fourier-Mukai.  This means that understanding the property of being Fourier-Mukai, respectively of not being Fourier-Mukai, is essential for understanding which functors between derived categories of sheaves may arise from geometric constructions and which do not. Another indicator of the geometric nature of Fourier-Mukai functors are the following results by V. Orlov and B. To\"en:

\begin{theorem*}\cite{Orlov1997}
Let $X$ and $Y$ be smooth projective schemes. Then every fully faithful exact functor $\Psi_M:\D^b\left(X\right)\to \D^b\left(Y\right)$ is a Fourier-Mukai functor for some Fourier-Mukai kernel $M\in \D^b\left(X\times Y\right)$.  
\end{theorem*}
\begin{theorem*}\cite{Toen2004}
Let $X$	and $Y$ be smooth projective schemes. Then a functor $\D^b\left(X\right)\to \D^b\left(Y\right)$ is precisely Fourier-Mukai if it is induced by a dg-functors between the canonical dg-enhancements.
\end{theorem*} 
The above results show that a lot of functors between derived categories of smooth projective schemes are Fourier-Mukai. So Bondal, Larsen and Lunts \cite{Bondal2004} conjectured nearly 20 years ago that every exact functor between such derived categories admits a description as a Fourier-Mukai functor.

This conjecture was disproven fifteen years later when A. Rizzardo, M. Van den Bergh and A. Neeman \cite{Rizzardo2019} constructed the first non-Fourier-Mukai functor
 $$\Psi_\eta:\D^b\left( Q_3\right)\hookrightarrow \D^b\left( \PP^4\right),$$
  where $Q_3$ denotes the smooth three dimensional quadric in $\PP^4$. Shortly thereafter V. Vologodsky constructed in a note \cite{Vologodsky2016} another class of non-Fourier-Mukai functors over a field of characteristic $p>0$. However, Vologodsky's functor turns out to be liftable to a $\Z_p$-linear dg-level, whereas the example from \cite{Rizzardo2019} can be proven to not even have a lift to the spectral level if one works over the rational numbers. 
  
In this work we generalize the result from \cite{Rizzardo2019} to higher dimensions. In particular we will work over a field of characteristic zero in order to show that even in the nicest possible case there is an abundance of non-Fourier-Mukai functors.

We then verify that in the case of a smooth odd dimensional quadric we can apply our result to get a non-Fourier-Mukai functors in arbitrary high dimensions. 

\begin{theorem*}
Let $Q \hookrightarrow \PP^{2k}$ be the embedding of a smooth odd dimensional quadric for $k>2$. Then we have an exact functor
$$\Psi_\eta :\D^b\left(Q\right)\to \D^b\left(\PP^n\right)$$ that cannot be Fourier-Mukai.
\end{theorem*}

\subsection{Proof strategy}

Generally we follow the ideas from \cite{Rizzardo2019}. In order to conclude that we can construct more non-Fourier-Mukai functors we include an auxiliary results on the kernel of the push forward in Hochschild cohomology. Furthermore we will use more general objects, degrees and indices. We need to do this as the proof in \cite{Rizzardo2019} is very specialized to the three dimensional quadric and one needs to take care when generalizing their strategy to a more general setting.

Recall that the construction in \cite{Rizzardo2019} proceeds in two steps:
\begin{enumerate}
\item\label{Introduction enumeration item 1} First the authors construct a prototypical non-Fourier-Mukai functor between not necessarily geometric dg-categories.
\item\label{Introduction enumeration item 2} Using behaviour of Hochschild cohomology under embeddings this functor is turned into a geometric functor.
\end{enumerate}

More precisely, in step (\ref{Introduction enumeration item 1}) \cite{Rizzardo2019} construct a functor 
$$L:\D^b\left( X\right)\to \D_\infty\left(\X_\eta\right)$$
 for a smooth scheme $X$ and $\eta \in \HHH^{\ge \dim X + 3}$, where $\D_\infty\left(\X_\eta\right)$ is the derived category of an $\Ainfty$-category arising as infinitesimal deformation in the $\eta$-direction.

In step (\ref{Introduction enumeration item 2}) the construction of $L$ is turned into a geometric one. In \cite{Rizzardo2019} this is achieved by showing that the canonical $\eta\in \HHH^{2 \dim Q_3}\left(Q_3, \omega^{\otimes 2}_{Q_3}\right)$ is annihilated by the embedding $Q_3 \hookrightarrow \PP^4$, which allows the passing from the algebraic world to the geometric world. The authors then define $\Psi_\eta$ to be $L$ composed with the pushforward into the geometric category $\D^b\left(\PP^4\right)$.

Although the construction in \cite{Rizzardo2019} is very general, it has two major drawbacks:

The first is that although $L$ is constructed to be prototypical non-dg it is not obvious that the composition with the pushforward is again non-Fourier-Mukai. One usually handles these complications by applying an inductive obstruction theory that gets unwieldy quickly as one needs to keep track of inductively chosen lifts. Indeed \cite{Rizzardo2019} only gives a single example of a non-Fourier-Mukai functor although the construction given in step (\ref{Introduction enumeration item 1}) and (\ref{Introduction enumeration item 2}) is very general in nature. 

We are able to solve this issue by restricting to Hochschild cohomology classes in degree $\dim\left(X\right)+3$, this leads to the first obstruction vanishing and so we do not need to control the previous lifts in order to conclude that the pushed forward obstruction does not vanish.

The second drawback is that the results in \cite{Rizzardo2019} rely heavily on the existence of a tilting bundle in order to conclude that the prototypical functor $L$ cannot be dg. Furthermore in \cite{Raedschelders2019} T. Raedschelders, A. Rizzardo and M. Van den Bergh construct an infinite amount of non-Fourier-Mukai functors using the prototypical $L$ mentioned above. However, to do this they apply a geometrification result by Orlov and hence lose control over the target space. In particular the above mentioned geometrification result relies even more on the existence of a tilting bundle. Although our concrete examples still require the existence of a tilting bundle we study the naturality of the characteristic morphism, which might in future allow results using more general generators. In particular we phrase our main result such that a non-vanishing characteristic morphism suffices, which is guaranteed for tilting objects.

Altogether this work improves on the construction from \cite{Rizzardo2019} to prove the existence of non-Fourier-Mukai functors
$$\Psi_\eta: \D^b\left(Q\right)\to \D^b\left(\PP^{n+1}\right),$$
for $Q$ a smooth quadric in arbitrary high dimension.
 
Furthermore one can use our results to calculate the dimensions of choices for constructing candidate non-Fourier-Mukai functors as entries in twisted Hodge diamonds. For instance, if one wants to deform a smooth degree $6$ hypersurface $f:X\hookrightarrow \PP^{n+1}$ along the Hochschild cohomology of $\OO_X\left(-8\right)$ in a way that might gives rise to a non-Fourier-Mukai functor, we may pick an $\eta$ in a $20993$-dimensional space:

$$\begin{matrix}
\;&\;&\;&\;&\;&0&\;&\;&\;&\;&\;\\
\;&\;&\;&\;&0&\;&0&\;&\;&\;&\;\\
\;&\;&\;&0&\;&0&\;&0&\;&\;&\;\\
\;&\;&0&\;&0&\;&0&\;&0&\;&\;\\
\;&0&\;&0&\;&0&\;&0&\;&0&\;\\
\mathllap{299}6&\;&\color{purple}\mathclap{20993}\color{purple}&\;&\mathclap{15267}&\;&\mathclap{917}&\;&0&\;&0\\
\;&\mathllap{157}5&\;&0&\;&0&\;&0&\;&0&\;\\
\;&\;&\mathllap{577}5&\;&0&\;&0&\;&0&\;&\;\\
\;&\;&\;&\mathllap{1039}5&\;&0&\;&0&\;&\;&\;\\
\;&\;&\;&\;&\mathllap{900}2&\;&0&\;&\;&\;&\;\\
\;&\;&\;&\;&\;&\mathllap{299}6.&\;&\;&\;&\;&\;\\
\end{matrix}$$

\section*{Acknowledgements}

I would like to thank my supervisor A. Rizzardo for giving me this awesome topic and supporting me throughout my work on this project.

\section*{Notation}

Throughout this work we consider $\kk$ to be a field of characteristic zero and all schemes, algebras, $\Ainfty$-categories and dg-categories are considered to be over $\kk$. We will assume all $\Ainfty$-structures to be strictly unital and graded cohomologically. 

Furthermore the bounded derived category of coherent sheaves over a scheme $X$ will be denoted by $\D^b\left(X\right)$ or $\D^b\left(\Coh\left(X\right)\right)$ depending on the context, wherever we need to pass to the category $\D^b_{\Coh X}\left(\Qch \left(X\right) \right)$ using \cite[Proposition~3.5]{HuybrechtsFm} we will indicate this. We denote the derived category of modules over a $\kk$-linear category $\X$ as $\D\left(\X\right)$ and also use the same notation for dg-categories. We refer to the  dg-category of $\Ainfty$-modules over an $\Ainfty$-algebra $\X_\eta$ with homotopic maps identified by $\D_\infty\left(\X_\eta\right)$, this is often also referred to as the derived category of $\X_\eta$-modules.

 The change of rings functor associated to a $\kk$-linear functor $f:\X \to \Y$ will be referred to by $f_*:\D\left(\Y\right)\to \D\left(\X\right)$ in order to be compatible with the notation for schemes. Also, wherever applicable, functors are intended as derived.

\section{Preliminaries: $\Ainfty$ Deformations of Schemes and Objects}

For the general notion of $\Ainfty$-structures we refer to \cite{Lefevre}, for an English reference consult \cite{kellerAinfty} or \cite{Efimov2010}.

We recall the following results from \cite{Rizzardo2019} which we will need to construct our candidate functors. In particular we refer the interested reader to \cite{Rizzardo2019} for a more in depth approach.

\begin{definition}
Let $\X$ be a $\kk$-linear category. An $\X$-bimodule $\M$ is called $\kk$-central, if the $\kk$-action induced by the left $\X$-action coincides with the $\kk$-action induced by the right $\X$-action.
\end{definition}

\begin{definition}\label{definition Hochschild cohomology category}
Let $\X$ be a small $\kk$-linear category and let $\M$ be a $\kk$-central $\X$-bimodule. The Hochschild complex $\CC^*\left(\X,\M\right)$ is defined as
$$\CC^n\left(\X,\M\right):=\prod_{\mathclap{X_0,...,X_n\in \obj\left(\X\right)}}\Hom\left(\X\left(X_{0},X_1\right)\otimes_\kk ... \otimes_\kk  \X\left(X_{n-1},X_n\right),\M\left(X_0,X_n\right)\right)$$
with differential given by
\begin{align*}
d f \left(x_1\otimes...\otimes x_{n+1}\right) :=&x_1 f\left(x_{2}\otimes...\otimes x_{n+1}\right)\\&+\sum^{n}_{i=1}\left(-1\right)^{i} f \left(x_1\otimes...\otimes x_{i}x_{i+1}\otimes ...\otimes x_{n+1}\right) \\
&+ \left(-1\right)^{n+1}   f \left(x_1\otimes ... \otimes x_n\right)x_{n+1}.
\end{align*}

The Hochschild cohomology $\HHH^*\left(\X,\M\right) $ is the cohomology of $\CC^*\left(\X,\M\right)$.
\end{definition}

\begin{definition}\label{definition tensor product with k-algebra}
Let $\X$ be an $\Ainfty$-category and let $\Gamma$ be a $\kk$-algebra. Then we define the $\Ainfty$-category $\X\otimes_\kk \Gamma$ to consist of the same objects as $\X$ and morphism spaces given by
\begin{align*}
\left(\X\otimes_\kk \Gamma\right) \left(a,b\right):= \X\left(a,b\right)\otimes_\kk \Gamma,
\end{align*}
with higher composition morphisms given by
\begin{align*}
\m_{i,\X\otimes_\kk \Gamma}\left(\left(x_1\otimes_\kk \gamma_1\right) ,..., \left(x_i\otimes_\kk \gamma_i\right) \right) := \m_{i,\X}\left(x_1,..., x_i\right) \otimes_\kk \gamma_1 ... \gamma_i
\end{align*} 
for composable arrows $x_j\otimes \gamma_j \in \left(\X\otimes_\kk\Gamma \right)\left(a_{j-1},a_{j}\right) $.
\end{definition}
Observe that there is no sign arising as we are considering $\Gamma$ to be a $\kk$-linear algebra and so all $\gamma_i$ are in degree $0$.

We now will define a version of $\X$ deformed along a Hochschild cocycle $\eta$. For a more in depth discussion of this construction we refer to \cite[\S 6]{Rizzardo2019}.

\begin{definition}\label{definition deformed A-infinity category}
Let $\X$ be a small $\kk$-linear category, $\M$ a $\kk$-central $\X$-bimodule and let $\eta\in \CC^{\ge 3}\left(\X,\M\right)$ such that $ d \eta =0$.

We define the $\Ainfty$-category $\X_\eta$ to have the same objects as $\X$, morphism spaces given by $$\X_\eta \left(a,b\right):= \X\left(a,b\right)\oplus \M\left(a,b\right)[n-2]$$
and non-zero composition morphisms
\begin{align*}
\m_2\left(\left(x,m\right),\left(x',m'\right)\right)&:=\left(xx',xm'+mx'\right)\\
\m_n\left(\left(x_1,m_1\right),...,\left(x_n,m_n\right)\right)&:=\left(0,\eta\left(x_1,...,x_n\right)\right)
\end{align*}
for composable arrows $x_1,...,x_n$.

The category $\X_\eta$ comes with a canonical $\kk$-linear functor $\pi:\X_\eta \to \X$ acting by the identity on objects and on morphisms by
\begin{align*}
\pi:\X_\eta \left(a,b\right)=\X\left(a,b\right)\oplus \M\left(a,b\right)[n-2] &\to \X\left(a,b\right)\\
\left(\varphi,\psi\right)&\mapsto \varphi
\end{align*}
\end{definition}

\begin{proposition}\cite[Lemma 6.1.1]{Rizzardo2019}
Let $\X$ be a $\kk$-linear category, let $\M$ be a $\kk$-central $\X$-bimodule and let $\eta, \mu \in \CC^n\left(\X,\M\right)$  with $n\geq 3$ such that $\overline{\eta}=\overline{\mu}\in \HHH^n\left(\X,\M\right)$.
Then we have
$$\X_\eta\cong \X_{\mu}.$$
\end{proposition}

\begin{remark}\cite[\S 6.2]{Rizzardo2019}
Let $\X$ be a $\kk$-linear category, $\M$ a $\kk$-central $\X$-bimodule and $\Gamma$ a $\kk$-algebra.  Then we define the morphism of Hochschild complexes
\begin{align*}
\CC^*\left(\X,\M\right)&\to \CC^*\left(\X\otimes \Gamma,\M\otimes\Gamma\right) \\
\eta &\mapsto \eta \cup \id
\end{align*}
to send a Hochschild cocycle $\eta$ in degree $n$ to the degree $n$ morphism
\begin{align*}
\eta\cup \id:\left(\X\otimes\Gamma\right)^{\otimes n} &\to \left(\M\otimes \Gamma\right)\\
\left(a_1\otimes \gamma_1\right) \otimes ...\otimes \left(a_n \otimes \gamma_n\right) &\mapsto \eta \left(a_1\otimes ...\otimes a_n\right)\otimes \left(\gamma_1...\gamma_n\right). 
\end{align*}
In particular we get a morphism
\begin{align*}
\HHH^*\left(\X ,\M \right)&\to \HHH^*\left(\X\otimes \Gamma , \M\otimes \Gamma\right)\\
\eta &\mapsto \eta \cup \id.
\end{align*}

One can compute that this morphism is compatible with deformations, i.e.
 $$\X_\eta \otimes \Gamma \cong \left(\X\otimes \Gamma\right)_{\eta\cup 1}.$$
\end{remark}

\begin{definition}\cite[\S 6.4]{Rizzardo2019}\label{definition colift}
Let $\X$ be a small $\kk$-linear category, $\M$ a $\kk$-central $\X$-bimodule, $\eta \in \HHH^* \left(\X,\M\right)$ and let $U \in \X\modules$. A colift of $U$ to $\X_\eta$ is a pair $\left(V,\phi\right) $, where $V \in \D_\infty\left(X_\eta\right)$ and $\phi$ is an isomorphism of graded $\HH^*\left(\X_\eta\right)$-modules $$V \xrightarrow{\sim}\Hom_\X\left(\HH^*\left(\X_\eta\right), U\right)$$ 
\end{definition}

Although we will later discuss the geometric characteristic morphism in depth we recall the next Proposition using the algebraic characteristic morphism from \cite{Rizzardo2019} here, as it introduces obstructions against the existence of colifts. Later in \S \ref{Section Nontrivial kernel in Hochschild cohomology give non-Fourier-Mukai functors} we will compare the geometric and algebraic characteristic morphisms.

\begin{proposition}\label{proposition obstruction colift}
Assume $\M$ is an invertible $\kk$-central $\X$-bimodule and $\X_\eta$ is as in Definition~\ref{definition deformed A-infinity category}. Then we have that the object $U\in\D\left(\X\right)$ has a colift if and only if $c_U\left(\eta\right)=0$, where $c_U$ is the (algebraic) characteristic morphism
$$c_U:\HHH^*\left(\X,M\right)\to \Ext^*\left(U, M\otimes U\right),$$
obtained by interpreting $\eta \in \HHH^*\left(\X,M\right)$ as a degree $n$ morphism $\X\to M$ in $\D\left(\X\otimes \X^\op\right)$ and applying $\_\otimes U$.
\begin{proof}
This is a combination of \cite[6.4.1]{Rizzardo2019} and \cite[6.3.1]{Rizzardo2019}.
\end{proof}
\end{proposition}

\section{Equivariant sheaves and the characteristic morphism}\label{Section Equivariant sheaves and the characteristic morphism}

 In this section we define $\Gamma$-equivariant sheaves on a scheme $X$, for a $\kk$-algebra $\Gamma$. We will use this in order to study the (geometric) $\Gamma$-equivariant characteristic morphism.
\subsection{Equivariant Sheaves and Fourier-Mukai Functors}

In this section we introduce equivariant sheaves and prove that the equivariant structure is compatible with Fourier-Mukai functors. In particular we can use this later to get a contradiction to being Fourier-Mukai.

\begin{definition}\cite[\S 4]{Lowen2006}\label{definition Gamma equivariant sheaves}
Let $\Gamma$ be a $\kk$-algebra and $\C$ a $\kk$-linear category. Then we define the category $\C_\Gamma$ to consist of objects
 $$\obj\left(\C_\Gamma\right):=\left(M, \psi: \Gamma \to \End_\C\left(M\right)\right),$$
where $\M \in \C$ and $\psi$ is a morphism of $\kk$-algebras, and morphisms
$$\C_\Gamma \left(\left(M,\psi\right),\left(N,\varphi\right)\right):=\left\{ \alpha \in \C\left(F,G\right)| \alpha \circ \psi\left(\gamma \right) = \varphi\left(\gamma\right) \circ \alpha \in \C\left(M,N\right) \quad\forall \gamma \in \Gamma \right\}.$$
\end{definition}
 
We will mostly denote $\left(T,\varphi\right) $ by $T$ if the action is clear from context to avoid clumsy notation.

\begin{example}\label{Examples equivariant sheaves} We give a few examples to illustrate the above Definition~\ref{definition Gamma equivariant sheaves}:

\begin{itemize}
\item Since $\C$ is required to be $\kk$-linear we have that $\End_\C\left(\_\right)$ comes with a canonical $\kk$-action and so we have $$\C_\kk\cong \C.$$
\item Let $\C$ be a $\kk$-linear category and $M\in \C$, then we have canonically
\begin{equation*}
M=\left(M,\id \right)\in \C_{\End_\C\left(M\right)}.
\end{equation*}
\item Let $F:\C\to \C'$ be a $\kk$-linear functor between $\kk$-linear categories and let $\Gamma$ be a (possibly non-commutative) $\kk$-algebra. Then we can extend $F$ canonically to a functor
\begin{align*}F: \C_\Gamma&\to \C'_\Gamma\\
M&\mapsto FM:=\left(FM,F\circ \gamma\right)\in \C'_\Gamma.
\end{align*}
\item Consider the point $*=\spec\left(\kk\right)$ and a $\kk$-algebra $\Gamma$. Then $\left(M,\varphi\right)\in \Coh\left(*\right)_\Gamma$ consists of $M\in \Coh\left(*\right)\cong \Vect_\kk$ and a $\kk$-algebra morphism $\varphi: \Gamma \to \End_\kk\left(M\right)$. Which means that 
$$\Coh\left(*\right)_\Gamma\cong \Gamma\modules .$$
\item Let $T\in \Coh\left(X\right)$ be tilting for $X$ smooth projective and set $\Gamma:=\End_X\left(T\right)$. Then we have
$$\D^b\left( X\right)\cong \D^b\left(\Gamma\right) \cong \D^b\left(\Coh\left(*\right)_\Gamma\right).$$
We will prove in Lemma~\ref{Tilting gives HT equivalence} that this equivalence is compatible with products of schemes under mild conditions.

We will use the next specific version of the second example throughout this work:
\item Let $\F\in \Coh\left(X\right)$ be a coherent sheaf on a scheme. Then we have canonically
\begin{equation}\label{Example coherent sheave with action by endomorphisms}
F=\left(\F,\id\right)\in \Coh\left(X\right)_{\End_X\left(\F\right)}.
\end{equation}
\end{itemize}
\end{example}

\begin{remark}\label{Remark difference homotopy coherent action or not}
The categories $\D\left(\C_\Gamma\right)$ and $\D\left(\C\right)_\Gamma$ may seem very similar in notion, however, they do not coincide. An object in $M\in \D\left(\C_\Gamma\right)$ can be interpreted as a complex of equivariant objects, i.e. it admits an action in every degree and a differential that is compatible with these actions. On the other hand an object in $\D\left(\C\right)_\Gamma$ can be interpreted as a complex of sheaves together with an action on the whole complex that suffices the relations given by the $\Gamma$-action up to homotopy. The difference between these two notions essentially boils down to the difference between commutative diagrams up to homotopy not coinciding with homotopy commutative diagrams, which also led to the development of derivators \cite{Groth2011}. For some more information on this interplay we refer to \cite{Rizzardo2014}.

By the above discussions there is a canonical forgetful functor
\begin{align*}
\pi: D\left(\C_\Gamma\right) &\to \D\left(\C\right)_\Gamma\\
\overline{\left(M^\cdot,\varphi^\cdot\right)} &\mapsto \left(\overline{M}^\cdot,\overline{\varphi}^\cdot\right),
\end{align*}
where we denote by $\overline{\left(\_\right)}$ an equivalence class of $\left(\_\right)$.
One can think of the above functor as forgetting that $\Gamma$ acts on every degree separately. However, this functor is neither essentially injective nor surjective in general, which we will use later.
\end{remark}

\begin{remark}[{\cite[Remark~2.51]{HuybrechtsFm}}]
Let $f:\C \to \C'$ be a left or right exact functor betweed abelian categories.
Recall that an object $M \in \A$ is called $f$-adapted if $\mathrm{R}^i f \left( M \right)\cong 0$ respectively $\mathrm{L}^i f\left(M\right)\cong 0$ for $i>0$.
\end{remark}

\begin{lemma}\label{Deriving is compatible with Gamma action}
Let $f:\C \to \C'$ be a right or left exact functor between abelian $\kk$-linear categories such that $\C$ has enough $f$-adapted objects and let $\Gamma$ be a $\kk$-algebra. Then the canonical functor
\begin{align*}
f:\D^\natural\left(\C_\Gamma \right)&\to \D^\natural\left(\C'\right)_\Gamma\\
\overline{\left(M,\psi\right)} &\mapsto \left(\overline{M},\overline{f \circ \psi}\right)
\end{align*}
admits a lift $$ f_{\Gamma}: \D^\natural\left(\C_\Gamma\right)\to \D^\natural\left(\C'_\Gamma\right),$$

with $\natural \in \left\{b,+,-,\;\right\}$. In the case $\natural = b$, respectively $\natural=-$ for left exact and $\natural=+$ for right exact functors, we assume that every $M\in \C$ admits a bounded $f$-adapted resolution.
\begin{proof}
We have by Example~\ref{Examples equivariant sheaves} a canonical functor
\begin{align*}
f_\Gamma: \C_\Gamma &\to \C'_\Gamma\\
\left(M,\psi\right)&\mapsto \left(fM,f\circ \psi\right).
\end{align*}

Now as $\C_\Gamma$ and $\C'_\Gamma$ are abelian with kernels and cokernels computed on objects we get that $f_\Gamma$ has the same exactness as $f$, and as cohomology also is computed on $M$ only we get that every $f$-adapted object is also $f_\Gamma$ adapted.

Since we have enough $f$-adapted objects we may consider for $M\in \D^\natural\left(\C\right)$ an $f$-adapted replacement, which by assumption is also finite for $\natural=b$ respectively if $f$ is left exact and $\natural=-$ or $f$ being right exact and $\natural=+$. In particular we may invoke \cite[Theorem~10.5.9]{Weibel1994} in order to find a well-defined derived functor:
 \begin{align*}
f_{\Gamma}:\D^\natural\left(\C_\Gamma\right)&\to \D^\natural\left(\C'_\Gamma\right)\\
\overline{\left(M,\psi\right)^\cdot}&\mapsto \overline{\left(f M,f\circ \psi\right)^\cdot}.
\end{align*}
Furthermore \cite[Theorem~10.5.9]{Weibel1994} allows us to freely use $f$-adapted resolutions to compute $f_\Gamma$ on the derived category, i.e. we will assume from now on that every $M^\cdot$ is $f$-adapted.

Recall the functor from Remark~\ref{Remark difference homotopy coherent action or not}
\begin{align*}
\pi: \D^\natural\left(\C'_\Gamma\right) &\to \D^\natural\left(\C'\right)_\Gamma\\
\overline{\left(F,\psi\right)^\cdot} &\mapsto \left(\overline{F^\cdot},\overline{\psi^\cdot}\right).
\end{align*}
 We now just need to verify that the diagram
$$\tikz[heighttwo,xscale=4,yscale=2,baseline]{
\node (A) at (0,1) {$\D^\natural\left(\C_\Gamma\right)$};
\node (B) at (1,1) {$\D^\natural\left(\C'_\Gamma\right)$};
\node (D) at (1,0) {$\D^\natural\left(\C'\right)_\Gamma$};
\draw[->]
(A) edge node[below left] {$f$} (D)
(B) edge node[left]{$\pi$} (D)
(A) edge node[above]{$f_{\Gamma}$}(B);
}
$$
commutes. 

We indeed get
\begin{align*}
\pi\circ f_{\Gamma} \overline{\left(M, f \circ \psi\right)^\cdot}&=\pi\overline{\left(f M,f\circ \psi\right)^\cdot}\\
&=\left(\overline{f M^\cdot},\overline{f \circ \psi^\cdot}\right)\\
&=f\left(\overline{M^\cdot},\overline{\psi}^\cdot\right)
\end{align*}
as claimed.
\end{proof}
\end{lemma}

We will drop the $\Gamma$ in $f_\Gamma$ if it is clear from context, respectively from the target or source categories.

\begin{lemma}\label{direct image, pullback and tensor product compatible with basechange}
Let $f:X \to Y$ be a morphism of finite-dimensional noetherian $\kk$-schemes, $\Gamma$ a $\kk$-algebra and $M\in \D^b\left(X\right)$. Then we have the following:
\begin{itemize}
\item If $f$ is proper, then the functor $f_*:\D^b\left(  \Coh\left(X\right)_\Gamma\right)\to \D^b\left(\Coh\left(Y\right)\right)_\Gamma$ admits a canonical lift
$$f_{*,\Gamma}:\D^b\left(\Coh\left(X\right)_\Gamma\right)\to \D^b\left(\Coh\left(Y\right)_\Gamma\right). $$ 
\item If $f$ is flat, then the functor $f^*:\D^b\left( \Coh\left(Y\right)_\Gamma\right)\to \D^b\left( \Coh\left(X\right)\right)_\Gamma$ admits a canonical lift
$$f_{\Gamma}^*:\D^b\left(\Coh\left(Y\right)_\Gamma\right)\to \D^b\left(\Coh\left(X\right)_\Gamma\right). $$ 
\item If $X$ is regular, then the functor $M\otimes \_:\D^b\left(\Coh\left(X\right)_\Gamma\right)\to \D^b\left(\Coh\left(X\right)\right)_\Gamma$ admits a canonical lift
$$M\otimes_\Gamma\_:\D^b\left(\Coh\left(X\right)_\Gamma\right)\to \D^b\left(\Coh\left(X\right)_\Gamma\right). $$ 
\end{itemize}
\begin{proof}
We check the cases separately:
\begin{itemize}
\item[$f_*$:] By Lemma~\ref{Deriving is compatible with Gamma action} it suffices to show that every coherent sheaf $M$ admits an $f_*$-adapted finite resolution in $\Coh\left(X\right)$.
 By \cite[Theorem 3.22]{HuybrechtsFm} the object $M$ admits an $f_*$-adapted resolution of finite length of quasi-coherent sheaves.
By \cite[Theorem 3.23]{HuybrechtsFm} these quasi-coherent sheaves can be picked to be coherent for $f$ proper.

So we can find by Lemma~\ref{Deriving is compatible with Gamma action} a lift
  $$f_{ * ,\Gamma}: \D^b\left(\Coh\left(X\right)_\Gamma\right) \to \D^b\left(\Coh\left(Y\right)_\Gamma\right).$$
\item[$f^*$:] As $f$ is flat $f^*$ is exact and does not need to be derived. In particular we get by Lemma~\ref{Deriving is compatible with Gamma action} immediately a lift 
$$f_{\Gamma}^*:\D^b\left(\Coh\left(Y\right)_\Gamma\right)\to \D^b\left(\Coh\left(X\right)_\Gamma\right). $$ 
\item[$M\otimes\left(\_\right)$:] By \cite[Proposition 3.26]{HuybrechtsFm} we have that every $\F\in \Coh\left(X\right)$ admits a bounded locally free resolution, which is in particular $M\otimes \left(\_\right) $ adapted. So we get by Lemma~\ref{Deriving is compatible with Gamma action} that $M\otimes\left(\_\right)$ admits a lift
$$M\otimes_\Gamma\left(\_\right):\D^b\left(\Coh\left(X\right)_\Gamma\right)\to \D^b\left(\Coh\left(X\right)_\Gamma\right) $$
\end{itemize}
as claimed.
\end{proof}
\end{lemma}

\begin{corollary}\label{Fourier-Mukai compatible with basechange}
Let $f:\D^b\left(X\right)\to \D^b \left( Y\right)$ be a Fourier-Mukai functor between finite-dimensional smooth projective $\kk$-schemes and let $\Gamma$ be a $\kk$-algebra. Then we have that the induced functor
$$f:\D^b \left( \Coh\left(X\right)_\Gamma\right) \to \D^b \left( \Coh\left(Y\right)\right)_\Gamma$$
admits a lift:
$$f_\Gamma:\D^b \left(\Coh\left(X\right)_\Gamma\right) \to \D^b \left(\Coh\left(Y\right)_\Gamma\right).$$
\begin{proof}
Observe first that $X$ and $Y$ being smooth projective immediately gives that
\begin{align*}
\pi_1:X \times Y \to X & \text{ is proper,}\\
\pi_2:X \times Y \to Y & \text{ is flat}\\
\text{and } X\times Y & \text{ is regular.}
\end{align*}
As $f$ is a Fourier-Mukai functor it has the form $\pi_{1,*}\left(M\otimes\pi_2^*\left(\_\right) \right)$ for some $M\in \D^b\left(Y\times X\right)$. So we get by Lemma~\ref{direct image, pullback and tensor product compatible with basechange} that  $\pi_{1,*}$, $\pi_2^*$ and $M\otimes \_$ admit canonical lifts. In particular $f$ admits a canonical lift $$f_\Gamma:=\pi_{1,*,\Gamma}\left(M\otimes_\Gamma \pi_{2,\Gamma}^*\left(\_\right)\right)$$
as claimed.
\end{proof}
\end{corollary}

\subsection{Hochschild Cohomology and the characteristic Morphism}

As we want to study the characteristic morphism we start by recalling the definition of the (geometric) Hochschild cohomology:

\begin{definition}\label{definition Hochschild cohomologoy scheme}\cite{Swan1996}
Let $X$ be a separated scheme and $M$ a sheaf on $X$. Then the Hochschild cohomology of $X$ with coefficients in $M$ is given by
$$\HHH^*\left(X,M\right):=\Ext^*_{X\times X}\left(\OO_\Delta, \Delta_*M\right),$$
where $\Delta: X \hookrightarrow X\times X$ is the diagonal embedding.
\end{definition}

For the definition of the (geometric) characteristic morphism below we follow \cite{Lowen2007} and \cite[\S 3.3]{Buchweitz2006}.

\begin{definition}\label{definition characteristic morphism}
Let $X,Y$ be regular schemes, $\Gamma$ a $\kk$-algebra and $M,T\in \Coh\left(X\right)$. Then the (geometric) characteristic morphism is defined to be
\begin{align*}
 c_T \left(M\right): \HHH^*\left(X,M\right)=\Ext_{X\times X}^* \left(\OO_\Delta, \Delta_* M\right) &\to \Ext^*_X\left(T,M\otimes T\right) \\
 \left(\alpha: \OO_{\Delta} \to \Sigma^n \Delta_* M\right) &\mapsto \left(T \xrightarrow{\pi_{1*}\left(\alpha\otimes \pi^*_{2}\id\right)} \Sigma^n M\otimes T	\right),
 \end{align*}
 where we use $T\cong \pi_{1*} \left(\OO_\Delta \otimes \pi^*_{2}T\right)$ and $M\otimes T \cong  \pi_{1*} \left(\Delta_*\Sigma^n M \otimes \pi^*_{2}T\right)$.
 If we have a $\Gamma$-action on $T$, i.e. $\left(T,\varphi\right)\in \Coh\left(X\right)_\Gamma$, there also exists a $\Gamma$-equivariant characteristic morphism 
\begin{align*}
 c_{T,\Gamma} \left(M\right): \HHH^*\left(X,M\right)=\Ext_{X\times X}^* \left(\OO_\Delta, \Delta_* M\right) &\to \Ext^*_{\Coh\left(X\right)_\Gamma}\left(T,M\otimes T\right) \\
 \alpha:  \left(\OO_{\Delta} \to \Sigma^n \Delta_* M\right) &\mapsto \left(T \xrightarrow{\pi_{1*}\left(\alpha\otimes \pi^*_{2}\id\right)} \Sigma^n M \otimes T\right),
\end{align*}
where we consider $M\otimes T$ as on object in $\Coh\left(X\right)_\Gamma$ via the functor $M\otimes \left(\_\right)$, i.e.
\begin{align*}
\psi: \Gamma &\to \End\left(M\otimes T\right)\\
\gamma &\mapsto \id \otimes \varphi\left( \gamma\right).
\end{align*}
\end{definition}

To study the characteristic morphism for special $T$ we will define the following functor realizing the characteristic morphism on a categorical level.

\begin{definition}\label{definition CT}
Let $X,Y$ be projective schemes and let $T=\left(T,\varphi\right)\in \Coh\left(Y\right)_\Gamma$. Then we define the functor:
\begin{align*}
C^X_T: D^b\left(X\times Y\right)&\to \D^b\left(\Coh\left(X\right)_\Gamma\right)\\
M &\mapsto \left(\pi_{1*} \left(M\otimes \pi_{2}^*T \right),\gamma \mapsto   \pi_{1*} \left(\id\otimes \pi_{2}^*\varphi\left(\gamma\right) \right)\right)\\
\left(\alpha: M\to N\right) &\mapsto C_T^X\left(\alpha\right)= \pi_{1*} \left(\alpha\otimes \pi_{2}^*T \right).
\end{align*}
\end{definition}

\begin{remark}\label{CT can be interpreted via FM functors}
One can think of the functor $C_T^X$ to send an object $M\in \D^b\left(X\times Y\right)$ to the image of $T$ under the Fourier-Mukai functor with kernel $M$, equipped with the action induced by $\Phi_{M,\Gamma}$, i.e.
\begin{align*}
C_T^X :D^b\left(X\times Y\right) &\to \D^b\left(\Coh\left(X\right)_\Gamma\right)\\
M&\mapsto \Phi_M\left(T\right).
\end{align*}
\end{remark}

The functor $C_T^X$ allows us to compute $c_{T,\Gamma}$ on a categorical level.

\begin{proposition}\label{cT=CT}
Let $X$ be a scheme and $T\in \Coh\left(X\right)_\Gamma$ and consider $$C_T^X: \D^b \left( X\times X\right)\to \D^b\left(\Coh\left(X\right)_\Gamma\right).$$
Then we have that the equivariant characteristic morphism $c_{T,\Gamma}\left(M\right)$ is given by evaluating the functor $C_T^X$ on the morphism space $\Ext_{X\times X} \left(\OO_\Delta, \Delta_*M\right)$,
$$c_{T,\Gamma} \left(M\right)=C^X_T : \Ext_{X\times X}^*\left(\OO_\Delta ,\Delta_* M\right)\to \Ext_{\Coh\left(X\right)_\Gamma}^*\left(T,M\otimes T\right).$$
\begin{proof}
By Definition~\ref{definition CT} we have
\begin{align*}
C^X_T: \Ext_{X\times X}^*\left(\OO_\Delta ,\Delta_* M\right) &\to \Ext_{\Coh\left(X\right)_\Gamma}\left(C^X_T\left(\OO_\Delta\right),C^X_T \left(\Delta_* M\right)\right)\\
\alpha &\mapsto \pi_{1*}\left(\alpha \otimes \pi_2^* T\right).
\end{align*}
We now have $C_T^X\left(\OO_\Delta\right)\cong T$ and $C_T^X\left(\Delta_*M\right)\cong M \times T$. So the above turns by Definition~\ref{definition characteristic morphism} into $$c_{T,\Gamma}\left(M\right): \Ext_{X\times X}^*\left(\OO_\Delta ,\Delta_* M\right)\to \Ext_{\Coh\left(X\right)_\Gamma}^*\left(T,M\otimes T\right)$$
as claimed.
\end{proof}
\end{proposition}

\begin{definition}[{\cite[§2.1]{Bondal2002}}]
Let $\C$ be a pointed category, i.e. a category admitting a zero object. An object $G\in \C$ is a generator if $\C\left(G,M\right)=0$ implies $M=0$.

Let $\T$ be a pointed graded category. An object $G \in \T$ is called a generator if $\T\left(G,M[i]\right)=0$ for all $i \in \Z$ implies $M=0$.
\end{definition}

\begin{remark}
In a pointed category $\C$ the equation $\C\left(M,N\right)=0$ for two objects $M,N\in \C$ means that the only morphism between $M$ and $N$ is the unique morphism factoring over $0$.

Furthermore we have for an object $M\in \C$ that if $\C\left(M,M\right)=0$ means that $M\cong0$ as in that case $0=\id$ and so the unique morphisms $0 \to M$ and $M\to 0$ define isomorphisms.
\end{remark}

\begin{proposition}\label{proposition pullback of generator is generator}
Let $f_*:\C \to \D$ be a faithful functor between pointed categories or graded functor between pointed graded categories with a left adjoint $f^*:\D\to \C$ and let $T \in \D$ be a generator. Then $f^* T$ is a generator.
\begin{proof}
We cover the case of a graded functor. Observe that the same argument holds for pointed categories by ignoring $[i]$.

Let $ M$ be such that $\C\left(f^*T[i],M\right)=0$. Then we have $\C\left(f^*T[i],M\right)=\D\left(T[i],f_*M\right)=0$, in particular $f_*M\cong 0$. Now $\D\left(f_*M,f_*M\right)=0$ and so $\C\left(M,M\right)=0$. This can only hold if $M\cong 0$ and so $f^*T$ is a generator.
\end{proof}
\end{proposition}

\begin{proposition}\label{G tensor Gamma generator}
Let $\C$ be a $\kk$-linear category or pointed graded category that admits a generator $G$ and let $\Gamma$ be a $\kk$-algebra. Then
\begin{align*}
\left(G\otimes \Gamma ,\psi:\gamma' \mapsto \left(G \otimes \Gamma \xrightarrow{g\otimes \gamma\mapsto g\otimes \gamma'\gamma} G\otimes \Gamma\right) \right)
\end{align*}
defines a generator of $\C_\Gamma$. Where we denote by $G \otimes \Gamma$ the sheaf arising by tensoring locally with the $\kk$-algebra $\Gamma$ as $\kk$-vectorspaces and acting exclusively on $\Gamma$.
\begin{proof}
Let $\left(X,\varphi\right) \in C_\Gamma$ and let $f: G \to X$ be a morphism. Then we have the following morphism in $\C_\Gamma$,
\begin{align*}
\widehat{f}: G\otimes \Gamma &\to X\\
g\otimes \gamma &\mapsto \varphi\left(\gamma\right) \circ f\left(g\right).
\end{align*}
This indeed defines a morphism in $\C_\Gamma$ as
\begin{align*}
\varphi\left(\gamma' \right)\circ \widehat{f}\left(g\otimes \gamma\right)&= \varphi\left(\gamma'\right) \circ \varphi\left(\gamma\right)\circ f\left(g\right)\\
&= \varphi\left(\gamma'\gamma\right)\circ f\left(g\right)\\
&= \widehat{f}\left(g\otimes \gamma'\gamma\right)\\
&= \widehat{f}\circ \psi\left(\gamma'\right)\left(g\otimes \gamma\right).
\end{align*}

We can compute that if $\widehat{f}$ vanishes then $f$ has to vanish as well since
\begin{align*}
0&=\widehat{f}\left(g\otimes \id\right)\\
&=\varphi\left(\id\right)\circ f\left(g\right)\\
&=\id \circ f\left(g\right)\\
&=f\left(g\right).
\end{align*}

This means that if the morphism space $\C_\Gamma\left(G\otimes \Gamma,\left(X,\varphi\right)\right)$ vanishes, that also $\C\left(G,X\right)$ vanishes. Observe that the discussion so far did not assume $G$ to be a generator.

We continue again by considering the graded pointed case. For the pointed case it again suffices to ignore the shift $[i]$.

Now assume that $\C_\Gamma\left(G\otimes \Gamma[i], \left(X,\varphi\right) \right)=0$. Then we have by the above discussion that $\C\left(G[i],X\right)=0$. As $G$ is a generator we get that $X$ has to be a zero object. And so $\left(X,\varphi\right)$ has to be a zero object as well. In particular we get that $G\otimes \Gamma$ is indeed a generator of $\C_\Gamma$.
\end{proof}
\end{proposition}

\begin{remark}
Proposition~\ref{G tensor Gamma generator} is a consequence of $M\otimes \Gamma$ being the free object in $\C_\Gamma$ over $M$.
\end{remark}

\begin{remark}
Recall that an object $\T$ in an abelian category $\A$ is called tilting if $\T$ is a generator in $\D\left(\A\right)$ and $\Ext^i\left(T,T\right)\cong 0$ for all $i>0$.
\end{remark}

\begin{lemma}\label{Tilting gives HT equivalence}
Let $X,Y$ be smooth projective schemes, such that $X$ admits a generator $G\in \D^b\left(X\right)$ with $\Rhom^i\left(G,G\right)$ finite-dimensional for all $i$, let $Y$ be such that it admits a tilting object $T \in \Coh\left(Y\right)$ and set $\Gamma:=\End\left(T\right)$. Then 
$$C^X_T: \D^b\left(X\times Y\right) \to \D^b\left(\Coh\left(X\right)_\Gamma\right)$$
is an equivalence of derived categories.
\begin{proof}
Throughout this proof we denote by $T^D:= \RcHom_Y\left(T, \OO_Y\right)$, the dual of $T$ and by $\pi_{S}:S \to \spec \left(\kk\right)$ the unique projection from a scheme $S$ to the point $\spec \left(\kk\right)$. Observe that by \cite[Proposition 3.26]{HuybrechtsFm} we have  $T^D\in \D^b\left(Y\right)$ as smooth schemes are in particular regular. Furthermore, we will use the following diagram for flat base change twice
$$\tikz[heighttwo,xscale=2,yscale=2,baseline]{
\node (A) at (0,1) {$X$};
\node (B) at (1,2) {$X\times Y$};
\node (C) at (2,1) {$Y$};
\node (D) at (1,0) {$\spec\left(\kk\right)$.};
\draw[->]
(A) edge node[below left]{$\pi_X$}(D)
(B) edge node[above left]{$\pi_1$} (A)
(C) edge node[below right]{$\pi_Y$}(D)
(B) edge node[above right] {$\pi_2$}(C);
}$$

Since $T$ is tilting, it is a generator of $\D^b\left( Y\right)$ and $T^D$ is generating $\D^b\left( Y\right)$ by \cite[Lemma~8.9.1]{Rizzardo2019}. So we get that $G\boxtimes T^D$ generates $\D^b\left(X\times Y\right)$ by \cite[Lemma~3.4.1]{Bondal2002}. Furthermore, we have by \cite[1.10]{Buchweitz2012} that $\Gamma=\End_Y\left(T\right)$ is finite-dimensional.

We first show that $C_T^X\left(G \boxtimes T^D \right)$ is isomorphic to $G \otimes \Gamma$:
\begin{align*}
C^X_T\left(G\boxtimes T^D \right)&= \pi_{1*}\left(\left(G\boxtimes T^D\right)\otimes \pi_2^*T\right)&\text{definition of } C^X_T\\
&\cong \pi_{1*}\left(\pi_1^*G\otimes \pi_2^* T^D\otimes \pi_2^*T \right)&\text{definition of } \boxtimes\\
&\cong \pi_{1*}\left(\pi_1^*G\otimes \pi_2^*\left(T\otimes T^D\right)\right)&\text{\cite[(3.12)]{HuybrechtsFm}}\\
&\cong \pi_{1*}\left(\pi_1^*G \otimes \pi_2^*\RcHom_Y\left(T,T\right) \right)& \text{definition of } T^D\\
&\cong  G\otimes \pi_{1*}\pi_2^*\left( \RcHom_Y\left(T,T\right)\right)&  \text{\cite[(3.11)]{HuybrechtsFm}}\\
&\cong G\otimes \pi_X^* \pi_{Y,*}  \RcHom_Y\left(T,T\right) &\text{flat base change}\\
&\cong G\otimes \Gamma. & T \text{ has no higher $\Ext$-groups} 
\end{align*}
The above computation is compatible with the $\Gamma$-action as all isomorphisms involved are natural isomorphism. In particular replacing $\pi_2^* T$ by $\pi_2^*\gamma$ yields multiplication with $\gamma $ in $\Gamma$.
 
As by Proposition~\ref{G tensor Gamma generator} $G\otimes \Gamma$ is a generator for $\D^b\left(\Coh\left(X\right)_\Gamma\right)$ the functor $C^X_T$ sends a generator to a generator. So it suffices to prove that 
$$\Rhom_{X\times Y}\left(G\boxtimes T^D, G \boxtimes T^D\right)\xrightarrow{C_X^T} \Rhom_{\Coh\left(X\right)_\Gamma}\left(C^X_T\left(G\boxtimes T^D\right),C^X_T\left(G\boxtimes T^D \right)\right) $$
is an isomorphism.

To do that we first compute the source and target spaces:
\begin{align*}
&\Rhom_{X\times Y}\left(G\boxtimes T^D,G\boxtimes T^D\right)\cong\\
&\cong \Rhom_{X\times Y}\left(\pi_1^*G\otimes \pi_2^*T^D,\pi_1^* G\otimes \pi_2^*T^D \right)&\text{definition of } \boxtimes\\
&\cong \Rhom_{X\times Y}\left(\pi_1^*G,\RcHom_{X\times Y}\left(\pi_2^*T^D,\pi_1^*G\otimes \pi_2^*T^D \right)\right)&\text{\cite[(3.14)]{HuybrechtsFm}}\\
&\cong \Rhom_{X\times Y}\left(\pi_1^*G,\pi_1^*G\otimes \RcHom_{X\times Y}\left(\pi_2^*T^D, \pi_2^*T^D \right)\right)&\text{(\cite[3.13)]{HuybrechtsFm}}\\
&\cong \Gamma_{X\times Y} \RcHom_{X\times Y}\left(\pi_1^*G,\pi_1^*G\otimes \RcHom_{X\times Y}\left(\pi_2^*T^D, \pi_2^*T^D \right)\right)&\text{\cite[p.85]{HuybrechtsFm}}\\
&\cong \pi_{\kk,*} \RcHom_{X\times Y}\left(\pi_1^*G,\pi_1^*G\otimes \RcHom_{X\times Y}\left(\pi_2^*T^D, \pi_2^*T^D \right)\right)&\Gamma_{X\times Y} \cong \pi_{X \times Y,*}\\ 
&\cong \pi_{X\times Y,*} \left(\RcHom_{X\times Y}\left(\pi_1^*G,\pi_1^*G\right)\otimes \RcHom_{X\times Y}\left(\pi_2^*T^D, \pi_2^*T^D \right)\right)&\text{\cite[(3.13)]{HuybrechtsFm}}\\
&\cong \pi_{X\times Y,*} \left(\pi_1^*\RcHom_{X}\left(G,G\right)\otimes \pi^*_2\RcHom_{Y}\left(T^D, T^D \right)\right)&\text{\cite[{(3.13)}]{HuybrechtsFm}}\\
&\cong \pi_{X,*}\circ\pi_{1,*} \left(\pi_1^*\RcHom_{X}\left(G,G\right)\otimes \pi^*_2\RcHom_{Y}\left(T^D, T^D \right)\right)&\pi_{X\times Y}=\pi_{X}\circ \pi_{1}\\
&\cong \pi_{X,*}\circ\left(\pi_{1,*} \left(\pi_1^*\RcHom_{X}\left(G,G\right)\otimes \pi^*_2\RcHom_{Y}\left(T^D, T^D \right)\right)\right)&\circ \text{ is associative} \\
&\cong \pi_{X,*} \left(\RcHom_{X}\left(G,G\right)\otimes \pi_{1,*}\pi^*_2\RcHom_{Y}\left(T^D, T^D \right)\right)&\text{\cite[(3.11)]{HuybrechtsFm}}\\
&\cong \pi_{X,*}\left( \RcHom_{X}\left(G,G\right)\otimes \pi_{X}^*\pi_{Y,*}\RcHom_{Y}\left(T^D, T^D \right)\right)&\text{flat base change} \\
&\cong \pi_{X,*}\left(\RcHom_{X}\left(G,G\right)\otimes \pi_{X}^*\Gamma^{\op} \right)&\text{\cite[p.85]{HuybrechtsFm}}\\
&\cong \pi_{X,*}\RcHom_{X}\left(G,G\right)\otimes \Gamma^{\op} &\text{\cite[(3.11)]{HuybrechtsFm}} \\
&\cong \Rhom_{X}\left(G,G\right)\otimes \Gamma^{\op}. &\text{\cite[p.85]{HuybrechtsFm}}
\end{align*}

Now for $\Rhom_{\Coh\left(X\right)_\Gamma}\left(G\otimes \Gamma, G\otimes \Gamma\right)$ we have
\begin{align*}
\Rhom_{\Coh\left(X\right)_\Gamma}\left(G\otimes \Gamma, G\otimes \Gamma\right)&=\Rhom_{\Coh\left(X\right)}\left(G,G\right)\otimes \Rhom_{\Gamma\modules}\left(\Gamma,\Gamma\right) \\
&\cong \Rhom_X\left(G,G\right)\otimes \Gamma^{\textrm{op}}.
\end{align*}
As the two spaces are isomorphic and in particular degree-wise isomorphic, it suffices to prove bijectivity on $\Rhom^i_{X\times Y}\left(G\boxtimes T^D,G\boxtimes T^D\right)$. Since 
$$\Rhom^i_{X\times Y}\left(G\boxtimes T^D,G\boxtimes T^D\right)\cong \Rhom^i_{X}\left(G,G\right)\otimes \Gamma$$
we know that $\Rhom^i_{X\times Y}\left(G\boxtimes T^D,G\boxtimes T^D\right)$ is finite-dimensional as tensor product of finite-dimensional vector spaces. So it suffices to check that $C^X_T$ is surjective. For this let $$\alpha \otimes \beta \in \Rhom_{\Coh\left(X\right)_\Gamma}^i \left(G\otimes \Gamma, G\otimes \Gamma\right)\cong \Rhom^i\left(G,G\right) \otimes \Gamma^\op.$$ Then we can pick $\alpha\boxtimes \beta \in \Rhom^i_{X\times Y}\left(G\boxtimes T^D,G\boxtimes T^D\right)$ and get
\begin{align*}
C^X_T\left(\alpha\boxtimes\beta\right)& \cong \pi_{1,*}\left(\alpha \boxtimes \beta \otimes \pi_2^* \id_T\right)\\
&\cong \pi_{1,*}\left(\alpha\boxtimes \beta \right)\\
&\cong \alpha \otimes \beta \in  \Rhom^i_{\Coh\left(X\right)_\Gamma}\left(G\otimes \Gamma, G\otimes \Gamma\right).
\end{align*}
This means that $C^X_T$ is surjective on the generating set of morphisms of the form $\alpha \otimes \beta$. In particular $C^X_T$ is surjective and an isomorphism as it is surjective between vector spaces of the same dimension which finishes the proof.
\end{proof}
\end{lemma}

\begin{lemma}\label{naturality C}
Let $f:X \to Y$ be a proper morphism of schemes, $\Gamma$ a $\kk$-algebra and $T\in \Coh\left(Y\right)_\Gamma$. We have $f^*T\in \Coh\left(X\right)_\Gamma$.
Consider the two functors:
\begin{align*}
C^X_{f^*T}: \D^b\left(X\times X\right)&\to \D^b\left(\Coh\left(X\right)_\Gamma\right)\\
M &\mapsto \left(\pi_{1*} \left(M\otimes \pi_{2}^*f^*T \right)\right)\\
\left(\alpha: M\to N \right)&\mapsto \pi_{1*} \left(\alpha\otimes \pi_{2}^*f^*T \right)
\end{align*}
and
\begin{align*}
C^X_T\circ \left(\id\times f\right)_*: \D^b\left(X\times X\right)\to \D^b\left(X\times Y\right) \to& \D^b\left(\Coh\left(X\right)_\Gamma\right)\\
M \mapsto \left(\id \times f\right)_* M \mapsto& \left(\pi_{1*} \left(\id\times f\right)_*M\otimes \pi_{2}^*T \right)\\
\left(\alpha: M\to N \right) \mapsto \left(\id \times f\right)_* \alpha \mapsto& \pi_{1*} \left(\left(\id\times f\right)_*\alpha\otimes \pi_{2}^*T \right).
\end{align*}
Then we have a natural isomorphism $C^X_{f^*T}\cong C_T^X \circ \left(\id\times f\right)_*$.
\begin{proof}
Observe that $\id \times f$ is proper as product of proper morphisms and so by \cite[Theorem 3.23]{HuybrechtsFm} 
$$\left(\id\times f\right)_*:\D^b\left(X\times X\right)\to \D^b\left(X\times Y\right)$$
is well-defined.
We will use the following two commutative diagrams in order to construct the isomorphism
$$\tikz[heighttwo,xscale=3.3,yscale=1.7,baseline]{
\node (A) at (0,1) {$X\times X$};
\node (B) at (1,1) {$X\times Y$};
\node (C) at (0,0) {$X$};
\node (D) at (1,0) {$X$};
\draw[->]
(A) edge node[above]{$\id\times f$}(B)
(A) edge node[left]{$\pi_1$} (C)
(C) edge node[below]{$\id$}(D)
(B) edge node[right] {$\pi_1'$}(D);
}
\qquad \; \qquad
\tikz[heighttwo,xscale=3,yscale=1.7,baseline]{
\node (A) at (0,1) {$X\times X$};
\node (B) at (1,1) {$X\times Y$};
\node (C) at (0,0) {$X$};
\node (D) at (1,0) {$Y$,};
\draw[->]
(A) edge node[above]{$\id \times f$}(B)
(A) edge node[left]{$\pi_2$} (C)
(C) edge node[below]{$f$}(D)
(B) edge node[right] {$\pi_2'$}(D);
}$$
where we distinguish between the projections from $X\times X$ and $X\times Y$ in order to avoid confusion. This means that in this notation $C^X_T=\pi'_{1*}\left(\left(\_\right)\otimes \pi_{2}'^* T\right)$ and $C^X_{f^*T}=\pi_{1*}\left(\left(\_\right)\otimes \pi_2^* f^*T\right)$.

On objects and morphisms we have the following sequence of natural isomophisms
\begin{align*}
C^X_{f^*T}\left(\_\right)&= \pi_{1*} \left(\left(\_\right)\otimes \pi_{2}^*f^*T \right) &\text{Definition~\ref{definition CT}}\\
&\cong \pi'_{1*}\left(\id\times f\right)_* \left(\left(\_\right)\otimes \left(\id\times f\right)^*\pi_{2}'^*T \right)& \begin{matrix}\pi_1=\pi'_1\circ \left(\id\times f\right)\\f\circ \pi_2 =\pi'_2\circ \left(\id \times f\right)
\end{matrix}\\
&\cong \pi'_{1*} \left(\left(\id\times f\right)_*\left(\_\right)\otimes {\pi'_2}^* T \right) &\text{projection formula}\\
&=C^X_T\circ \left(\id\times f\right)_*\left(\_\right). &\text{Definition~\ref{definition CT}}
\end{align*}
 Both functors also induce the same $\Gamma$-action as we get analogously
\begin{align*}
\pi_{1*} \left(\id\otimes \pi_{2}^*f^*\gamma \right)&\cong \pi'_{1*}\left(\id\times f\right)_* \left(\id \otimes \left(\id\times f\right)^*\pi_{2}'^*\gamma \right)& \begin{matrix}\pi_1=\pi'_1\circ \left(\id\times f\right)\\f\circ \pi_2 =\pi'_2\circ \left(\id \times f\right)
\end{matrix}\\
&\cong \pi'_{1*} \left(\left(\id\times f\right)_*\id\otimes {\pi'_{2}}^*\gamma \right). &\text{projection formula}\\
\end{align*}
This means that the actions match up along the same natural isomorphisms and so $$C^X_{f^*T}\cong C^X_T \circ \left(\id\times f\right)_*$$ as claimed.
\end{proof}
\end{lemma}

\begin{remark}\label{remark lemma naturality C}
The above Lemma~\ref{naturality C} can be interpreted very naturally using Remark~\ref{CT can be interpreted via FM functors}. As $C^X_T$ sends an $M$ to the image of $T$ under the Fourier-Mukai functor $\Phi_M$ and we have by \cite[Exercise~5.12]{HuybrechtsFm} $\Phi_{\left(\id\times f\right)_* M}\cong \Phi_{M}\circ f^*$. In particular the two functors $C_{f^*T}^X$ and $C^X_T\circ \left(f \times \id\right)_*$ should be isomorphic.
\end{remark}

\begin{proposition}\label{composition cfT,Gamma}
Let $f:X\to Y$ be a proper morphism of schemes and $T\in \Coh\left(Y\right)_\Gamma$. Then we have 
$$ c_{f^*T,\Gamma}\left(M\right)=  C^X_T \circ\left(\id\times f\right)_*:\Ext^*_{X\times X}\left(\OO_\Delta, \Delta_* M\right)\to \Ext^*_{\Coh\left(X\right)_\Gamma}\left(f^*T,M\otimes f^*T\right). $$
\begin{proof}
By Proposition~\ref{cT=CT} we have 
$$c_{f^*T,\Gamma}\left(M\right)=C^X_{f^*T}: \Ext_{X\times X}^*\left(\OO_\Delta ,\Delta_* M\right)\to \Ext_{\Coh\left(X\right)_\Gamma}^*\left(T,M\otimes T\right)$$
 and by Lemma~\ref{naturality C} we get
$$c_{f^*T,\Gamma}\left(M\right)=C^X_{f^*T}=C^X_T\circ \left(\id\times f\right)_*:\Ext^*_{X\times X}\left(\OO_\Delta, \Delta_* M\right)
 \to \Ext^*_{\Coh\left(X\right)_\Gamma}\left(f^*T,M\otimes f^*T\right) $$
as claimed.
\end{proof}
\end{proposition}

\begin{remark}
The above result could be used to compute injectivity of the characteristic morphism if one can find an $\left(f,\id\right):X\times X \to X\times Y$ that is injective on $\Ext^*_{X\times X}\left(\OO_\Delta,\Delta_* M\right)$ and such that $Y$ admits a tilting bundle. However, the existence of such a morphism is not straight forward. In particular a closed immersion of a divisor $f: X \hookrightarrow \PP^n$ is in general not injective on $\Ext_X^i\left(\_,\_\right)$ as by the Grothendieck-Serre spectral sequence there might be correction terms arising in degrees $i>1$. 
\end{remark}

\section{Twisted Hodge diamonds give kernels in Hochschild cohomology}\label{Section Twisted Hodge diamonds}

We will show how twisted Hodge diamonds, and in particular their interior, can be used to understand the pushforward of Hochschild cohomology under the closed embedding of a smooth projective hypersurface of degree $d$.

Throughout this chapter we will follow Br\"uckmann`s paper "Zur Kohomologie von projektiven Hyperfl\"achen" \cite{Brueckmann1974} for computations.

\begin{definition}\label{twisted hodge diamond} Let $X$ be a projective scheme of dimension $n$ and let $\OO_X\left(1\right)$ be a very ample line bundle. Then we define the twisted Hodge numbers of $X$ to be $$\h^{i,j}_p \left(X\right):=\dim \HH^{j}\left(X, \Omega_X^i\left(p\right)\right).$$
Similarly to ordinary Hodge numbers the twisted Hodge numbers can be arranged in a twisted Hodge diamond:
$$
\tikz[heighttwo,xscale=1.5,yscale=1.5,baseline]{
\node (hnn) at (2,4) {$\h^{n,n}_{p}\left(X\right)$};
\node (hn0) at (0,2) {$\h^{n,0}_{p}\left(X\right)$};
\node (h00) at (2,0) {$\h^{0,0}_{p}\left(X\right).$};
\node (h0n) at (4,2) {$\h^{0,n}_{p}\left(X\right)$};
\draw[-]
(h00)  edge (hn0)
(h00) edge (h0n)
(hn0) edge (hnn)
(h0n) edge (hnn);}
$$
We will drop the $X$ if the space is clear from context.
\end{definition}

\begin{lemma}\label{decomposition Hochschild into twisted hodge numbers}
Let $X$ be a smooth projective scheme of dimension $n$ with canonical sheaf of form $\OO_{X}\left(t\right)$. Then we have $$\HHH^{m}\left(X, \OO_X\left(p\right)\right)\cong \bigoplus_{i=0}^{n} \HH^{i-m+n}\left(X, \Omega_X^i\right(t-p\left)\right).$$ 
In particular this gives
$$\dim\HHH^{m}\left(X, \OO_X\left(p\right)\right)=\sum_{i=0}^{n} \h^{i,i-m+n}_{t-p}\left(X\right).$$
\begin{proof}
We compute, using $\omega_X\cong \OO_X\left(t\right)$ and the Hochschild-Kostant-Rosenberg (HKR) isomorphism \cite{Swan1996}
\begin{align*}
\HHH^{m}\left(X, \OO_X\left(p\right)\right)&\cong \bigoplus_{i=0}^{n} \Ext_{X}^{m-i} \left(\Omega_X^i,\OO_X\left(p\right)\right)&\text{HKR}\\
&\cong \bigoplus_{i=0}^{n} \Ext_{X}^{n-m+i} \left(\OO_X\left(p\right),\Omega_X^i\left(t\right)\right)^*&\text{Serre Duality }\\
&\cong \bigoplus_{i=0}^{n} \Ext_{X}^{n-m+i} \left(\OO_{X},\Omega_X^i\left(t-p\right)\right)^*&\text{twisting on both sides}\\
&\cong \bigoplus_{i=0}^{n} \HH^{n-m+i}\left( X,\Omega_X^i\left(t-p\right)\right)^{\mathrlap{*}}. & \Ext^j_{X}\left(\OO_{X}, \_\right)\cong \HH^j\left(X, \_ \right)
\end{align*}
Applying dimension on both sides gives 
$$\dim\HHH^{m}\left(X, \OO_X\left(p\right)\right)=\sum_{i=0}^{n} \h^{i,i-m+n}_{t-p}\left(X\right)$$ 
as desired.
\end{proof}
\end{lemma}

\begin{remark}\label{Hochschild cohomology via twisted Hodge diamonds}
By Lemma \ref{decomposition Hochschild into twisted hodge numbers} one can compute $\dim\HHH^{m}\left(X, \OO_X\left(p\right)\right)$ as the sum over the $m$-th column in the $t-p$ twisted Hodge diamond,
$$
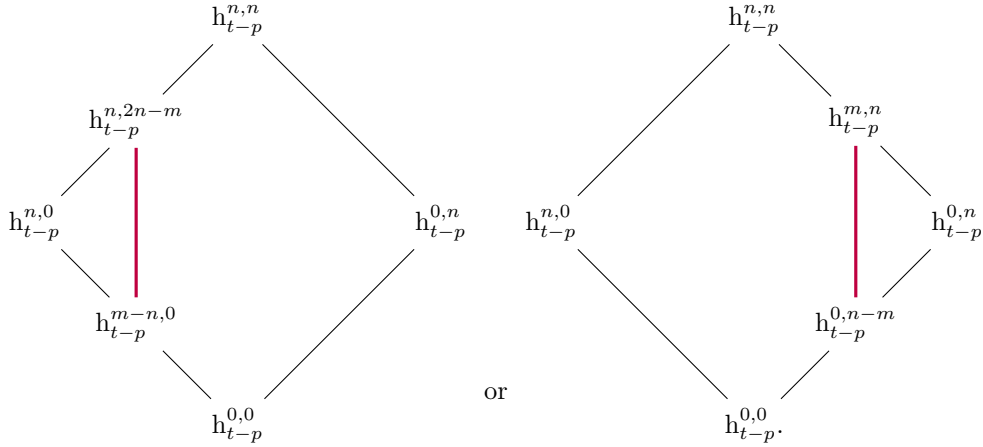
\tikz[heighttwo,xscale=1.5,yscale=1.5,baseline]{
\node (hnn) at (2,4) {$\h^{n,n}_{t-p}$};
\node (hn0) at (0,2) {$\h^{n,0}_{t-p}$};
\node (h00) at (2,0) {$\h^{0,0}_{t-p}$};
\node (h0n) at (4,2) {$\h^{0,n}_{t-p}$};
\node (hn2n-m) at (1,3) {$\h^{n,2n-m}_{t-p}$};
\node (hm-n0) at (1,1) {$\h^{m-n,0}_{t-p}$};

\draw[-]
(h00) edge (hm-n0)
(hm-n0) edge (hn0)
(h00) edge (h0n)
(hn0) edge (hn2n-m)
(hn2n-m)edge (hnn)
(h0n) edge (hnn);
\draw[-,very thick,purple]
(hn2n-m)edge (hm-n0);

}\text{ or }\tikz[heighttwo,xscale=1.5,yscale=1.5,baseline]{
\node (hnn) at (2,4) {$\h^{n,n}_{t-p}$};
\node (hn0) at (0,2) {$\h^{n,0}_{t-p}$};
\node (h00) at (2,0) {$\h^{0,0}_{t-p}.$};
\node (h0n) at (4,2) {$\h^{0,n}_{t-p}$};
\node (hmn) at (3,3) {$\h^{m,n}_{t-p}$};
\node (h0n-m) at (3,1) {$\h^{0,n-m}_{t-p}$};

\draw[-]
(h00) edge (hn0)
(h0n-m) edge (h0n)
(h00) edge (h0n-m)
(h0n) edge (hmn)
(hmn)edge (hnn)
(hn0) edge (hnn);
\draw[-,very thick,purple]
(h0n-m)edge (hmn);
}$$
\end{remark}

\subsection{The Hochschild Cohomology of a smooth Hypersurface}
We will use the computations in \cite{Brueckmann1974} and Lemma~\ref{Hochschild cohomology via twisted Hodge diamonds} to compute the Hochschild cohomology of $X$.

\begin{lemma}\label{upper, middle and bottom part}
Let $X\hookrightarrow \PP^{n+1}$ be a smooth degree $d$ hypersurface. Then
 $$\h_p^{i,j}\left(X\right)=0$$
if $\left(i,j\right)$ is not of the form $\left(i,0\right),\left(i,n\right),\left(i,n-i\right),\left(i,i\right)$, with $0\le i\le n$. And we have for $\left(i,i\right)$:
$$\h^{i,i}_p\left(X\right)=\delta_{p,0} \quad \text{ if } i\notin\left\{0,\frac{n}{2},n\right\}.$$
Moreover we get
\begin{equation}\label{formula middle part}
\h^{i,n-i}_p\left(X\right)=\sum_{\mu=0}^{n+2}\left(-1\right)^\mu\binom{n+2}{\mu}\binom{-p+id -\left(\mu -1\right)\left(d-1\right)}{n+1}+\delta_{p,0}\delta_{i,n-i}.
\end{equation}
\begin{proof}
 First of all we can assume $0\le i,j\le n$ as outside of that range we have $\Omega_X^i\left(p\right) = 0$, respectively $\HH^j\left(X,\Omega_X^i\left(p\right)\right)=0$ for dimension reasons. 

By \cite[Satz 2,(42),(40),(38) and (39)]{Brueckmann1974} we have for $0<i<n $
\begin{align*}
\h^{i,j}_p\left(X\right)= \begin{cases} \binom{-p-1}{n-i}\binom{-p+1+i}{1+i} +\sum_{\mu =1}^{n-i+1} \left(-1\right)^\mu \binom{n+2}{\mu} \binom{-p-\mu\left(d-1\right)+i}{n+1}  &\quad\text{if }j=n\\
 \sum_{\mu=0}^{n+2}\left(-1\right)^{\mu}\binom{n+2}{\mu}\binom{-p+id -\left(\mu -1\right)\left(d-1\right)}{n+1} +\delta_{p,0}\delta_{i,j}&\quad\text{if } i+j=n \\
\binom{p-1}{i}\binom{p+n+1-i}{n+1-i} +\sum_{\mu =1}^{i+1} \left(-1\right)^\mu \binom{n+2}{\mu} \binom{p+n-\mu\left(d-1\right)-i}{n+1} &\quad\text{if } j=0\\
\delta_{p,0}  &\quad\text{if } i=j\notin \left\{0,n\right\}\\
0 &\quad \text{else}.
\end{cases} 
\end{align*}
So only the cases for $i\in \left\{0,n\right\}$ remain. Now \cite[Lemma 5]{Brueckmann1974} gives for $j\notin \left\{0,n\right\}$
$$\h^{0,j}_p\left(X\right)=0=\h^{n,j}_p\left(X\right).$$
Which finishes the claim.

\end{proof}
\end{lemma}

\begin{remark}
By Lemma~\ref{upper, middle and bottom part} the $p$-twisted Hodge diamond of a smooth degree $d$ hypersurface has the shape:

$$
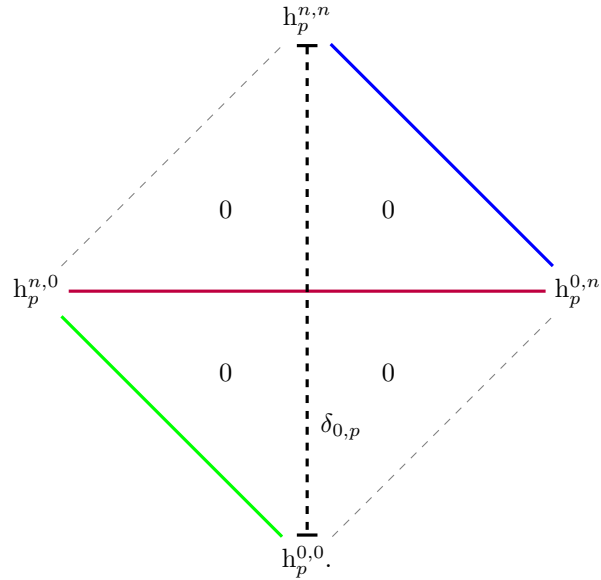
\tikz[heighttwo,xscale=2,yscale=2,baseline]{
\node (hnn) at (2,4) {$\h^{n,n}_{p}$};
\node (hn0) at (0,2) {$\h^{n,0}_{p}$};
\node (h00) at (2,0) {$\h^{0,0}_{p}.$};
\node (h0n) at (4,2) {$\h^{0,n}_{p}$};
\node (delta) at (2.25,1){$\delta_{0,p}$};
\node (01) at (1.4,1.4){$0$};
\node (01) at (2.6,2.6){$0$};
\node (01) at (1.4,2.6){$0$};
\node (01) at (2.6,1.4){$0$};
\draw[dashed,gray]
(h00) edge (hn0)
(h0n) edge (hnn)
(h00) edge (h0n)
(hn0) edge (hnn);
\draw[very thick,purple]
(hn0) edge  (h0n);
\draw[very thick, blue]
(hnn) edge  (h0n);
\draw[very thick, green]
(h00) edge  (hn0);
\draw[|-|,very thick,dashed]
(hnn) edge  (h00)
;
}$$

In particular the only non-trivial entries appear along the indicated lines. More precisely we have along the blue lines the values for $\h^{i,n}\left(X\right)$, along the red lines the values $\h^{i,n-i}\left(X\right)$ and along the green line $\h^{i,0}\left(X\right)$. Furthermore the dashed line disappears if $p\neq 0$ as these are the Kronecker deltas $\delta_{p,0}$.

\end{remark}

\begin{proposition}\label{formulas for upper, middle and bottom part}
Let $X\hookrightarrow \PP^{n+1}$ be the embedding of a smooth degree $d$ hypersurface. Then the following formulas hold.
\begin{align*}
\h^{i,0}_{p}\left(X\right)&=\h^{n-i,n}_{-p}\left(X\right)\\
\h^{i,n-i}_{p}\left(X\right)&=\h^{i-1,n+1-i}_{p-d}\left(X\right)&i\notin \left\{0,1,n\right\},\;p\neq 0 \\
\h^{i,n+1}_{p-d}\left(\PP^{n+1}\right)-\h^{i,n+1}_{p}\left(\PP^{n+1}\right)&=\h^{i,n}_p\left(X\right) +\h^{i-1,n}_{p-d} \left(X\right)&i\notin \left\{0,1,n\right\} \\
\h^{i,0}_{p}\left(\PP^{n+1}\right)-\h^{i,0}_{p-d}\left(\PP^{n+1}\right)&=\h^{i,0}_p\left(X\right) +\h^{i-1,0}_{p-d} \left(X\right)&i\notin \left\{0,1,n\right\}. 
\end{align*}
\begin{proof}
We compute for the first equation:
\begin{align*}
\h^{i,0}_p\left(X\right)&=\dim\HH^0\left(X,\Omega_X^i\left(p\right)\right)&\text{definition}\\
&=\dim\Ext^0\left(\OO_X,\Omega_X^i\left(p\right)\right)&\Ext^* \left(\OO_X,\_\right)\cong \HH^*\left(X,\_\right)\\
&=\dim \Ext^n\left(\Omega_X^i\left(p\right),\Omega_X^n\right)&\text{Serre Duality}\\
&=\dim \Ext^n\left(\OO_X,\Omega_X^{n-i}\left(-p\right)\right)\\
&=\dim \HH^{n}\left(X,\Omega_X^{n-i}\left(-p\right)\right)&\Ext^* \left(\OO_X,\_\right)\cong \HH^*\left(X,\_\right)\\
&=\h^{n-i,0}_p\left(X\right).&\text{definition}
\end{align*}
For the second equation we we have by $\left(\ref{formula middle part}\right)$ the following identity
\begin{align*}
\h^{i,n-i}_{p}\left(X\right)&= \sum_{\mu=0}^{n+2}\left(-1\right)\binom{n+2}{\mu}\binom{-p+id -\left(\mu -1\right)\left(d-1\right)}{n+1}\\
 &=\sum_{\mu=0}^{n+2}\left(-1\right)\binom{n+2}{\mu}\binom{-p+d-d+id -\left(\mu -1\right)\left(d-1\right)}{n+1}\\
 &=\sum_{\mu=0}^{n+2}\left(-1\right)\binom{n+2}{\mu}\binom{-p+d+\left(i-1\right)d -\left(\mu -1\right)\left(d-1\right)}{n+1}\\
 &=\h^{i-1,n+1-i}_{p-d}\left(X\right).
\end{align*}
And for the last two Br\"uckmann gives the formula \cite[(31)]{Brueckmann1974}, which together with \cite[Satz 2]{Brueckmann1974} gives both
\begin{align*}
\h^{i,n}_{p-d}\left(X\right) &=\h^{i+1,n}_{p-d}\left(\PP^{n+1}\right) - \h^{i+1,n}_{p}\left(\PP^{n+1}\right)-\h^{i+1,n}_{p} \left(X\right)\\
\h^{i,0}_p\left(X\right) &= \h^{i,0}_{p}\left(\PP^{n+1}\right)- \h^{i,0}_{p}\left(\PP^{n+1}\right)-\h^{i-1,0}_{p-d} \left(X\right).
\end{align*}
After rearranging, these are
\begin{align*}
\h^{i,n}_{p-d}\left(X\right)+\h^{i+1,n}_{p}\left(X\right) &=\h^{i+1,n}_{p-d}\left(\PP^{n+1}\right) - \h^{i+1,n}_{p}\left(\PP^{n+1}\right) \\
\h^{i,0}_p\left(X\right)+\h^{i-1,0}_{p-d} \left(X\right) &= \h^{i,0}_{p}\left(\PP^{n+1}\right)- \h^{i,0}_{p}\left(\PP^{n+1}\right).
\end{align*}
Index shifting in the first equation gives
\begin{align*}
\h^{i,n}_{p-d}\left(\PP^{n+1}\right)-\h^{i,n}_{p}\left(\PP^{n+1}\right)&=\h^{i,n}_p\left(X\right) +\h^{i-1,n}_{p-d} \left(X\right)\\
\h^{i,0}_{p}\left(\PP^{n+1}\right)-\h^{i,0}_{p-d}\left(\PP^{n+1}\right)&=\h^{i,0}_p\left(X\right) +\h^{i-1,0}_{p-d} \left(X\right),
\end{align*}
as claimed.
\end{proof}
\end{proposition}

We can use Lemma~\ref{upper, middle and bottom part} together with Lemma~\ref{decomposition Hochschild into twisted hodge numbers} to compute the dimensions of $\HHH^m\left(X,\OO_X\left(p\right)\right)$:

\begin{corollary}\label{dimension Hochschild cohomology as upper middle and bottom part}
Let $X$ be a smooth $n$-dimensional hypersurface of degree $d$ and let $t=d-n-2$. Then we have $$\dim\HHH^m\left(X,\OO_X\left(p\right)\right)=\begin{cases}
\h^{0,n}_{t-p}\left(X\right) &\text{for } m=0\\
\h^{m-n,0}_{t-p}\left(X\right)+\h_{t-p}^{\frac{m}{2},n-\frac{m}{2}}\left(X\right)+\h_{t-p}^{m,n}\left(X\right) + \left(n-2\right)\delta_{t,p}\delta_{m,n}&\text{for } 0< m < 2n \text{ even}\\
\h_{t-p}^{m-n,0}\left(X\right)+\h_{t-p}^{m,n}\left(X\right) + \left(n-1\right)\delta_{t,p}\delta_{m,n}&\text{for }0< m < 2n \text{ odd} \\
\h^{n,0}_{t-p}\left(X\right) &\text{for } m=2n\\
0 &\text{else.}
\end{cases}$$
\begin{proof}
Observe that we have $\OO_X\left(t\right)\cong \omega_X$ and that we can assume $0\le m\le 2n$ for dimension reasons. By Lemma~\ref{decomposition Hochschild into twisted hodge numbers} we have $$\dim\HHH^{m}\left(X, \OO_X\left(p\right)\right)=\sum_{i=0}^{n} \h^{i,i-m+n}_{t-p}\left(X\right).$$
So we can use Lemma~\ref{upper, middle and bottom part} to compute every summand. In particular we get for $n=0$ and $n=2n$
\begin{align*}
\dim\HHH^0\left(X,\OO_X\left(p\right)\right)&=\h^{0,n}_{t-p}\left(X\right) \\
\dim\HHH^{2n}\left(X,\OO_X\left(p\right)\right)&=\h^{n,0}_{t-p}\left(X\right)
\end{align*}
as there only one summand appears.

Now for $0<m<2n$ we can use Lemma~\ref{upper, middle and bottom part} to get for
$m \neq n$
\begin{align*}
\dim\HHH^{m}\left(X, \OO_X\left(p\right)\right)&=\sum_{i}^{\mathclap{m-n}} \h^{i,i-n+m}_{t-p}\left(X\right)\\
&=\sum_{\mathclap{i-j=m-n}} \h^{i,j}_{t-p}\left(X\right)\\
&=\begin{cases}\h^{m-n,0}_{t-p}\left(X\right) + \h^{\frac{2n-m}{2},\frac{m}{2}}_{t-p}\left(X\right) + \h^{m,n}_{t-p}\left(X\right)  &\text{for }m \text{ even} \\
\h^{m-n,0}_{t-p}\left(X\right) + \h^{m,n}_{t-p}\left(X\right)  &\text{for }m \text{ odd}.
\end{cases}
\end{align*}
For $m=n$ all the above calculations still hold, however, we get for all $\left(i,i\right)$ with $0< i < n$ and $\left(i,i\right)\neq\left(\frac{n}{2},\frac{n}{2}\right)$ an additional $\delta_{0,t-p}=\delta_{t,p}$, which means that 
$$\dim\HHH^n\left(X,\OO_X\left(p\right)\right)=\begin{cases}
\h_{t-p}^{0,0}\left(X\right)+\h_{t-p}^{\frac{n}{2},\frac{n}{2}}\left(X\right)+\h_{t-p}^{n,n}\left(X\right) + \left(n-2\right)\delta_{t,p}\delta_{m,n}&\text{for } m \text{ even}\\
\h_{t-p}^{0,0}\left(X\right)+\h_{t-p}^{n,n}\left(X\right) + \left(n-1\right)\delta_{t,p}\delta_{m,n}&\text{for }m \text{ odd} 
\end{cases}$$
as claimed.
\end{proof}
\end{corollary}

\subsection{The Hochschild Cohomology of the direct Image}
Since we want to control the pushforward in Hochschild cohomology we will use computations by \cite{Brueckmann1974} to understand the Hochschild dimensions of the direct image of a line bundle under a smooth embedding.

\begin{lemma}\label{computation Hochschild cohomology direct image}
Let $f:X\hookrightarrow Y$ be an embedding of a smooth $n$-dimensional degree $d$ hypersurface and assume $\omega_X\cong \OO\left(t\right)$. Then we have
$$\HHH^{m}\left(Y,f_*\OO_X\left(p\right)\right)\cong \bigoplus_{i=0}^{\dim Y} \HH^{n-m+i}\left(X,f^*\Omega^i_Y\left(t-p\right)\right).$$
\begin{proof}
We can compute, using  $\omega_X \cong\OO_X\left(t\right)$ and the Hochschild-Kostant-Rosenberg Isomorphism (HKR) \cite{Swan1996}
\begin{align*}
\HHH^{m}\left(Y,f_*\OO_X\left(p\right)\right)&\cong \bigoplus_{i=0}^{\dim Y} \Ext_{Y}^{m-i} \left(\Omega^i_Y,f_*\OO_X\left(p\right)\right)&\text{HKR}\\
&\cong \bigoplus_{i=0}^{\dim Y} \Ext_{X}^{m-i} \left(f^*\Omega^i_Y,\OO_X\left(p\right)\right)&f^* \dashv f_*\\
&\cong \bigoplus_{i=0}^{\dim Y} \Ext_{X}^{n-m+i} \left(\OO_X\left(p\right),f^*\Omega^i_Y\left(t\right)\right)^*&\text{Serre Duality }\\
&\cong \bigoplus_{i=0}^{\dim Y} \Ext_{X}^{n-m+i} \left(\OO_{X},f^*\Omega^i_Y\left(t-p\right)\right)^*&\text{twisting on both sides}\\
&\cong \bigoplus_{i=0}^{\dim Y} \HH^{n-m+i}\left( X,f^*\Omega^i_Y\left(t-p\right)\right)^*& \Ext^j_{X}\left(\OO_{X}, \_\right)\cong \HH^j\left(X, \_ \right)
\end{align*}
as desired.
\end{proof}
\end{lemma}

\begin{lemma}\label{computation pullback}
Let $f:X\hookrightarrow\PP^{n+1}$ be a closed embedding of a smooth degree $d$ hypersurface. Then we have for $\left(i,j\right)\notin\left\{ \left(0,0\right),\left(0,n\right),\left(0,n+1\right),\left(n,1\right)\left(n,n\right),\left(n,n+1\right)\right\}$
\begin{align*}
\dim\HH^j\left(X, f^*\Omega_{\PP^{n+1} }^i\left(p\right)\right)&=\begin{cases}\h_p^{i,0}\left(X\right)+\h_{p-d}^{i-1,0}\left(X\right)&\text{if j=0} \\
\h^{i,n}_p\left(X\right) +\h^{i-1,n}_{p-d}\left(X\right) &\text{if } j=n\\
\delta_{p,0} &\text{if } i=j\notin \left\{0,n\right\} \\
\delta_{p,d} &\text{if } i-1=j\notin \left\{0,n\right\}\\
0 &\text{else.}
\end{cases}
\end{align*}
Moreover, we get 
\begin{align*}
\dim\HH^0\left(X, f^*\Omega_{\PP^{n+1} }^n\left(p\right)\right)&= \h_p^{n,0}\left(X\right) +\h_{p-d}^{n-1,0}\left(X\right)-\h_{p-d}^{n-1,1}\left(X\right) \\
\dim\HH^n\left(X, f^*\Omega_{\PP^{n+1} }^1\left(p\right)\right)&= \h_p^{1,n}\left(X\right) +\h_{p-d}^{0,n}\left(X\right)-\h_{p}^{1,n-1}\left(X\right) \\
\dim \HH^0\left(X,f^*\Omega_{\PP^{n+1} }^0 \left(p\right)\right)&=\h_p^{0,0}\left(X\right)\\
\dim\HH^n\left(X, f^*\Omega_{\PP^{n+1} }^0\left(p\right)\right)&= \h_p^{0,n}\left(X\right) \\
\dim\HH^0\left(X, f^*\Omega_{\PP^{n+1}}^{n+1}\left(p\right)\right)&= \h_{p-d}^{n,0}\left(X\right)\\
\dim\HH^n\left(X, f^*\Omega_{\PP^{n+1}}^{n+1}\left(p\right)\right)&= \h_{p-d}^{n,n}\left(X\right).
\end{align*}
\begin{proof}
First observe that $\Omega_{\PP^{n+1}}^{i}=0$ for $i>n+1$ and $i<0$. In particular we can assume $0\le i \le n+1$ and for dimension reasons we can additionally assume $0\le j\le n$. Furthermore, \cite[Lemma 5 and Lemma 6]{Brueckmann1974} give for $0<j<n$
\begin{align*}
\dim\HH^j\left(X,f^*\Omega^0\left(p\right)\right)&=0\\
\dim\HH^j\left(X,f^*\Omega^{n+1}\left(p\right)\right)&=0.
\end{align*}

By \cite[Satz 1, Lemma 6, (21) and (25) ]{Brueckmann1974} we have for $0<i< n+1$
\begin{align*}
\dim\HH^j\left(X, f^*\Omega_{\PP^{n+1}}^i\left(p\right)\right)&=\begin{cases}\h_p^{i,0}\left(\PP^{n+1}\right)-\h_{p-d}^{i,0}\left(\PP^{n+1}\right)&\text{if j=0} \\
\h_{p}^{i,n+1}\left(\PP^{n+1}\right) -\h_{d-p}^{i,n+1}\left(\PP^{n+1}\right) &\text{if } j=n\\
\delta_{p,0}\delta_{i,j}+\delta_{p,d}\delta_{i-1,j} &\text{if } j\notin \left\{0,n\right\}. \\
\end{cases}
\end{align*}
Which gives after applying Proposition~\ref{formulas for upper, middle and bottom part} for $i\notin \left\{1,n\right\}$
\begin{align*}
\dim\HH^j\left(X, f^*\Omega_{\PP^{n+1}}^i\left(p\right)\right)&=\begin{cases}\h_p^{i,0}\left(X\right)+\h_{p-d}^{i-1,0}\left(X\right)&\text{if j=0} \\
\h_{p}^{i,n}\left(X\right) +\h_{d-p}^{i-1,n}\left(X\right)  &\text{if } j=n\\
\delta_{p,0}\delta_{i,j}+\delta_{p,d}\delta_{i-1,j} &\text{if } j\notin \left\{0,n\right\}. \\
\end{cases}
\end{align*}

Now for the special cases:

We start with the case of $i\in \left\{1,n\right\}$. By the discussion above we have
\begin{align*}
\dim\HH^0\left(X,f^*\Omega_{\PP^{n+1}}^n\left(p\right)\right)&=\h^{n,0}_p\left(\PP^{n+1}\right)-\h_{p-d}^{n,0}\left(\PP^{n+1}\right)\\
\dim\HH^n\left(X,f^*\Omega_{\PP^{n+1}}^1\left(p\right)\right)&=\h^{1,n+1}_p\left(\PP^{n+1}\right)-\h_{p-d}^{1,n+1}\left(\PP^{n+1}\right).
\end{align*} 
This turns, using \cite[(31), (33), Satz 2 and Lemma 5]{Brueckmann1974} into
\begin{align*}
\dim\HH^0\left(X, f^*\Omega_{\PP^{n+1}}^n\left(p\right)\right)&= \h_p^{n,0}\left(X\right) +\h_{p-d}^{n-1,0}\left(X\right)-\h_{p-d}^{n-1,1}\left(X\right) \\
\dim\HH^n\left(X, f^*\Omega_{\PP^{n+1}}^1\left(p\right)\right)&= \h_p^{1,n}\left(X\right) +\h_{p-d}^{0,n}\left(X\right)-\h_{p}^{1,n-1}\left(X\right).
\end{align*} 
So only the cases for $i=0$ and $i=n+1$ remain:

For $i=0$ we have $f^*\Omega_{\PP^{n+1}}^0\left(p\right)\cong\OO_X\left(p\right)$, so we can apply Lemma~\ref{twisted hodge diamond} to get:
\begin{align*}
\dim\HH^0\left(X,\OO_X\left(p\right)\right)&=\h^{0,0}_p\left(X\right)\\
\dim\HH^n\left(X,\OO_X\left(p\right)\right)&=\h^{0,n}_p\left(X\right).
\end{align*}

Now for $i=n+1$ we have
\begin{align*}
f^*\Omega_{\PP^{n+1}}^{n+1}\left(p\right)&\cong f^* \OO_X\left(p-n-2\right)\\
&\cong \OO_X\left( d-n-2+p-d\right)\\
&\cong \Omega_X^{n}\left(p-d\right).
\end{align*}
And so we get by Definition~\ref{twisted hodge diamond}
\begin{align*}
\dim\HH^0\left(X, \Omega_X^{n}\left(p-d\right)\right)&= \h^{n,0}_{p-d}\left(X\right)\\
\dim\HH^n\left(X, \Omega_X^{n}\left(p-d\right)\right)&= \h^{n,n}_{p-d}\left(X\right)
\end{align*}
as claimed.
\end{proof}
\end{lemma}

\begin{remark}
If we arrange the computation of the cohomology dimensions from Lemma~\ref{computation pullback} analogously to a twisted Hodge diamond, we get that it is of the shape
\begin{small}
$$
\tikz[heighttwo,xscale=0.85,yscale=0.85,baseline]{
\node (hnn) at (5,10) {$\h^n\left(X, f^*\Omega^n\left(p\right)\right)$};
\node[red] (hn0) at (0,5) {$\mathllap{\text{\begin{tiny}$\h_p^{n,0}\left(X\right) +\h_{p-d}^{n-1,0}\left(X\right)-\h_{p-d}^{n-1,1}\left(X\right)$\end{tiny}}=}\h^0\left(X, f^*\Omega^n\left(p\right)\right)$};
\node (h00) at (5,0) {$\mathllap{\text{\begin{tiny}$\h^{0,0}_p\left(X\right)$\end{tiny}}=}\h^0\left(X, f^*\Omega^0\left(p\right)\right),$};
\node (h0n) at (10,5) {$\h^n\left(X, f^*\Omega^0\left(p\right)\right) \mathrlap{=\text{\begin{tiny}$\h^{n,0}_{p-d}\left(X\right)$\end{tiny}}}$};
\node[red] (h1n) at (8.4,6.6) {$\h^n\left(X, f^*\Omega^1\left(p\right)\right)\mathrlap{=\text{\begin{tiny}$ \h_p^{1,n}\left(X\right) +\h_{p-d}^{0,n}\left(X\right)-\h_{p}^{1,n-1}\left(X\right)$\end{tiny}}}$};
\node (h10) at (3.4,1.6) {$\h^0\left(X, f^*\Omega^1\left(p\right)\right)$};
\node (hn+10) at (-1.6,6.6) {$\mathllap{\text{\begin{tiny}$\h^{n,0}_{p-d}\left(X\right)$\end{tiny}}=}\h^{n}\left(X, f^*\Omega^{n+1}\left(p\right)\right)$};
\node (hn+1n) at (3.4,11.6) {$\h^{n}\left(X, f^*\Omega^{n+1}\left(p\right)\right)\mathrlap{=\text{\begin{tiny}$\h^{n,n}_{p-d}\left(X\right)$\end{tiny}}}$};

\draw[dashed,gray]
(hnn) edge (hn0)
(h00) edge (h0n)
(h10) edge (h1n)
(hn+1n) edge (hn+10);
\draw[very thick,green]
(hn0) edge (hn+10)
(h00) edge (h10)
(h10) edge node[below left]{$\h^{*,0}_p+\h^{*-1,0}_{p-d}$} (hn0);
\draw[|-|,very thick,dashed]
(hnn) edge node[right]{$\delta_{p,0}$} (h00)
;
\draw[-,very thick,blue]
(hnn) edge node[above right]{$\h^{*,n}_{p}+\h^{*-1,n}_{p-d}$} (h1n)
(hnn) edge (hn+1n)
(h1n) edge (h0n);
\draw[|-|,very thick,dashed]
(hn+1n) edge node[left]{$\delta_{p,d}$} (h10)
;
}$$
\end{small}
where apart from the two special cases \color{red}
\begin{align*}\h^0\left(X, f^*\Omega_{\PP^{n+1}}^{n}\left(p\right)\right) &=\h^{n,0}_p\left(X\right)+\h^{n-1,0}_{p,d}\left(X\right)-\h^{n-1,1}_{d}\left(X\right)  \\
 \h^n\left(X, f^*\Omega_{\PP^{n+1}}^{1}\left(p\right)\right) &=\h^{1,n}_p\left(X\right)+\h^{0,n}_{p-d}\left(X\right)-\h^{1,n-1}_{p-d}\left(X\right)
\end{align*}  
\color{black}
the only non-trivial entries are along the indicated lines. There we have
 \begin{align*}
 \color{green}\h^0\left(X,f^*\Omega_{\PP^{n+1}}^i\right)&\color{green}=\h^{*,0}_p\left(X\right)+\h^{*-1,0}_{p-d}\left(X\right)\\
\color{blue}\h^n\left(X,f^*\Omega_{\PP^{n+1}}^i\right)&\color{blue}=\h^{*,n}_p\left(X\right)+\h^{*-1,n}_{p-d}\left(X\right)
\end{align*}
and along the two vertical diagonals we have 
\begin{align*}
\h^i\left(X,f^*\Omega_{\PP^{n+1}}^i\right)&=\delta_{p,0}\\
\h^i\left(X,f^*\Omega_{\PP^{n+1}}^{i+1}\right)&=\delta_{p,d}.
\end{align*}
Observe that this has the shape of the $p$ and $p-d$ twisted Hodge diamond for $X$ laid on top of each other with the interior middle line removed.

Since we will focus on the case $p>d$ we will be able to ignore the dashed lines.
\end{remark}

\begin{proposition}\label{Hochschild dimension direct image explicit}
Let $f:X\hookrightarrow \PP^{n+1}$ be a smooth $n$-dimensional hypersurface of degree $d$ and set $t=d-n-2$. Then we have for $ m\notin\left\{1,2n\right\}$:
\begin{align*}
\dim\HHH^m\left(\PP^{n+1},f_*\OO_X\left(p\right)\right)&=\h^{0,m-n}_{t-p}\left(X\right) +\h^{0,m-n-1}_{t-p-d}\left(X\right) +\h^{n,m}_{t-p}\left(X\right) +\h^{n,m-1}_{t-p-d}\left(X\right) + \left(n-1\right)\left(\delta_{d,p}\delta_{m,n+1}+\delta_{0,p}\delta_{m,n}\right)
\end{align*}
and for $m=1$, $m=2n$:
\begin{align*}
\dim\HHH^1\left(\PP^{n+1},f_*\OO_X\left(p\right)\right)&=\h_{t-p}^{1,n}\left(X\right)+\h_{t-p-d}^{0,n}\left(X\right)- \h_{t-p}^{1,n-1}\left(X\right)\\
\dim\HHH^{2n}\left(\PP^{n+1},f_*\OO_X\left(p\right)\right)&=\h_{t-p}^{n,0}\left(X\right)+\h_{t-p-d}^{n-1,0} \left(X\right) - \h_{t-p-d}^{n-1,1}\left(X\right).
\end{align*}
\begin{proof}
We will compute the cases separately using $\OO_X\left(t\right)\cong \omega_X$:

We can use Lemma \ref{computation Hochschild cohomology direct image} to get for $m\notin\left\{1,2n\right\}$
\begin{align*}
\dim\HHH^m\left(\PP^{n+1},f_*\OO_X\left(p\right)\right)&=\sum_{i=0}^{n}  \dim\HH^{n-m+i}\left(X,f^*\Omega_{\PP^{n+1}}^{i}\left(t-p\right)\right)\\
&=\sum_{\mathclap{i-j=m-n}} \dim\HH^{j}\left(X,f^*\Omega_{\PP^{n+1}}^{i}\left(t-p\right)\right)\\
&=\h^{0,m-n}_{t-p}\left(X\right) +\h^{0,m-n-1}_{t-p-d}\left(X\right) +\h^{n,m}_{t-p}\left(X\right) +\h^{n,m-1}_{t-p-d}\left(X\right) + \left(n-1\right)\left(\delta_{d,p}\delta_{m,n+1}+\delta_{0,p}\delta_{m,n}\right).
\end{align*}
For $m=1$ we similarly get
\begin{align*}
\dim\HHH^1\left(\PP^{n+1},f_*\OO_X\left(p\right)\right)&=\sum_{i=0}^{n+1} \dim\HH^{n-1+i}\left(X,f^*\Omega_{\PP^{n+1}}^{i}\left(t-p\right)\right)\\
&= \dim\HH^n\left(X,f^*\Omega_{\PP^{n+1}}^{1}\left(t-p\right)\right)\\
&=\h_{t-p}^{1,n}\left(X\right)+\h_{t-p-d}^{0,n}\left(X\right)- \h_{t-p}^{1,n-1}\left(X\right).
\end{align*}
And for $m=2n$
\begin{align*}
\dim\HHH^{2n}\left(\PP^{n+1},f_*\OO_X\left(p\right)\right)&=\sum_{i=0}^{n+1} \dim\HH^{-n-1+i}\left(X,f^*\Omega_{\PP^{n+1}}^{i}\left(t-p\right)\right)\\
&= \dim\HH^0\left(X,f^*\Omega_{\PP^{n+1}}^{n+1}\left(t-p\right)\right)\\
&=\h_{t-p}^{n,0}\left(X\right)+\h_{t-p-d}^{n-1,0} \left(X\right) - \h_{t-p-d}^{n-1,1}\left(X\right)
\end{align*}
which finishes the claim.
\end{proof}
\end{proposition}

\begin{remark}
Since we will be able to assume $p\notin \left\{0,d\right\}$ in the next section we will exclude these cases. However, all of the following proofs and arguments still hold in these cases, one just needs to keep track of the Kronecker deltas in $\dim \HHH^n\left(X,\OO_X\left(p\right)\right)$.
\end{remark}

\begin{proposition}\label{Dimensions match up apart from middle part}
Let $f:X\hookrightarrow \PP^{n+1}$ be a smooth $n$-dimensional hypersurface of degree $d$, let $t=d-n-2$. Then we have for all $p \in \Z$ such that $t-p \notin \left\{0,d\right\}$ that $\dim\HHH^m\left(\PP^{n+1},f_*\OO_X\left(p\right)\right)$ is given by:
$$
\begin{cases}
\dim \HHH^0\left(X,\OO_X\left(p\right)\right) & m=0\\
\dim \HHH^1\left(X,\OO_X\left(p\right)\right) + \dim\HHH^{0}\left(X,\OO_X\left(p+d\right) \right)-  \h_{t-p}^{1,n-1}\left(X\right)  & m=1\\
\dim \HHH^m\left(X,\OO_X\left(p\right)\right) + \dim\HHH^{m-1}\left(X,\OO_X\left(p+d\right) \right)-  \h_{t-p}^{\frac{m}{2},n-\frac{m}{2}}\left(X\right)&\begin{matrix}1<m<2n \\\text{even}\end{matrix}\\
\dim \HHH^m\left(X,\OO_X\left(p\right)\right) + \dim\HHH^{m-1}\left(X,\OO_X\left(p+d\right) \right)-  \h_{t-p}^{\frac{m+1}{2},n-\frac{m+1}{2}}\left(X\right)&\begin{matrix}1<m<2n \\ \text{odd}\end{matrix}\\
\dim \HHH^{2n} \left(X,\OO_X\left(p\right)\right) + \dim\HHH^{2n-1}\left(X,\OO_X\left(p+d\right) \right)-  \h_{t-p-d}^{n-1,1}\left(X\right)& m=2n \\
\dim \HHH^{2n}\left(X,\OO_X\left(p+d\right)\right) & m=2n+1\\
0 &\text{else.}
\end{cases}
$$
\begin{proof}
For dimension reasons we immediately get $\dim\HHH^m\left(\PP^{n+1},f_*\OO_X\left(p\right)\right)=0$ for $m<0$ respectively $2n+1<m$.
Now for the computations:

\begin{itemize}
\item[$m=0$:] We compute
\begin{align*}
\dim\HHH^0\left(\PP^{n+1},f_*\OO_X\left(p\right)\right)&=\h_{t-p}^{0,n} \left(X\right)  &\text{Proposition~\ref{Hochschild dimension direct image explicit}}\\
&=\dim\HHH^0\left(X,\OO_X\left(p\right)\right).&\text{ Corollary~\ref{dimension Hochschild cohomology as upper middle and bottom part}}
\end{align*}

\item[$m=1$:] We get by Proposition~\ref{Hochschild dimension direct image explicit} and Corollary~\ref{dimension Hochschild cohomology as upper middle and bottom part}
\begin{align*}
\dim\HHH^{1}\left(\PP^{n+1},f_*\OO_X\left(p\right)\right)&=\h_{t-p}^{1,n}\left(X\right)+\h_{t-p-d}^{0,n}\left(X\right)-\h_{t-p}^{1,n-1}  \left(X\right) \\
&=\dim\HHH^1\left(X,\OO_X\left(p\right) \right) + \dim\HHH^0\left(X,\OO_X\left(p+d\right)\right)-\h_{t-p}^{1,n-1}\left(X\right). 
\end{align*}

\item[$1<m<2n$:] In this case we have by Corollary~\ref{dimension Hochschild cohomology as upper middle and bottom part} and Proposition~\ref{Hochschild dimension direct image explicit}
\begin{align*}
\dim\HHH^m\left(\PP^{n+1},f_*\OO_X\left(p\right)\right)&=\h^{0,m-n}_{t-p}\left(X\right) +\h^{0,m-n-1}_{t-p-d}\left(X\right) +\h^{n,m}_{t-p}\left(X\right) +\h^{n,m-1}_{t-p-d}\left(X\right) \\
&=\h^{0,m-n}_{t-p}\left(X\right)+\h^{n,m}_{t-p}\left(X\right)  +\h^{0,m-n-1}_{t-p-d}\left(X\right) +\h^{n,m-1}_{t-p-d}\left(X\right)  \\
&=\begin{cases}
\dim \HHH^m\left(X,\OO_X\left(p\right)\right) + \dim\HHH^{m-1}\left(X,\OO_X\left(p+d\right) \right)-  \h_{t-p}^{\frac{m}{2},n-\frac{m}{2} }\left(X\right)&m \text{ even}\\
\dim \HHH^m\left(X,\OO_X\left(p\right)\right) + \dim\HHH^{m-1}\left(X,\OO_X\left(p+d\right) \right)-  \h_{t-p}^{\frac{m}{2},n-\frac{m+1}{2}}\left(X\right) &m \text{ odd}.\\
\end{cases}
\end{align*}

\item[$m=2n$:] Here we get by Proposition~\ref{Hochschild dimension direct image explicit} and Corollary~\ref{dimension Hochschild cohomology as upper middle and bottom part}
\begin{align*}
\dim\HHH^{2n}\left(\PP^{n+1},f_*\OO_X\left(p\right)\right)&=\h_{t-p}^{n,0}\left(X\right)+\h_{t-p-d}^{n-1,0}\left(X\right)-\h_{t-p}^{n-1,1}  \left(X\right) \\
&=\dim\HHH^{2n}\left(X,\OO_X\left(p\right) \right) + \dim\HHH^{2n-1}\left(X,\OO_X\left(p+d\right)\right)-\h_{t-p-d}^{n-1,1}\left(X\right). 
\end{align*}

\item[$m=2n+1$:]We compute
\begin{align*}
\dim\HHH^{2n+1}\left(\PP^{n+1},f_*\OO_X\left(p\right)\right)&=\h_{t-p-d}^{n,0} \left(X\right) &\text{Proposition~\ref{Hochschild dimension direct image explicit}}\\
&=\dim\HHH^n\left(X,\OO_X\left(p+d\right)\right)&\text{ Corollary~\ref{dimension Hochschild cohomology as upper middle and bottom part}.}
\end{align*}
\end{itemize}
So we covered all cases and the statement holds.
\end{proof}
\end{proposition}

\begin{proposition}[{\cite[Proposition 9.5.1]{Rizzardo2019}}]\label{long exact sequence Hochschild cohomology}
Consider the embedding of a smooth $n$-dimensional degree $d$ hypersurface $X \xhookrightarrow{f} \PP^{n+1}$. Then we have a long exact sequence of the form:
\begin{equation*}
\cdots \to \HHH^{i-2}\left(X,\OO_{X}\left(p+d\right)\right)\to \HHH^{i}\left(X,\OO_{X}\left(p\right)\right)\xrightarrow{f_*} \HHH^{i}\left(\PP^{n+1},f_*\OO_X\left(p\right)\right)\to \cdots.
\end{equation*}
\end{proposition}

\begin{theorem}\label{kernels are precisely the middle part}
Let $f:X \hookrightarrow \PP^{n+1}$ be the embedding of a smooth degree $d$ hypersurface and set $t=d-n-2$. Then we have for all $p \in \Z$ such that $t-p \notin \left\{0,d\right\}$
$$\dim \ker\left(f_*:\HHH^m\left(X,\OO_X\left(p\right)\right)\to \HHH^m\left(\PP^{n+1},f_*\OO_X\left(p\right)\right)\right)=\begin{cases}
\h^{\frac{m}{2},n-\frac{m}{2} }_{t-p} \left(X\right)&\begin{matrix} 0< m< 2n  \text{ even}\end{matrix}\\
\h^{n-1,1}_{t-p-d}\left(X\right) & m=2n\\
0& \text{else.}
\end{cases}$$

\begin{proof} For dimension reasons we may assume $0\le m\le 2n$. In the diagrams for this prove we will denote $\OO_X$ by $\OO$ and $\PP^{n+1}$ by $\PP$ in order to avoid clumsy notation.

We will proceed by induction over $l$ with $2l=m$ using the long exact sequence from Proposition~\ref{long exact sequence Hochschild cohomology}:
$$\cdots \to \HHH^{m-2}\left(X,\OO_{X}\left(p+d\right)\right)\to \HHH^{m}\left(X,\OO_{X}\left(p\right)\right)\xrightarrow{f_*} \HHH^{m}\left(\PP^{n+1},f_*\OO_X\left(p\right)\right)\to \cdots .$$
 This way we can cover the odd case $2l-1$ and even case $2l$ in the induction step simultaneously:

We will start with $l=1$ as induction start and include the case of $m=0$ to cover the cases for $m=0,1,2$:

We compute all the dimensions in the long exact sequence in Proposition~\ref{long exact sequence Hochschild cohomology} using Proposition~\ref{Dimensions match up apart from middle part} and proceed by diagram chase. Consider the following diagram, where we denote the spaces on the left and their dimensions on the right. We will use the arrows on the right-hand side to indicate that the dimensions to the right of their tail are the dimensions to the left of their tip.

\begin{small}
$$\tikz[xscale=3.3,yscale=1.05,baseline]{
\node (ghost) at (6,0){$\;$};

\node (zero) at (0,7) {$0$};
\node (HH2l-2p) at (0,6) {$\HHH^{0}\left(X,\OO\left(p\right)\right)$};
\node (HH2l-2fp) at (0,5) {$\HHH^{0}\left(\PP,f_*\OO\left(p\right)\right)$};
\node (HH2l-3p+d) at (0,4) {$\HHH^{-1}\left(X,\OO\left(p+d\right)\right)$};
\node (HH2l-1p) at (0,3) {$\HHH^{1}\left(X,\OO\left(p\right)\right)$};
\node (HH2l-1fp) at (0,2) {$\HHH^{1}\left(\PP,f_*\OO\left(p\right)\right)$};
\node (HH2l-2p+d) at (0,1) {$\HHH^{0}\left(X,\OO\left(p+d\right)\right)$};
\node (HH2lp) at (0,0) {$\HHH^{2}\left(X,\OO\left(p\right)\right)$};
\begin{small}
\node (k2l-2p) at (2.4,6) {$\mathllap{0}$};
\node (k2l-2fp) at (2.4,5) {$\mathllap{\dim\HHH^{0}\left(X,\OO\left(p\right)\right)}$};
\node (k2l-3p+d) at (2.4,4) {$\mathllap{0}$};
\node (k2l-1p) at (2.4,3) {$\mathllap{0}$};
\node (k2l-1fp) at (2.4,2) {$\mathllap{\dim\HHH^{1}\left(X,\OO\left(p\right)\right)}$};
\node (k2l-2p+d) at (2.4,1) {$\mathllap{\dim\HHH^{0}\left(X,\OO\left(p+d\right)\right)-\h_{t-p}^{1,n-1}}$};
\node (k2lp) at (2.4,0) {$\mathllap{\h_{t-p}^{1,n-1}}$};
\node (i2l-3p) at (2.6,7) {$\mathrlap{0}$};
\node (i2l-2p) at (2.6,6) {$\mathrlap{\dim\HHH^{0}\left(X,\OO\left(p\right)\right)}$};
\node (i2l-2fp) at (2.6,5) {$\mathrlap{0}$};
\node (i2l-3p+d) at (2.6,4) {$\mathrlap{0}$};
\node (i2l-1p) at (2.6,3) {$\mathrlap{\dim \HHH^{1}\left(X,\OO\left(p\right)\right)}$};
\node (i2l-1fp) at (2.6,2) {$\mathrlap{\dim\HHH^{0}\left(X,\OO\left(p+d\right)\right)-\h_{t-p}^{1,n-1}}$};
\node (i2l-2p+d) at (2.6,1) {$\mathrlap{\h_{t-p}^{1,n-1}}$};
\node (i2lp) at (2.6,0) {$\mathrlap{\h_{t-p}^{1,n-1}+\dim \HHH^{2}\left(X,\OO\left(p\right)\right).}$};
\node (+) at (2.5,6) {$+$};
\node (+) at (2.5,5) {$+$};
\node (+) at (2.5,4) {$+$};
\node (+) at (2.5,3) {$+$};
\node (+) at (2.5,2) {$+$};
\node (+) at (2.5,1) {$+$};
\node (+) at (2.5,0) {$-$};
\node (ghost) at (4,0) {$\;$};
\end{small}
\draw[->]
(zero) edge (HH2l-2p)
(HH2l-2p)edge node[right]{\begin{tiny}$f_*$\end{tiny}} (HH2l-2fp)
(HH2l-2fp)edge (HH2l-3p+d) 
(HH2l-3p+d) edge (HH2l-1p)
(HH2l-1p) edge node[right]{\begin{tiny}$f_*$\end{tiny}} (HH2l-1fp)
(HH2l-1fp) edge (HH2l-2p+d)
(HH2l-2p+d) edge (HH2lp);
\draw[->,dashed]
(i2l-3p)edge (k2l-2p)
(i2l-2p)edge (k2l-2fp)
(i2l-2fp)edge (k2l-3p+d)
(i2l-1fp)edge (k2l-2p+d)
(i2l-3p+d)edge (k2l-1p)
(i2l-1p)edge (k2l-1fp)
(i2l-2p+d)edge (k2lp)
; 
}$$
\end{small}
By the above diagram chase we get that the image of the last arrow on the left has dimension $\h_{t-p}^{1,n-1}\left(X\right)$. So by the exactness of the sequence from Proposition~\ref{long exact sequence Hochschild cohomology} we get that this is also the dimension of the kernel of 
$$f_*:\HHH^2\left(X,\OO_X\left(p\right)\right)\to \HHH^2\left(\PP^n,f_*\OO_X\left(p\right)\right).$$
And so we get:
\begin{align*}
&\dim \ker\left(f_*:\HHH^0\left(X,\OO_X\left(p\right)\right)\to \HHH^0\left(\PP^n,f_*\OO_X\left(p\right)\right)\right)=0\\
&\dim \ker\left(f_*:\HHH^1\left(X,\OO_X\left(p\right)\right)\to \HHH^1\left(\PP^n,f_*\OO_X\left(p\right)\right)\right)=0\\
&\dim \ker\left(f_*:\HHH^2\left(X,\OO_X\left(p\right)\right)\to \HHH^2\left(\PP^n,f_*\OO_X\left(p\right)\right)\right)=\h_{t-p}^{1,n-1}\left(X\right)
\end{align*}
as expected.

For the induction step we will cover the cases $m=2l-1$ and $m=2l$ simultaneously. Assume that
$$\dim \ker\left(f_*:\HHH^{2l-2}\left(X,\OO_X\left(p\right)\right)\to \HHH^{2l-2}\left(\PP^n,f_*\OO_X\left(p\right)\right)\right)=\h^{n-l+1,l-1}_{t-p} \left(X\right).$$

We compute again the dimensions in the long exact sequence from Proposition~\ref{long exact sequence Hochschild cohomology} using our computations in Proposition~\ref{Dimensions match up apart from middle part}. We write the long exact sequence on the left and the dimensions on the right. We draw the arrows on the right hand side from left to right to indicate that the dimensions to the right of their tail are the dimensions to the left of their tip.
\begin{small}

$$\tikz[xscale=3.3,yscale=1.5,baseline]{
\node (ghost) at (6,0){$\;$};

\node (ker) at (0,7) {$\ker\left(f_*\right)^{2l-2} $};
\node (HH2l-2p) at (0,6) {$\HHH^{2l-2}\left(X,\OO\left(p\right)\right)$};
\node (HH2l-2fp) at (0,5) {$\HHH^{2l-2}\left(\PP,f_*\OO\left(p\right)\right)$};
\node (HH2l-3p+d) at (0,4) {$\HHH^{2l-3}\left(X,\OO\left(p+d\right)\right)$};
\node (HH2l-1p) at (0,3) {$\HHH^{2l-1}\left(X,\OO\left(p\right)\right)$};
\node (HH2l-1fp) at (0,2) {$\HHH^{2l-1}\left(\PP,f_*\OO\left(p\right)\right)$};
\node (HH2l-2p+d) at (0,1) {$\HHH^{2l-2}\left(X,\OO\left(p-d\right)\right)$};
\node (HH2lp) at (0,0) {$\HHH^{2l}\left(X,\OO\left(p\right)\right)$};
\begin{small}
\node (k2l-2p) at (2.4,6) {$\mathllap{\h^{l-1,n-l+1}_{t-p} }$};
\node (k2l-2fp) at (2.4,5) {$\mathllap{-\h_{t-p}^{l-1,n-l+1}+\dim\HHH^{2l-2}\left(X,\OO\left(p\right)\right)}$};
\node (k2l-3p+d) at (2.4,4) {$\mathllap{\dim\HHH^{2l-2}\left(X,\OO\left(p+d\right)\right)}$};
\node (k2l-1p) at (2.4,3) {$\mathllap{0}$};
\node (k2l-1fp) at (2.4,2) {$\mathllap{\dim\HHH^{2l-1}\left(X,\OO\left(p\right)\right)}$};
\node (k2l-2p+d) at (2.4,1) {$\mathllap{\dim\HHH^{2l-2}\left(X,\OO\left(p+d\right)\right)-\h_{t-p}^{l,n-l}}$};
\node (k2lp) at (2.4,0) {$\mathllap{\h_{t-p}^{l,n-l}}$};
\node (i2l-2p) at (2.6,6) {$\mathrlap{\h^{l-1,n-l+1}_{t-p}+\dim\HHH^{2l-2}\left(X,\OO\left(p\right)\right)}$};
\node (i2l-2fp) at (2.6,5) {$\mathrlap{\dim\HHH^{2l-2}\left(X,\OO\left(p+d\right)\right)}$};
\node (i2l-3p+d) at (2.6,4) {$\mathrlap{0}$};
\node (i2l-1p) at (2.6,3) {$\mathrlap{\dim \HHH^{2l-1}\left(X,\OO\left(p\right)\right)}$};
\node (i2l-1fp) at (2.6,2) {$\mathrlap{\dim\HHH^{2l-2}\left(X,\OO\left(p+d\right)\right)-\h^{l-1,n-l+1}_{t-p} }$};
\node (i2l-2p+d) at (2.6,1) {$\mathrlap{\h_{t-p}^{l,n-l}\left(X\right)}$};
\node (i2lp) at (2.6,0) {$\mathrlap{\h_{t-p}^{l,n-l}+\dim \HHH^{2l}\left(X,\OO\left(p\right)\right).}$};
\node (dimker) at (2.6,7) {$\mathrlap{\h^{l-1,n-l+1}_{t-p} }$};
\node (+) at (2.5,0) {$-$};
\node (+) at (2.5,6) {$-$};
\node (+) at (2.5,5) {$+$};
\node (+) at (2.5,4) {$+$};
\node (+) at (2.5,3) {$+$};
\node (+) at (2.5,2) {$+$};
\node (+) at (2.5,1) {$+$};
\node (ghost) at (4,0) {$\;$};
\end{small}
\draw[right hook ->]
(ker) edge (HH2l-2p);
\draw[->]
(HH2l-2p)edge node[right]{\begin{tiny}$f_*$\end{tiny}} (HH2l-2fp)
(HH2l-2fp)edge (HH2l-3p+d) 
(HH2l-3p+d) edge (HH2l-1p)
(HH2l-1p) edge node[right]{\begin{tiny}$f_*$\end{tiny}} (HH2l-1fp)
(HH2l-1fp) edge (HH2l-2p+d)
(HH2l-2p+d) edge (HH2lp);
\draw[->,dashed]
(dimker) edge (k2l-2p)
(i2l-2p)edge (k2l-2fp)
(i2l-2fp)edge (k2l-3p+d)
(i2l-1fp)edge (k2l-2p+d)
(i2l-3p+d)edge (k2l-1p)
(i2l-1p)edge (k2l-1fp)
(i2l-2p+d)edge (k2lp)
; 
}$$
\end{small}

By the exactness of the sequence this means that
\begin{align*}
&\dim \ker\left(f_*:\HHH^{2l-1}\left(X,\OO_X\left(p\right)\right)\to \HHH^{2l-1}\left(\PP^n,f_*\OO_X\left(p\right)\right)\right)=0\\
&\dim \ker\left(f_*:\HHH^{2l}\left(X,\OO_X\left(p\right)\right)\to \HHH^{2l}\left(\PP^n,f_*\OO_X\left(p\right)\right)\right)=\h^{n-l,l}_{t-p} \left(X\right).
\end{align*}

Now finally for the case of $l=n$:

By the above induction we have
$$\dim \ker\left(f_*:\HHH^{2n-2}\left(X,\OO_X\left(p\right)\right)\to \HHH^{2n-2}\left(\PP^n,f_*\OO_X\left(p\right)\right)\right)=\h^{1,n-1}_{t-p} \left(X\right).$$

We apply again diagram chase along long exact sequence from Proposition~\ref{long exact sequence Hochschild cohomology} using the computations in Proposition~\ref{Dimensions match up apart from middle part}. We continue to write the long exact sequence on the left and the dimensions on the right. The diagonal arrows on the right again symbolize that the dimensions to the right of their tail are the dimensions of the kernel to the left of their tip:
\begin{small}

$$\tikz[xscale=3.45,yscale=1.4,baseline]{
\node (ghost) at (6,0){$\;$};

\node (ker) at (0,7) {$\ker\left(f_*\right)^{2n-2} $};
\node (HH2l-2p) at (0,6) {$\HHH^{2n-2}\left(X,\OO\left(p\right)\right)$};
\node (HH2l-2fp) at (0,5) {$\HHH^{2n-2}\left(\PP,f_*\OO\left(p\right)\right)$};
\node (HH2l-3p+d) at (0,4) {$\HHH^{2n-3}\left(X,\OO\left(p+d\right)\right)$};
\node (HH2l-1p) at (0,3) {$\HHH^{2n-1}\left(X,\OO\left(p\right)\right)$};
\node (HH2l-1fp) at (0,2) {$\HHH^{2n-1}\left(\PP,f_*\OO\left(p\right)\right)$};
\node (HH2l-2p+d) at (0,1) {$\HHH^{2n-2}\left(X,\OO\left(p-d\right)\right)$};
\node (HH2lp) at (0,0) {$\HHH^{2n}\left(X,\OO\left(p\right)\right)$};
\begin{small}
\node (k2l-2p) at (2.3,6) {$\mathllap{\h^{n-1,1}_{t-p} }$};
\node (k2l-2fp) at (2.3,5) {$\mathllap{-\h_{t-p}^{n-1,1}+\dim\HHH^{2n-2}\left(X,\OO\left(p\right)\right)}$};
\node (k2l-3p+d) at (2.3,4) {$\mathllap{\dim\HHH^{2n-2}\left(X,\OO\left(p+d\right)\right)}$};
\node (k2l-1p) at (2.3,3) {$\mathllap{0}$};
\node (k2l-1fp) at (2.3,2) {$\mathllap{\dim\HHH^{2n-1}\left(X,\OO\left(p\right)\right)}$};
\node (k2l-2p+d) at (2.3,1) {$\mathllap{\dim\HHH^{2n-2}\left(X,\OO\left(p+d\right)\right)-\h_{t-p-d}^{n-1,1}}$};
\node (k2lp) at (2.3,0) {$\mathllap{\h_{t-p-d}^{n-1,1}}$};
\node (i2l-2p) at (2.5,6) {$\mathrlap{\h^{n-1,1}_{t-p} +\dim\HHH^{2n-2}\left(X,\OO\left(p\right)\right)}$};
\node (i2l-2fp) at (2.5,5) {$\mathrlap{\dim\HHH^{2n-3}\left(X,\OO\left(p+d\right)\right)}$};
\node (i2l-3p+d) at (2.5,4) {$\mathrlap{0}$};
\node (i2l-1p) at (2.5,3) {$\mathrlap{\dim \HHH^{2n-1}\left(X,\OO\left(p\right)\right)}$};
\node (i2l-1fp) at (2.5,2) {$\mathrlap{\dim\HHH^{2n-2}\left(X,\OO\left(p+d\right)\right)-\h^{n-1,1}_{t-p} }$};
\node (i2l-2p+d) at (2.5,1) {$\mathrlap{\h_{t-p-d}^{n-1,1}}$};
\node (i2lp) at (2.5,0) {$\mathrlap{\dim \HHH^{2n}\left(X,\OO\left(p\right)\right)-\h_{t-p-d}^{n-1,1}.}$};
\node (dimker) at (2.5,7) {$\mathrlap{\h^{n-1,1}_{t-p}}$};
\node (+) at (2.4,0) {$+$};
\node (+) at (2.4,6) {$-$};
\node (+) at (2.4,5) {$+$};
\node (+) at (2.4,4) {$+$};
\node (+) at (2.4,3) {$+$};
\node (+) at (2.4,2) {$+$};
\node (+) at (2.4,1) {$+$};
\node (ghost) at (4,0) {$\;$};
\end{small}
\draw[right hook ->]
(ker) edge (HH2l-2p);
\draw[->]
(HH2l-2p)edge node[right]{\begin{tiny}$f_*$\end{tiny}} (HH2l-2fp)
(HH2l-2fp)edge (HH2l-3p+d) 
(HH2l-3p+d) edge (HH2l-1p)
(HH2l-1p) edge node[right]{\begin{tiny}$f_*$\end{tiny}} (HH2l-1fp)
(HH2l-1fp) edge (HH2l-2p+d)
(HH2l-2p+d) edge (HH2lp);
\draw[->,dashed]
(dimker) edge (k2l-2p)
(i2l-2p)edge (k2l-2fp)
(i2l-2fp)edge (k2l-3p+d)
(i2l-1fp)edge (k2l-2p+d)
(i2l-3p+d)edge (k2l-1p)
(i2l-1p)edge (k2l-1fp)
(i2l-2p+d)edge (k2lp)
; 
}$$
\end{small}
This diagram gives us 

\begin{align*}
&\dim \ker\left(f_*:\HHH^{2n-1}\left(X,\OO_X\left(p\right)\right)\to \HHH^{2n-1}\left(\PP^n,f_*\OO_X\left(p\right)\right)\right)=0\\
&\dim \ker\left(f_*:\HHH^{2n}\left(X,\OO_X\left(p\right)\right)\to \HHH^{2n}\left(\PP^n,f_*\OO_X\left(p\right)\right)\right)=\h^{n-1,1}_{t-p-d} \left(X\right).
\end{align*}

So we covered the case for $m=2n$ and are done as for $m<0$ and $m>0$ the source space is trivial.
\end{proof}
\end{theorem}

We now finally state the following in order to guarantee the existence of non-trivial kernels of pushforwards of Hochschild cohomology.

\begin{proposition}\label{proposition guaranteed non-trivial kernel}
Let $f:X\hookrightarrow \PP^{2k}$ be an embedding of a smooth odd dimensional degree $d>1$ hypersurface of dimension $n=2k-1$ for $k>2$ and let $p=-kd-d$. Then we have
$$ \ker\left(\HHH^{n+3}\left(X,\OO_X\left(p\right)\right)\to \HHH^{n+3}\left(\PP^{n+1},f_*\OO_X\left(p\right)\right)\right)\cong \kk.$$
\begin{proof}
By Theorem~\ref{kernels are precisely the middle part} we have 
$$\dim \ker\left(\HHH^{n+3}\left(X,\OO_X\left(p\right)\right)\to \HHH^{n+3}\left(\PP^{n+1},f_*\OO_X\left(p\right)\right)\right)=\h^{k+1,k-2}_{t-p}\left(X\right) $$
with $t=d-n-2=d-2k-1$.

So it suffices to compute that $\h^{k+1,k-2}_{t-p}\left(X\right)=1$, with $t-p=kd+2d -2k -1$. By $\left(\ref{formula middle part}\right)$ this is
\begin{align*}
\h^{k+1,k-2}_{t-p}\left(X\right)&= \sum_{\mu=0}^{2k-1+2}\left(-1\right)^\mu \binom{2k-1+2}{\mu}\binom{-kd-2d+2k+1+\left(k+1\right)d-\left(\mu-1\right)\left(d-1\right)}{2k-1+1}\\
&= \sum_{\mu=0}^{2k+1}\left(-1\right)^\mu \binom{2k+1}{\mu}\binom{-kd-2d+2k+1+kd+d-\mu d + d +\mu-1}{2k}\\
&= \sum_{\mu=0}^{2k+1}\left(-1\right)^\mu \binom{2k+1}{\mu}\binom{2k-\mu d +\mu}{2k}\\
&=  \binom{2k+1}{0}\binom{2k }{2k}\\
&=1.
\end{align*}
Here we used that for $\mu>1$ we have $2k-\mu d+\mu <2k$ as $d>1$, which means that the terms $\binom{2k+1}{\mu}\binom{2k+\mu d +\mu}{2k}$ vanish for $\mu\geq 1$.

So we get
$$\dim \ker\left(\HHH^{n+3}\left(X,\OO_X\left(p\right)\right)\to \HHH^{n+3}\left(\PP^{n+1},f_*\OO_X\left(p\right)\right)\right)=1$$
as claimed.
\end{proof}
\end{proposition}

\subsection{Examples}
We collect a few examples of twisted Hodge diamonds that were computed using the Sage package by Pieter Belmans and Piet Glas \cite{TwistedHodgeDiamondsSage}. 

The first two examples illustrate the general shape as given in Lemma~\ref{upper, middle and bottom part} and the third will be an explicit example of Proposition~\ref{proposition guaranteed non-trivial kernel}.
\begin{example}
Let $f:X\hookrightarrow \PP^{6}$ be a smooth degree $7$ hypersurface then the $8$-twisted Hodge diamond is

$$\begin{matrix}
\;&\;&\;&\;&\;&0&\;&\;&\;&\;&\;\\
\;&\;&\;&\;&0&\;&0&\;&\;&\;&\;\\
\;&\;&\;&0&\;&0&\;&0&\;&\;&\;\\
\;&\;&0&\;&0&\;&0&\;&0&\;&\;\\
\;&0&\;&0&\;&0&\;&0&\;&0&\;\\
\mathllap{299}6&\;&\color{purple}\mathclap{20993}\color{purple}&\;&\color{purple}\mathclap{15267}\color{blue}&\;&\color{purple}\mathclap{917}\color{purple}&\;&0&\;&0\\
\;&\mathllap{157}5&\;&0&\;&0&\;&0&\;&0&\;\\
\;&\;&\mathllap{577}5&\;&0&\;&0&\;&0&\;&\;\\
\;&\;&\;&\mathllap{1039}5&\;&0&\;&0&\;&\;&\;\\
\;&\;&\;&\;&\mathllap{900}2&\;&0&\;&\;&\;&\;\\
\;&\;&\;&\;&\;&\mathllap{299}6.&\;&\;&\;&\;&\;\\
\end{matrix}$$

And so we have by Theorem~\ref{kernels are precisely the middle part}, since $t-8=-8$,
\begin{align*}
\dim \ker\left(f_*:\HHH^{4}\left(X,\OO_X\left(-8\right)\right)\to \HHH^{4}\left(X,f_*\OO_X\left(-8\right)\right)\right)&=\color{purple}917\color{black}\\
\dim \ker\left(f_*:\HHH^{6}\left(X,\OO_X\left(-8\right)\right)\to \HHH^{6}\left(X,f_*\OO_X\left(-8\right)\right)\right)&=\color{purple}15267\color{black}
\\\dim \ker\left(f_*:\HHH^{8}\left(X,\OO_X\left(-8\right)\right)\to \HHH^{8}\left(X,f_*\OO_X\left(-8\right)\right)\right)&=\color{purple}20993\color{black}.\\
\end{align*}
\end{example}

\begin{example}
Let $f:X\hookrightarrow \PP^{8}$ be smooth degree $5$ hypersurface then the $-7$-twisted Hodge diamond is

$$\begin{matrix}
\;&\;&\;&\;&\;&\;&\;&6\mathrlap{390}&\;&\;&\;\\
\;&\;&\;&\;&\;&\;&0&\;&2\mathrlap{0511}&\;&\;\\
\;&\;&\;&\;&\;&0&\;&0&\;&2\mathrlap{5704}&\;\\
\;&\;&\;&\;&0&\;&0&\;&0&\;&1\mathrlap{6840}\\
\;&\;&\;&0&\;&0&\;&0&\;&0&\;&4\mathrlap{950}&\;&\;\\
\;&\;&0&\;&0&\;&0&\;&0&\;&0&\;&7\mathrlap{20}&\;\\
\;&0&\;&0&\;&0&\;&0&\;&0&\;&0&\;&3\mathrlap{6}&\;\\
0&\;&0&\;&0&\;&\color{purple}\mathclap{486}&\;&\color{purple}\mathclap{13051}&\;&\color{purple}\mathclap{30276}&\;&\color{purple}\mathclap{8451}&\;&1\mathrlap{65}&\\
\;&0&\;&0&\;&0&\;&0&\;&0&\;&0&\;&0\\
\;&\;&0&\;&0&\;&0&\;&0&\;&0&\;&0&\;&\;\\
\;&\;&\;&0&\;&0&\;&0&\;&0&\;&0&\;\\
\;&\;&\;&\;&0&\;&0&\;&0&\;&0\\
\;&\;&\;&\;&\;&0&\;&0&\;&0&\;\\
\;&\;&\;&\;&\;&\;&0&\;&0&\;\\
\;&\;&\;&\;&\;&\;&\;&0.&\;\\
\end{matrix}$$

And so we have by Theorem~\ref{kernels are precisely the middle part}, since $t+7=3$,
\begin{align*}
\dim \ker\left(f_*:\HHH^{2}\left(X,\OO_X\left(3\right)\right)\to \HHH^{2}\left(X,f_*\OO_X\left(3\right)\right)\right)&=\color{purple}8451\\
\dim \ker\left(f_*:\HHH^{4}\left(X,\OO_X\left(3\right)\right)\to \HHH^{4}\left(X,f_*\OO_X\left(3\right)\right)\right)&=\color{purple}15267
\\\dim \ker\left(f_*:\HHH^{6}\left(X,\OO_X\left(3\right)\right)\to \HHH^{6}\left(X,f_*\OO_X\left(3\right)\right)\right)&=\color{purple}13051\\
\dim \ker\left(f_*:\HHH^{8}\left(X,\OO_X\left(3\right)\right)\to \HHH^{8}\left(X,f_*\OO_X\left(3\right)\right)\right)&=\color{purple}486.\\
\end{align*}
\end{example}

The next example illustrates a case of Proposition~\ref{proposition guaranteed non-trivial kernel}:
\begin{example}
Let $f:X \hookrightarrow \PP^{10}$ be a smooth degree $5$ hypersurface and consider $\OO_X\left(-30\right)$. Then we can compute, using Theorem~\ref{kernels are precisely the middle part}
$$\dim \ker\left(f_*:\HHH^{m}\left(X,\OO_X\left(-30\right)\right)\to \HHH^{m}\left(X,f_*\OO_X\left(-30\right)\right)\right).$$
To do this we need to compute the $t-p=24$ twisted Hodge-diamond
$$\begin{matrix}
\;&\;&\;&\;&\;&\;&\;&\;&\;&0&\;&\;&\;&\;&\;&\;&\;&\;&\;&\;&\;\\
\;&\;&\;&\;&\;&\;&\;&\;&0&\;&0&\;&\;&\;&\;&\;&\;&\;&\;&\;&\;\\
\;&\;&\;&\;&\;&\;&\;&0&\;&0&\;&0&\;&\;&\;&\;&\;&\;&\;&\;&\;\\
\;&\;&\;&\;&\;&\;&0&\;&0&\;&0&\;&0&\;&\;&\;&\;&\;&\;&\;&\;\\
\;&\;&\;&\;&\;&0&\;&0&\;&0&\;&0&\;&0&\;&\;&\;&\;&\;\\
\;&\;&\;&\;&0&\;&0&\;&0&\;&0&\;&0&\;&0&\;&\;&\;&\;\\
\;&\;&\;&0&\;&0&\;&0&\;&0&\;&0&\;&0&\;&0&\;&\;&\;&\;&\;&\;&\\
\;&\;&0&\;&0&\;&0&\;&0&\;&0&\;&0&\;&0&\;&0&\;&\;&\;&\;&\;&\;&\;\\
\;&0&\;&0&\;&0&\;&0&\;&0&\;&0&\;&0&\;&0&\;&0&\;&\;&\;\\
\mathllap{1197904}4&\;&\mathclap{100298}&\;&\mathclap{2882}&\;&\color{purple}1\color{black}&\;&0&\;&0&\;&0&\;&0&\;&0&\;&0&\\
\;&\mathllap{10743961}8&\;&0&\;&0&\;&0&\;&0&\;&0&\;&0&\;&0&\;&0&\;\\
\;&\;&\mathllap{52344510}9&\;&0&\;&0&\;&0&\;&0&\;&0&\;&0&\;&0&\;\\
\;&\;&\;&\mathllap{158002079}4&\;&0&\;&0&\;&0&\;&0&\;&0&\;&0&\;\\
\;&\;&\;&\;&\mathllap{314953851}3&\;&0&\;&0&\;&0&\;&0&\;&0&\;\\
\;&\;&\;&\;&\;&\mathllap{423631847}1&\;&0&\;&0&\;&0&\;&0&\;\\
\;&\;&\;&\;&\;&\;&\mathllap{381562624}3&\;&0&\;&0&\;&0&\;&\\
\;&\;&\;&\;&\;&\;&\;&\mathllap{220962657}3&\;&0&\;&0&\;&\;&\;\\
\;&\;&\;&\;&\;&\;&\;&\;&\mathllap{74465034}6&\;&0&\;&\;&\;&\;&\;\\
\;&\;&\;&\;&\;&\;&\;&\;&\;&\mathllap{11109813}0\mathrlap{.}&\;&\;&\;&\;&\;\\
\end{matrix}$$
And as expected by Proposition~\ref{proposition guaranteed non-trivial kernel} we get 
$$\dim \ker\left(f_*:\HHH^{12}\left(X,\OO_X\left(-30\right)\right)\to \HHH^{12}\left(X,f_*\OO_X\left(-30\right)\right)\right)=\color{purple}{1}\color{black}.$$
\end{example}

\section{Nontrivial kernel in Hochschild cohomology give non-Fourier-Mukai functors}\label{Section Nontrivial kernel in Hochschild cohomology give non-Fourier-Mukai functors}

In this section we follow the ideas from \cite{Rizzardo2019} to construct candidate non-Fourier-Mukai functors for hypersurfaces of arbitrary degree. We then verify that under assumptions on the characteristic morphisms and some concentrated $\Ext$-groups these indeed cannot be Fourier-Mukai. We finish the chapter by computing that these assumptions are satisfied when the source category is the derived category of an odd dimensional quadric, which gives concrete non-Fourier-Mukai functors between well behaved spaces in arbitrary high dimensions.

Since we follow the approach from \cite{Rizzardo2019} we will consider functors of a similar form:
\begin{equation}\label{Non-Fourier-Mukai functor}
\Psi_\eta:\D^b\left( X\right)\xrightarrow{L} \D^b_{w \Coh \left(X\right)}\left(\X_{\eta}^{dg}\right)\xrightarrow{\psi_{{\X,\eta},*}} \D^b_{w \Coh \left(X\right)}\left(\X_\eta\right)\xrightarrow{\widetilde{f}_*} \D^b\left( \PP^{n+1}\right),
\end{equation}
where $\X_\eta^{dg}$ denotes the dg-hull of $\X_\eta$ and $\psi_{{\X,\eta},*}$ is the induced comparison functor.

\subsection{Constructing candidate non-Fourier-Mukai functors}

We start by collecting a few results from \cite{Rizzardo2019}, which are central for our construction. We refer the interested reader to \cite{Rizzardo2019} for an in depth discussion.

In order to apply Definition~\ref{definition curly notation}, Lemma~\ref{lemma w and W}, Lemma~\ref{lemma equivariant w} and Lemma~\ref{lemma naturality curly notation} we assume that every quasi-projective scheme $X$ comes equipped with an open affine cover $X=\bigcup_{i=0}^m U_i$.

The following construction was originally introduced by W. Lowen and M. Van den Bergh in \cite{Lowen2011}. 

\begin{definition}\cite[Definition 4.2]{Raedschelders2019}\label{definition curly notation}
Let $X=\bigcup_{i=1}^m U_i$ be an open affine covering of a quasi-projective scheme. Consider for $I\subset \left\{1,...,m \right\}$ the sets $U_I:=\bigcap_{i \in I} U_i$ indexed by $I\in \I:= \mathrm{P}\left(X\right)\setminus \emptyset$. Then $\X$ is the category with objects $\I$ and morphisms:
$$\X\left(I,J\right):=\begin{cases}\OO_X\left(U_J\right)&I \subset J \\
0 &\text{else,}\end{cases}$$
where composition is induced by composing with the restriction morphism.
\end{definition}

Roughly $\X\modules$ acts as the category of presheaves associated to an affine covering. This means that it comes with the following useful properties:

\begin{lemma}\label{lemma w and W}\cite{Lowen2011}
Let $X$ be quasi-projective. Then there is a fully faithful embedding
\begin{equation*}
w: \D\left(\Qch X\right)\xrightarrow{\sim}\D_{w\Qch\left(X\right)}\left(\X\right)\hookrightarrow \D\left(\X\right)
\end{equation*}
and a fully faithful embedding 
\begin{equation*}
W:\Delta_*\D\left(\Qch X\right) \to \D\left(\X\otimes_\kk \X^\op\right),
\end{equation*}
where $\Delta_*\left(\Qch X\right)$ is the essential image of the direct image of the diagonal embedding $\Delta: X \to X\times X$. In particular we have for quasi-coherent $M$
$$\HHH^*\left(X,M\right)\cong \HHH^*\left(\X,WM\right).$$
\begin{proof} 
The construction of $w$ can be found at \cite[$\left(8.5\right)$]{Rizzardo2019}. The functor $W$ gets constructed in the following paragraph of \cite{Rizzardo2019}.
For the Hochschild cohomology comparison we can use that $W$ is a fully faithful embedding to get
$$W: \HHH^*\left(X,M\right):= \Ext^*_{X\times X}\left(\OO_\Delta,\Delta_* M\right) \xrightarrow{\sim}\Ext^*_{\X\otimes \X^{op}}\left(\X,WM\right)=:\HHH^*\left(\X,WM\right)$$
as desired.
\end{proof}
\end{lemma}

\begin{lemma}\label{lemma equivariant w}
Let $X$ be quasi-projective and let $\Gamma$ be a $\kk$-algebra. Then there is an embedding
$$w: \D\left(\Coh\left(X\right)_\Gamma\right)\hookrightarrow \D\left(\X\otimes \Gamma\right).$$
\begin{proof}
By \cite[\S 8.5]{Rizzardo2019} we have an embedding $w: \D\left(\Qch\left(X\right)_\Gamma\right)\hookrightarrow \D\left(\X\otimes \Gamma\right)$. There also is a canonical embedding
$$\D\left(\Coh\left(X\right)_\Gamma\right)\hookrightarrow \D\left(\Qch\left(X\right)_\Gamma\right).$$ In particular we get the desired embedding by composition.
\end{proof}
\end{lemma}

\begin{lemma}\cite[\S 8.7]{Rizzardo2019}\label{lemma naturality curly notation}
Let $f:X \to Y$ be a closed embedding of quasi-projective schemes. Then we have an induced functor
$$\mathfrak{f}:\Y \to \X,$$
such that the diagram
$$\tikz[heighttwo,xscale=3,yscale=2,baseline]{
\node (A) at (0,1) {$\D\left(X\right)$};
\node (B) at (0,0) {$\D\left(\X\right)$};
\node (C) at (1,1) {$\D\left(Y\right)$};
\node (D) at (1,0) {$\D\left(\Y\right)$};

\draw[->]
(A) edge node[left] {$w$} (B)
(C) edge node[right] {$w$}(D)
(A) edge node[above] {$f_*$}(C)
(B) edge node[below] {$\mathfrak{f}_*$}(D);
}$$
commutes.
\begin{proof}
The construction of $\mathfrak{f}$ is done in \cite[\S 8.7]{Rizzardo2019}. By \cite[Proposition 3.5]{HuybrechtsFm} we have an inclusion $\D\left(X\right)\hookrightarrow \D\left(\Qch\left(X\right)\right)$. So we can restrict the diagram from \cite[Lemma 8.7.1]{Rizzardo2019} to $\D\left(X\right)$.
\end{proof}
\end{lemma}

We use the following construction from \cite{Rizzardo2019} as the core of our candidate functors:
\begin{proposition}\label{proposition construction L}
Let $X$ be smooth projective of dimension $n$ and let $\eta \in \HHH^{\geq n+3}\left(X, M\right)$. Then there exists an exact functor 
$$\D^b\left(X\right)\xrightarrow{L}\D^b_{w\Coh\left(X\right)}\left(\X_\eta^{dg}\right)$$
such that $\underline{\Rhom}_{\X_\eta} \left(\X ,L \left(\_\right)\right) \cong w$.
\begin{proof}
First observe that by Lemma~\ref{lemma w and W} we have an isomorphism $$\HHH^*\left(X,M\right)\cong \HHH^*\left(\X,W M\right)$$ 
and so we may consider $\eta \in \HHH^{\geq n+3}\left(\X,WM\right)$.

By \cite[Lemma 10.1]{Rizzardo2019} and since $\Qch\left(X\right)$ has global dimension $n$ and $\HH^i\X_\eta$ vanishes in the right degrees we can apply \cite[Proposition~5.3.1]{Rizzardo2019} with $\A=w\Qch\left(X\right)$ and $\mathfrak{c}=\X_\eta$ to get a functor 
$$L': \D^b\left(\Qch\left(X\right)\right)\cong\D^b\left(w\Qch X\right)\to \D^b_{w\Qch\left(X\right)}\left(\X_\eta^{dg}\right).$$
Now we can use \cite[Proposition 3.5]{HuybrechtsFm} to turn this into a functor:
$$L: \D^b\left(X\right)\xhookrightarrow{\sim}\D^b_{\Coh\left(X\right)}\left(\Qch\left(X\right)\right)\xrightarrow{L'}\D^b\left(\X^{dg}_\eta\right)$$
with the desired property.

Finally, by \cite[Corollary 10.4]{Rizzardo2019} we know that the essential image of this functor is contained in $\D^b_{w\Coh{X}}\left(\X^{dg}_\eta\right)$.
\end{proof}
\end{proposition}

We will also use the following notation from \cite{Rizzardo2019} for $\widetilde{f}$.

\begin{proposition}\label{proposition widetilde f}\cite[Proposition 7.2.6]{Rizzardo2019}
Let $f:\PPP^{n+1} \to \X$ be a functor of $\kk$-linear categories and $\eta\in \HHH^{k}\left(\X,\M\right)$ such that $f_*\eta = 0$. Then there exists an $\Ainfty$-functor $\widetilde{f}$ making the diagram
$$\tikz[heighttwo,xscale=2,yscale=2,baseline]{
\node (A) at (0,0) {$\X$};
\node (B) at (2,0) {$\PPP^{n+1}$};
\node (C) at (1,1.3) {$\X_\eta$};
\draw[->]
(B) edge node[below] {$f$} (A)
(C) edge node[above left] {$\pi$}(A);
\draw[->, dashed]
(B) edge node[above right] {$\widetilde{f}$}(C);
}$$
commute. In particular we have
\begin{equation}\label{pi circ widetilde f is f}
\pi\circ \widetilde{f} = f.
\end{equation}
\end{proposition}

Now we construct a candidate functor $\Psi_\eta$ for $\eta \in \HHH^{\geq n+3}\left(X,\OO_X\left(p\right)\right)$:

\begin{construction}\label{construction Psi eta}

Let $X\hookrightarrow \PP^{n+1}$ be the embedding of a smooth $n$-dimensional scheme with $n\geq 3$ and let $$0 \neq\eta\in\ker\left(f_*:\HHH^{n+3} \left(X,\OO_X\left(p\right)\right)\to \HHH^{n+3}\left(\PP^{n+1},f_*\OO_X\left(p\right)\right)\right).$$
Then a functor of the form $\left(\ref{Non-Fourier-Mukai functor}\right)$ is constructed to be,
$$\Psi_\eta:\D^b\left(\Coh\left(X\right)\right)\xrightarrow{L} \D^b_{w \Coh \left(X\right)}\left(\X_{\eta}^{dg}\right)\xrightarrow{\psi_{{\X_\eta},*}} \D^b_{w \Coh \left(X\right)}\left(\X_\eta\right)\xrightarrow{\widetilde{f}_*}\D^b_{w\Coh\left(\PP^{n+1}\right)}\left(\PPP^{n+1}\right)\cong \D^b\left(\Coh \left(\PP^{n+1}\right)\right),$$
where we have the functor $L$ by Proposition~\ref{proposition construction L}, $\psi_{\X_\eta,*}$ is the functor constructed in \cite[§ D.1]{Rizzardo2019} and $\widetilde{f}_*$ exists by Proposition~\ref{proposition widetilde f}.
\end{construction}

\begin{corollary}
Let $f:X \to \PP^{n+1}$ be the embedding of a degree $d$ hypersurface and let $m>n+2$ then we have a $\h^{\frac{m}{2},n-\frac{m}{2}}_{p}\left(X\right)$-dimensional space of choices to construct a candidate functor 
$$\Psi_\eta : \D^b\left(X\right)\to \D^b\left(\PP^{n+1}\right).$$
\begin{proof}
In order for Construction~\ref{construction Psi eta} to work we need 
$$0\neq\eta\in\ker\left(f_*:\HHH^{m} \left(X,\OO_X\left(p\right)\right)\to \HHH^{m}\left(\PP^{n+1},f_*\OO_X\left(p\right)\right)\right).$$
By Theorem~\ref{kernels are precisely the middle part} $\ker\left(f_*:\HHH^{m} \left(X,\OO_X\left(p\right)\right)\to \HHH^{m}\left(\PP^{n+1},f_*\OO_X\left(p\right)\right)\right)$ has dimension $\h^{\frac{m}{2},n-\frac{m}{2}}_{p}\left(X\right)$ which finishes the claim.
\end{proof}
\end{corollary}

Now we can state our main Theorem, which we will prove throughout \S~\ref{Section Verifying that the functors are not just candidates}.

\begin{theorem}\label{Main Theorem}
Let $f:X\hookrightarrow \PP^{n+1}$ be an embedding of a smooth degree $d$ hypersurface of dimension $n\geq 3$ and let $$ 0\neq\eta\in\ker\left(f_*:\HHH^{n+3} \left(X,\OO_X\left(p\right)\right)\to \HHH^{n+3}\left(\PP^{n+1},f_*\OO_X\left(p\right)\right)\right)$$ such that there exists a $\kk$-algebra $\Gamma$ and $G \in \D^b\left(\Coh\left(X\right)_\Gamma\right)$ with 
\begin{align*}
&c_{G,\Gamma}\left(\eta\right)\neq 0&\\
&\Ext^i_{X}\left(G\left(-p\right),T\right)=0 &\mathllap{\quad\text {for } i\neq n}\\
 &\Ext^{n-1}_{X}\left(G,G\left(p+d\right)\right)\cong \Ext_X^{n-2}\mathrlap{\left(G,G\left(p+d\right)\right)\cong 0.}\quad&
\end{align*}
Then we have that the functor 
$$\Psi_\eta:\D^b\left(\Coh \left(X\right)\right)\to\D^b\left(\Coh \left(\PP^{n+1}\right)\right)$$
is well-defined and not a Fourier-Mukai functor.
\end{theorem}

\begin{remark}
By the same proof as \cite{Rizzardo2019}[Proposition B.2.1] the functors $\psi_\eta$ do not admit a lift to the spectral level in the case of $\kk=\mathbb{Q}$.
\end{remark}

\subsection{Proving Theorem~\ref{Main Theorem}}\label{Section Verifying that the functors are not just candidates}

We fix for the rest of this section an embedding of a smooth degree $d$ hypersurface $f:X\hookrightarrow \PP^{n+1}$, a non-vanishing Hochschild cohomology class $\eta\in \HHH^{n+3}\left(X,\OO\left(p\right)\right)$, such that $f_*\eta = 0$, $\Gamma$ a $\kk$-algebra and $G\in \D\left(\Coh\left(X\right)_\Gamma\right)$ such that
\begin{align}
& c_{G,\Gamma}\left(\eta\right)\neq 0 \tag{I}\label{assumption characteristic morphism}\\
&\Ext^i_{X}\left(G\left(-p\right),T\right)=0 &\mathllap{\text {for } i\neq n}\tag{II}\label{assumption Ext concentrated in one degree}\\
&\Ext^{n-1}_{X}\left(G,G\left(p+d\right)\right)\cong \Ext^{n-2}\mathrlap{\left(G,G\left(p+d\right)\right)\cong 0.}\quad & \tag{III}\label{assumption Ext vanishes in n-1 n-2}
\end{align} 

Observe first that by Construction~\ref{construction Psi eta} $\Psi_\eta$ is well-defined and even unique up to a choice of $\widetilde{f}$. So we may focus for the rest of this section on verifying that $\Psi_\eta$ cannot be Fourier-Mukai.

We follow mostly the ideas from \cite{Rizzardo2019}.

We start by recalling the following Lemma~\ref{Lemma obstruction lift of objects} from \cite{Rizzardo2019} in order to have obstructions against lifts of $\widetilde{\G}$ to $\D\left(\X_\eta \otimes \Gamma\right)$. These obstructions and their naturality will be later used in order to conclude that $\Psi_\eta$ cannot be Fourier-Mukai.

\begin{lemma}{\cite[Lemma~7.3.1]{Rizzardo2019}}\label{Lemma obstruction lift of objects}
\begin{enumerate}
\item Let $\X$ be a dg-category, $\Gamma$ a $\kk$-algebra and $\G\in \D\left(\X\right)_\Gamma$. Then there is a sequence of obstructions 
$$o_{i+2}\left(\G\right)\in \HHH^{i+2}\left(\Gamma,\Ext_\X^{-i}\left(\G,\G\right)\right)$$
for $i\geq 1$ such that $\G$ lifts to an object in $\D\left(\X\otimes_\kk \Gamma\right)$ if and only if all obstructions vanish. More precisely $o_{i+1}\left(\G\right)$ is only defined if $o_3\left(\G\right),...,o_i\left(\G\right)$ vanish and it depends on choices.

\item if $f:\Y\to \X$ is a dg-functor and $f_*:\D\left(\X\right)\to\D\left(\Y\right)$ is the corresponding change of rings functor, then after having made choices for $\G$ we may make corresponding choices for $f_*\left(\G\right) $ in such a way that
$$f_*\left(o_{i+2}\left(\G\right) \right)=o_{i+2}\left(f_*\left(\G\right)\right)$$
\end{enumerate}
\end{lemma}

We now use the assumptions on $G$ to prove that the negative part of $\Ext_{\X_\eta}^*\left(LG,LG\right)$ is concentrated in degree $-1$ which allows us to control which $\Ainfty$-obstruction does not vanish. This obstruction we will then push forward to prove that $\Psi_\eta$ cannot be Fourier-Mukai.
In order to avoid clumsy notation we start by setting

\begin{equation}\label{definition curly G and widetilde curly G}
\G:=wG\in \X\modules \text{ and } \widetilde{\G}:=L\left(G\right).
\end{equation}

\begin{remark}
We have by \cite[§ D.1]{Rizzardo2019} an equivalence $\psi_{\X_\eta}:\X_\eta^{dg}\xrightarrow{\sim}\X_\eta$ and by Definition~\ref{definition deformed A-infinity category} a canonical functor $\pi: \X_\eta \to \X$. So we will denote the functor
$$\psi^{-1}_{\X_\eta,*} \circ \pi_*: \D\left(\X\right) \to \D_\infty\left(\X_\eta\right) \to \D\left(\X_\eta^{dg}\right)$$
simply by $\pi_*$ and 
$$\psi_{\X_\eta,*} \circ \pi^*: \D\left(\X_\eta^{dg}\right) \to \D_\infty\left(\X_\eta\right) \to \D\left(\X\right)$$
by $\pi^*$ to avoid clumsy and confusing notation. 
\end{remark}

\begin{definition}\label{morphism negative Ext}
Consider the distinguished triangle in $\D\left(\X^{dg}_\eta\right)$ \cite[Lemma 10.3]{Rizzardo2019}: 
\begin{equation}\label{triangle negative Ext}
\G\xrightarrow{\alpha}\widetilde{\G}\xrightarrow{\beta}\Sigma^{-n-1}\G\otimes w\OO_X\left(-p\right)\xrightarrow{\gamma}\Sigma\G,
\end{equation}
where $\G$ is considered as an $\X^{dg}_\eta$-module via $ \pi_*:\D^b\left(\X\right) \to \D\left(\X^{dg}_\eta \right) $.
Then define the morphism $\varphi$ by:
\begin{align*}
\varphi: \Ext_X^{n+1+i}\left(G \left(-p\right),G\right) &\to \Ext_{\X^{dg}_\eta}^{i}\left(\widetilde{\G},\widetilde{\G}\right)\\
\left( g: \Sigma^{-n-1-i}G\left(-p\right) \to G\right) &\mapsto \alpha \circ \pi_*\left(w\left(g\right)\right) \circ \Sigma^{-i}\beta: \left(\Sigma^{-i}\widetilde{\G} \to \widetilde{\G} \right).
\end{align*}
\end{definition}

\begin{lemma}\label{negative Ext computation}
For $i<0$ the morphism
$$\varphi: \Ext_X^{n+1+i}\left(G\left(-p\right),G\right)\cong \Ext_{\X^{dg}_\eta}^{i}\left(\widetilde{\G},\widetilde{\G}\right)$$
is an isomorphism.
\begin{proof}
We will check that for $i<0$ the morphisms involved in the definition of 
\begin{align*}
\varphi: \Ext_X^{n+1+i}\left(G \left(-p\right),G\right) &\to \Ext_{\X^{dg}_\eta}^{i}\left(\widetilde{\G},\widetilde{\G}\right)\\
\left( g: \Sigma^{-n-1-i}G\left(-p\right) \to G\right) &\mapsto \left( \alpha \circ \pi_*\right)\left(w\left(g\right)\right) \circ \Sigma^{-i}\beta: \left(\Sigma^{-i}\widetilde{\G} \to \widetilde{\G} \right)
\end{align*}
 are isomorphisms.
\begin{itemize}
\item[$w$:] By Lemma~\ref{lemma w and W} $w: \D^b\left( X\right)\to \D^b\left(\X\right)$ is a fully faithful embedding, in particular, using $\G=wG$ $\left(\ref{definition curly G and widetilde curly G}\right)$ we have
$$w:\Ext^{n+1+i}_X\left(G,G\right)\xrightarrow{\sim}\Ext^{n+1+i}_{\X}\left(\G,\G\right).$$ 
\item[$\alpha\circ \pi_*\left(\_\right)$:]
  We have by \cite[Corollary 5.3.2]{Rizzardo2019} an adjunction:
  $$\Rhom_{\X_\eta^{dg}}\left(\G\otimes w \OO_X\left(-p\right),\widetilde{\G}\right)\cong\Rhom_\X\left(\G\otimes w \OO_X \left(-p\right),\G\right).$$
This isomorphism can be computed explicitly to be:
   $$\alpha\circ \pi_*:\Rhom_{\X_\eta^{dg}}\left(\G\otimes w\OO_X\left(-p\right),\widetilde{\G}\right)\cong\Rhom_\X\left(\G \otimes w \OO_X \left(-p\right),\G\right),$$
   see \cite[$\left(11.6\right)$]{Rizzardo2019}.
\item[$\_\circ\beta$:] Consider the distinguished triangle $\left(\ref{triangle negative Ext}\right)$ in $\D\left(\X^{dg}_\eta\right)$:
$$\G \xrightarrow{\alpha} \widetilde{\G}\xrightarrow{\beta} \Sigma^{-n-1}\G\otimes w\OO_X\left(-p\right).$$
Apply $\Rhom_{\X^{dg}_\eta}\left(\_,\widetilde{\G}\right)$ to get the distinguished triangle:
\begin{small}$$ \Rhom_{\X^{dg}_\eta}\left(\Sigma^{-n-1}\G\otimes w\OO_X\left(-p\right),\widetilde{\G}\right)\xrightarrow{\_\circ\beta}\Rhom_{\X^{dg}_\eta}\left(\widetilde{\G},\widetilde{\G} \right)\to \Rhom_{\X^{dg}_\eta}\left(\G,\widetilde{\G} \right).$$
\end{small}
Now we may use \cite[Corollary 5.3.2]{Rizzardo2019} and Proposition~\ref{definition curly notation}, 
$$\Rhom_{\X_\eta^{dg}}\left(\G,\widetilde{\G}\right)\cong\Rhom_{\X_\eta}\left(\G,\widetilde{\G}\right)\cong \Rhom_\X\left(\G,\G\right) \cong \Rhom_X\left(G,G\right),$$
  to get
\begin{small}
$$ \Rhom_{\X^{dg}_\eta}\left(\Sigma^{-n-1}\G\otimes w\OO_X\left(-p\right),\widetilde{\G}\right)\xrightarrow{\_\circ\beta}\Rhom_{\X^{dg}_\eta}\left(\widetilde{\G},\widetilde{\G} \right)\to \Rhom_{X}\left(G,G \right).$$
\end{small}
Applying $\HH^i$ turns this into the long exact sequence:
$$ \cdots\to \Ext^{i-1}_{X}\left(G,G \right)\to \Ext^{n+1+i}_{\X^{dg}_\eta}\left(\G\otimes w\OO_X\left(-p\right),\widetilde{\G}\right)\xrightarrow{\_\circ\beta}\Ext^i_{\X^{dg}_\eta}\left(\widetilde{\G},\widetilde{\G} \right)\to  \cdots.$$
And as $G$ is a sheaf on $X$, specializing to $i<0$ yields the long exact sequence:
$$ \cdots\to 0\to \Ext^{n+1+i}_{\X^{dg}_\eta}\left(\G\otimes w\OO_X\left(-p\right),\widetilde{\G}\right)\xrightarrow{\_\circ\beta}\Ext^i_{\X^{dg}_\eta}\left(\widetilde{\G},\widetilde{\G} \right)\to 0 \to \cdots.$$
In particular 
$$\_\circ \beta:  \Ext^{n+1+i}_{\X^{dg}_\eta}\left(\G\otimes w\OO_X\left(-p\right),\widetilde{\G}\right)\xrightarrow{\sim}\Ext^i_{\X^{dg}_\eta}\left(\widetilde{\G},\widetilde{\G} \right) $$
is an isomorphism for $i<0$.
\end{itemize}
So altogether we get that 
\begin{align*}
\varphi: \Ext_X^{n+1+i}\left(G \left(-p\right),G\right) &\xrightarrow{\sim} \Ext_{\X^{dg}_\eta}^{i}\left(\widetilde{\G},\widetilde{\G}\right)\\
\left( g: \Sigma^{-n-1-i}G\left(-p\right) \to G\right) &\mapsto \left(\alpha \circ \pi_*\left(w\left(g\right)\right) \circ \Sigma^{-i}\beta: \Sigma^{-i}\widetilde{\G} \to \widetilde{\G} \right)
\end{align*}
is indeed an isomorphism for $i<0$ as it is a composition of isomorphisms.
\end{proof}
\end{lemma}

\begin{corollary}\label{negative Ext concentrated}
Let $p<-n-1$ and $i>1$. Then $\Ext_{\X_\eta}^{-i}\left(\widetilde{\G},\widetilde{\G}\right)=0$.
\begin{proof}
By \eqref{assumption Ext concentrated in one degree} we have that $\Ext_{X}^{*}\left(G\left(-p\right),G\right)$ is concentrated in degree $n$ and so we have by Lemma~\ref{negative Ext computation}
$$\Ext^{-i}_{\X^{dg}_\eta}\left(\widetilde{\G},\widetilde{\G}\right)\cong \Ext_{X}^{n+1-i}\left(G\left(-p\right),G\right)\cong 0 $$
for $i>1$.

And since we have a quasi-equivalence $\X_\eta\cong \X_\eta^{dg}$ we get 
$$\Ext^{-i}_{\X_\eta}\left(\widetilde{\G},\widetilde{\G}\right)\cong \Ext_X^{n+1-i}\left(G\left(-p\right),G\right)$$
as claimed.
\end{proof}
\end{corollary}

\begin{definition}[{\cite[Lemma~11.4]{Rizzardo2019}}]\label{definition algebraic equivariant characteristic morphism}
Let $\X$ be a $\kk$-linear category, $\Gamma$ a $\kk$-algebra and let $\M$ be a $\kk$-central $\X$-bimodule. Then we have for a $\Gamma$-equivariant $\X$-module $\G$ ,i.e. $\G\in \X \modules_\Gamma$, the (algebraic) $\Gamma$-equivariant characteristic morphism
\begin{align*}
c_{\G,\Gamma}: \HHH^*\left(\X,\M\right)=\Ext^*_{\X\otimes \X^\op}\left(\X,\M\right)&\to \Ext^{*}_{\X\otimes \Gamma}\left(\G,\G\otimes \M\right)\\
\eta &\mapsto \G \otimes_\X \eta
\end{align*}
Observe that this morphism factors naturally as $$c_{\G,\Gamma}:\HHH^*\left(\X,\M\right)\xrightarrow{\eta\mapsto \eta\cup 1}\HHH^*\left(\X\otimes \Gamma, \M \otimes\Gamma\right)\xrightarrow{c_\G}\Ext_{\X\otimes \Gamma}^*\left(\G,\G\otimes \M\right),$$
where $c_\G: \HHH^*\left(\X\otimes \Gamma,\M\otimes \Gamma\right) \to \Ext_{\X\otimes \Gamma}^*\left(\G,\G\otimes \M\right)$ is the (algebraic) characteristic morphism for $\G\in \D\left(\X\otimes \Gamma\right)$, see Proposition~\ref{proposition obstruction colift}.
\end{definition}

\begin{lemma}\label{Lemma Commutative Diagram Gamma-equivariant characteristic morphisms}
There is a commutative diagram:
\begin{equation*}
\tikz[heighttwo,xscale=8,yscale=3,baseline]{
\node (A) at (0,1) {$\HHH^{n+3}\left(X,\OO_X\left(p\right)\right)$};
\node (B) at (0,0) {$\HHH^{n+3}\left(\X,w\OO_X\left(p\right)\right)$};
\node (C) at (1,1) {$\Ext_{\Coh\left(X\right)_\Gamma}^{n+3}\left(G, G\left(p\right)\right)$};
\node (D) at (1,0) {$\Ext_{\X\otimes \Gamma}^{n+3}\left(\G, \G\otimes w \OO_X\left(p\right)\right)$,};

\draw[->]
(A) edge node[left] {\rotatebox{90}{$\sim$}} (B)
(C) edge node[right] {\rotatebox{90}{$\sim$}}(D)
(A) edge node[above] {$c_{G,\Gamma}$}(C)
(B) edge node[above] {$c_{\G,\Gamma}$}(D);
}
\end{equation*}
where $c_{G,\Gamma}$ is the (geometric) equivariant characteristic morphism discussed in \S~\ref{Section Equivariant sheaves and the characteristic morphism} and $c_{\G,\Gamma}$ is the (algebraic) characteristic morphism from Definition~\ref{definition algebraic equivariant characteristic morphism}.
\begin{proof}
By \cite[$\left(8.13\right)$]{Rizzardo2019} we have the commutative diagram
$$\tikz[heighttwo,xscale=8,yscale=3,baseline]{
\node (A) at (0,1) {$\Delta_*\D\left(X\right)$};
\node (B) at (0,0) {$\D\left(\X\otimes \X^\op \right)$};
\node (C) at (1,1) {$\D\left(\Coh\left(X\right)_\Gamma\right)$};
\node (D) at (1,0) {$\D\left(\X\otimes\Gamma\right),$};

\draw[->]
(A) edge node[right] {$W$} node[left] {\rotatebox{90}{$\sim$}} (B)
(C) edge node[left] {$w$} node[right] {\rotatebox{90}{$\sim$}}(D)
(A) edge node[above] {$\pi_{1,*}\left( \_\otimes \pi_2^* G \right)$}(C)
(B) edge node[below] {$ \G \otimes_\X \_$}(D);
}
$$
where we denote by $\Delta_*\D\left(X\right) \subset \D \left(X\times X\right)$ the essential image of the direct image along the diagonal embedding $\Delta: X \to X\times X$.

Considering the induced diagram on morphism spaces for $$\HHH^{n+3}\left(X,\OO\left(p\right)\right)=\Ext^{n+3}_{X\times X}\left(\OO_\Delta,\OO_\Delta\left(p\right)\right)=\Ext^{n+3}_{\Delta_*\D\left( X\right)}\left(\OO_\Delta,\OO_\Delta\left(p\right)\right)$$ gives that the diagram
$$
\tikz[heighttwo,xscale=8,yscale=3,baseline]{
\node (A) at (0,1) {$\HHH^{n+3}\left(X,\OO_X\left(p\right)\right)$};
\node (B) at (0,0) {$\HHH^{n+3}\left(\X,w\OO_X\left(p\right)\right)$};
\node (C) at (1,1) {$\Ext_{\Coh\left(X\right)_\Gamma}^{n+3}\left(G, G\left(p\right)\right)$};
\node (D) at (1,0) {$\Ext_{\X\otimes \Gamma}^{n+3}\left(\G, \G\otimes w \OO_X\left(p\right)\right)$};

\draw[->]
(A) edge node[left] {\rotatebox{90}{$\sim$}} (B)
(C) edge node[right] {\rotatebox{90}{$\sim$}}(D)
(A) edge node[above] {$c_{G,\Gamma}$}(C)
(B) edge node[above] {$c_{\G,\Gamma}$}(D);
}$$
commutes.
\end{proof}
\end{lemma}

Since $G \in \D^b\left(\Coh\left(X\right)_\Gamma\right)$ we get a $\Gamma$-action on $\G$ and $\widetilde{\G}$ via the functors $w$ and $L$, i.e. $\widetilde{\G}\in \D\left(\X_\eta^{dg} \right)_\Gamma$. So Lemma~\ref{Lemma obstruction lift of objects} gives well-defined obstructions against $\widetilde{\G}\in \D_\infty\left(\X_\eta\right)_\Gamma\cong\D\left(\X^{dg}_\eta\right)_\Gamma$ admitting a lift to an $\Ainfty$-module in $\D_\infty\left(\X_\eta\otimes \Gamma\right)$:
$$o_i\left(\widetilde{\G}\right)\in \HHH^{i}\left(\Gamma, \Ext_{\X_{\eta}}^{2-i}\left(\widetilde{\G},\widetilde{\G}\right)\right) \quad \text{for } i>2. $$

\begin{remark}
The next Lemma will use the obstruction obtained from the equivariant characteristic morphism \eqref{assumption characteristic morphism} in order to conclude that the first $\Ainfty$-obstruction against an equivariant lift of $\widetilde{\G}$ cannot vanish. We do this by observing that a colift of $\G$ to $\X_\eta$ would also give an equivariant lift of $\widetilde{\G}$. The control of $o_3 \left(\widetilde{\G}\right)$ is necessary as we want to push forward the obstruction from $\X_\eta$ to $\PPP^{n+1}$ which cannot be done with the obstruction arising by the characteristic morphism.
\end{remark}

\begin{lemma}\label{obstruction}
We have:
$$0 \neq o_3\left(\widetilde{\G}\right)\in \HHH^{3}\left(\Gamma,\Ext^{-1}_{\X_{\eta}}\left(\widetilde{\G},\widetilde{\G}\right)\right).$$
\begin{proof}
Assume $o_3\left(\widetilde{\G}\right)$ vanishes. Then by Corollary~\ref{negative Ext concentrated} $\Ext^{-i}_{\X_{\eta}}\left(\widetilde{\G},\widetilde{\G}\right)=0$ for $i>1$ and so: $$o_i\left(\widetilde{\G}\right)\in \HHH^{i}\left(\Gamma,\Ext^{2-i}_{\X_\eta}\left(\widetilde{\G},\widetilde{\G}\right)\right)=0$$
for all $i>2$.

So $\widetilde{\G}$ would admit a lift, i.e. an object $$\widehat{\G} \in \D\left(\X_\eta^\text{dg}\otimes\Gamma\right)\cong \D_\infty\left(\X_{\eta}\otimes \Gamma\right)$$
with $\widehat{\G}\cong \widetilde{\G}$ in $\D\left(\X_\eta\right)_\Gamma$.

Consider the triangle $\left(\ref{triangle negative Ext}\right)$ in $\D_\infty\left(\X_\eta\right)\cong\D\left(\X_\eta^{dg}\right)$
$$\G \to \widetilde{\G}\cong \widehat{\G}\to \Sigma^{n+1}\G\otimes w\OO_X\left(-p\right)\to \Sigma \G,$$
where we use the shorthand $\G$ for $\pi_*\G$. 
This gives:
$$\HH^*\left(\widehat{G}\right)\cong \G\oplus  \Sigma^{n+1}\G\otimes w\OO_X\left(-p\right). $$
By the construction of the triangle \eqref{triangle negative Ext} in \cite[§ 10]{Rizzardo2019}, the above isomorphism is compatible with the $\X_\eta$-action. So by Definition~\ref{definition colift} $\widehat{\G}$ is a colift of $\G\in \D\left(\X\otimes_\kk \Gamma\right)$ to $\D_\infty\left(\left(\X_\eta\otimes_\kk \Gamma\right)_{\eta\cup 1}\right)$.

By  Proposition~\ref{proposition obstruction colift}, the obstruction against such a colift is the image of $\eta\cup 1$ under the characteristic morphism
$$\HHH^{n+3}\left(\X\otimes \Gamma,w\OO_X\left(p\right)\otimes \Gamma\right) \to \Ext_{\X\otimes\Gamma}^{n+3}\left(\G,\G \otimes w \OO_X\left(p\right)\right).$$
However, this obstruction cannot vanish. As if we consider the equivariant characteristic morphism 
$$c_{\G,\Gamma}:\HHH^{n+3}\left(\X,w\OO_X\left(p\right)\right)\xrightarrow{\mu \mapsto \mu \cup 1}\HHH^{n+3}\left(\X\otimes \Gamma,w\OO_X\left(p\right)\otimes \Gamma\right) \xrightarrow{c_\G} \Ext_{\X\otimes\Gamma}^{n+3}\left(\G,\G \otimes w \OO_X\left(p\right)\right),$$
we have the commutative diagram from Lemma~\ref{Lemma Commutative Diagram Gamma-equivariant characteristic morphisms}:
$$\tikz[heighttwo,xscale=8,yscale=2,baseline]{
\node (A) at (0,1) {$\HHH^{n+3}\left(X,\OO_X\left(p\right)\right)$};
\node (B) at (0,0) {$\HHH^{n+3}\left(\X,w\OO_X\left(p\right)\right)$};
\node (C) at (1,1) {$\Ext_{\Coh\left(X\right)_\Gamma}^{n+3}\left(G, G\left(p\right)\right)$};
\node (D) at (1,0) {$\Ext_{\X\otimes \Gamma}^{n+3}\left(\G, \G\otimes w \OO_X\left(p\right)\right).$};

\draw[->]
(A) edge node[left] {\rotatebox{90}{$\sim$}} (B)
(C) edge node[right] {\rotatebox{90}{$\sim$}}(D)
(A) edge node[above] {$c_{G,\Gamma}$}(C)
(B) edge node[above] {$c_{\G,\Gamma}$}(D);
}$$
By assumption \eqref{assumption characteristic morphism} we have that $c_{G,\Gamma}\left(\eta\right)\neq 0$. So $c_\G\left(\eta\cup 1\right)\neq 0$, which means that such a colift of $\widetilde{\G}$ to $\left(\X \otimes \Gamma\right)_{\eta\cup 1}$ cannot exist. Now by the discussion above this means that a lift of $\widetilde{\G}$ to $\D_\infty\left(\X_\eta\otimes \Gamma\right)$ cannot exist and so $o_3\left(\widetilde{\G}\right)$ cannot be zero.
\end{proof}
\end{lemma}

\begin{lemma}\label{square embedding}
There is a commutative diagram
\begin{equation}\label{commutative square embedding}
\tikz[heighttwo,xscale=8,yscale=2,baseline]{
\node (A) at (0,1) {$\Ext^{n}_{X}\left(G\left(-p\right), G\right)$};
\node (B) at (0,0) {$\Ext^{n}_{\PP^{n+1}}\left(f_*\left(G\left(-p\right)\right), f_*\left(G\right)\right)$};
\node (C) at (1,1) {$\Ext^{-1}_{\X_{\eta}}\left(\widetilde{\G}, \widetilde{\G}\right)$};
\node (D) at (1,0) {$\Ext^{-1}_{\PPP^{n+1}}\left(\widetilde{f}_*\left(\widetilde{\G}\right), \widetilde{f}_*\left(\widetilde{\G}\right)\right)$,};

\draw[->]
(A) edge node[left] {$f_*$} (B)
(C) edge node[right] {$\widetilde{f}_* \circ \psi_{\X_{\eta},*}$}(D)
(A) edge node[above] {$\varphi$}(C)
(B) edge node[below] {$\widetilde{f}_* \varphi$}  (D);
}
\end{equation} 
where the lower morphism is given by
\begin{align*}
\widetilde{f}_*\varphi :\Ext^{n}_{\PP^{n+1}}\left(f_*G\left(-p\right), f_*G\right))&\to  \Ext^{-1}_{\PPP^{n+1}}\left(\widetilde{f}_*\widetilde{\G}, \widetilde{f}_*\widetilde{\G}\right)\\
g&\mapsto \widetilde{f}_*\alpha \circ w\left(g\right)\circ \widetilde{f}_*\beta.
\end{align*}
\begin{proof}
Recall that by Definition~\ref{morphism negative Ext} the morphism $\varphi$ is given by
\begin{align*}
\varphi: \Ext_X^{n+1+i}\left(G \left(-p\right),G\right) &\to \Ext^{-1}_{\X^{dg}_\eta}\left(\widetilde{\G}, \widetilde{\G}\right)\\
\left( g: \Sigma^{-n-1-i}G\left(-p\right) \to G\right) &\mapsto \left(\alpha \circ \pi_*\left(wg\right) \circ \beta: \Sigma^i\widetilde{\G} \to \widetilde{\G} \right),
\end{align*}
where $\alpha$ and $\beta$ are the first and second morphisms in the distinguished triangle $\left(\ref{triangle negative Ext}\right)$ in $\D\left(\X^{dg}_\eta\right)$
$$
\G\xrightarrow{\alpha}\widetilde{\G}\xrightarrow{\beta}\Sigma^{-n-1}\G\otimes w\OO_X\left(-p\right)\xrightarrow{\gamma}\Sigma\G.
$$

Applying the exact functor $\widetilde{f}_*\circ \psi_{\X_\eta,*}$ gives the distinguished triangle in $\D\left(\PPP^{n+1} \right)$
$$\widetilde{f}_*\G \xrightarrow{\widetilde{f}_*\alpha} \widetilde{f}_*\widetilde{\G}\xrightarrow{\widetilde{f}_*\beta} \Sigma^{-n-1}\widetilde{f}_*\G\otimes w\OO_X\left(-p\right)\xrightarrow{\widetilde{f}_*\gamma} \Sigma \widetilde{f}_*\G,$$
 which is a shorthand for
$$\widetilde{f}_*\pi_*\G \xrightarrow{\widetilde{f}_*\alpha} \widetilde{f}_*\widetilde{\G}\xrightarrow{\widetilde{f}_*\beta} \Sigma^{-n-1}\widetilde{f}_*\pi_*\G\otimes w\OO_X\left(-p\right)\xrightarrow{\widetilde{f}_*\gamma} \Sigma \widetilde{f}_*\pi_*\G.$$

So we may use $\pi\circ \widetilde{f}=f$  to get
\begin{equation}\label{triangle widetilde}
{f}_*\G \xrightarrow{\widetilde{f}_*\alpha} \widetilde{f}_*\widetilde{\G}\xrightarrow{\widetilde{f}_*\beta} \Sigma^{-n-1}{f}_*\G\otimes w\OO_X\left(-p\right)\xrightarrow{\widetilde{f}_*\gamma} \Sigma {f}_*\G.
\end{equation}

In particular
\begin{align*}
\widetilde{f}_*\varphi :\Ext^{n}_{\PP^{n+1}}\left(f_*\left(w\OO_X\left(p\right)\otimes \G\right), f_*\left(\G\right)\right)&\to  \Ext^{-1}_{\PPP^{n+1}}\left(\widetilde{f}_*\left(\widetilde{\G}\right), \widetilde{f}_*\left(\widetilde{\G}\right)\right)\\
g&\mapsto \widetilde{f}_*\alpha \circ wg\circ \widetilde{f}_*\beta
\end{align*}
is well-defined.

Now we compute
\begin{align*}
\widetilde{f}_*\circ \varphi \left(g\right)&= \widetilde{f_*}\left(\alpha\circ \pi_*\left(w g\right)\circ \beta\right)&\text{Definition of }\varphi\\
&=\widetilde{f}_* \alpha \circ \left(\widetilde{f}_* \circ \pi_*\right) \left(wg\right) \circ \widetilde{f}_* \beta &\widetilde{f}_*\text{ is a functor}\\
&=\widetilde{f}_* \alpha \circ\left( f_*\circ w\right)\left( g\right) \circ \widetilde{f}_* \beta & \widetilde{f}_* \circ \pi_* \left(g\right) = f_* g\\ 
&=\widetilde{f}_* \alpha \circ \left(w \circ f_*\right) \left(g\right) \circ \widetilde{f}_* \beta & \text{\cite[Lemma~8.7.1]{Rizzardo2019} }\\ 
&=\widetilde{f}_*\varphi\left( f_*g\right)  & \text{Definition of }\widetilde{f}_* \varphi
\end{align*}
and the diagram indeed commutes.
\end{proof}
\end{lemma}

\begin{corollary}\label{iso in degree -1}
The right map in the diagram \eqref{commutative square embedding}
$$\widetilde{f}_*\circ \psi_{\eta,*}:\Ext^{-1}_{\X^{dg}_\eta}\left(\widetilde{\G},\widetilde{\G}\right)\xrightarrow{\sim}\Ext^{-1}_{\PPP^{n+1}}\left(\widetilde{f}_*\left(\widetilde{\G}\right), \widetilde{f}_*\left(\widetilde{\G}\right)\right)$$
 is an isomorphism.
\begin{proof}

Since $G,G\left(-p\right)$ are coherent sheaves on $X$ and $f_*$ is exact we have
\begin{align*}
\Ext^{1}_{\PPP^{n+1}}\left(f_*\left(\Sigma^{-n-1}\G\otimes w\OO_X\left(-p\right)\right), f_*\left(\G\right)\right)
&\cong \Ext^{1}_{\PP^{n+1}}\left(f_*\left(\Sigma^{-n-1}G\otimes \OO_X\left(-p\right)\right), f_*\left(G\right)\right)&\text{Lemma~\ref{lemma w and W} }\\
&\cong \Ext^{n+2}_{\PP^{n+1}}\left( f_*\left(G\left(-p\right)\right), f_*\left(G\right)\right)&\mathllap{\Ext^i\left(\Sigma^{-j}\_,\_\right)\cong\Ext^{i+j}\left(\_,\_\right)}\\
&=0. &\dim\PP^{n+1}=n+1&
\end{align*}

So in the distinguished triangle (\ref{triangle widetilde})
$${f}_*\G \xrightarrow{\widetilde{f}_*\alpha} \widetilde{f}_*\widetilde{\G}\xrightarrow{\widetilde{f}_*\beta} \Sigma^{-n-1}{f}_*\G\otimes w\OO_X\left(-p\right)\xrightarrow{\widetilde{f}_*\gamma} \Sigma {f}_*\G$$
$\widetilde{f}_* \gamma$ vanishes, and we have:
$$\widetilde{f}_*\left(\widetilde{\G}\right)\cong f_*\left(\G\right)\oplus f_*\left(\Sigma^{-n-1}w\OO_X\left(-p\right)\otimes_{\X} \G\right)$$
via the splitting morphisms $\widetilde{f}_*\alpha$ and $\widetilde{f}_*\beta$. 

This means that both, the top morphism, by Lemma~\ref{negative Ext computation}, and the lower morphism, by splitting, in \eqref{commutative square embedding} are isomorphsims. So by Lemma~\ref{square embedding} it suffices to prove that 
$$f_*:\Ext^{n}_{X}\left(\OO_X\left(-p\right)\otimes G, G\right) \to \Ext^{n}_{\PP^{n+1}}\left(f_*\left(\OO_X\left(-p\right)\otimes G\right), f_*\left(G\right)\right)$$
is an isomorphism.

As tensoring with $\OO_X\left(p\right)$ is an autoequivalence this is equivalent to 
 $$f_*:\Ext^{n}_{X}\left(G, G\left(p\right)\right) \to \Ext^{n}_{\PP^{n+1}}\left(f_* G, f_*G\left(p\right)\right)$$
 being an isomorphism. Consider the long exact sequence associated to a divisor \cite[(9.13)]{Rizzardo2019}:
\begin{small}
$$\cdots \to \Ext^{n-2}_{X}\left( G, G\left(p+d\right)\right)\to \Ext^{n}_{X}\left(G, G\left(p\right)\right)\xrightarrow{f_*} \Ext^{n}_{\PP^{n+1}}\left(f_*\left(G\right), f_*G\left(p\right)\right)\to\cdots .$$
\end{small}

By assumption \eqref{assumption Ext vanishes in n-1 n-2} we have
$$\Ext^{n-2}_{X}\left( G, G\left(p+d\right)\right)\cong 0 \quad \text{and} \quad \Ext^{n-1}_{X}\left(G, G\left(p+d\right)\right)\cong 0, $$
 so the long exact sequence has the shape
$$\cdots \to 0 \to \Ext^{n}_{X}\left(G, G\left(p\right)\right)\xrightarrow{f_*} \Ext^{n}_{\PP^{n+1}}\left(f_*\left(G\right), f_*G\left(p\right)\right)\to 0\to\cdots .$$

By exactness that immediately gives that
 $$f_*:\Ext^{n}_{X}\left(G\left(-p\right), G\right)\xrightarrow{\sim} \Ext^{n}_{\PP^{n+1}}\left(f_*G\left(-p\right), f_*G\right)$$
is an isomorphism,  which finishes the proof.
\end{proof}
\end{corollary}

\begin{lemma}\label{obstruction on P^n does not vanish}
The obstruction $o_3\left(\Psi_\eta\left(G\right)\right)\in \HHH^3\left(\Gamma,\Ext^{-1}_{\PP^{n+1}}\left(\Psi\left(G\right),\Psi\left(G\right)\right)\right)$ against lifting to $\D\left(\PPP^{n+1}\otimes \Gamma\right)$ from Lemma~\ref{Lemma obstruction lift of objects} does not vanish.
\begin{proof}
By part (2) of Lemma~\ref{Lemma obstruction lift of objects} we have 
$$o_3\left(\Psi\left(G\right)\right)=\left(\widetilde{f}_*\circ \psi_{\X_{\eta},*}\right)o_3\left(\widetilde{\G} \right)\in\HHH^3\left(\Gamma,\Ext^{-1}_{\PP^{n+1}}\left(\Psi\left(G\right),\Psi\left(G\right)\right)\right).$$
Furthermore, as $o_3$ is the first obstruction we do not need to keep track of any choices.
So we can use Corollary \ref{iso in degree -1} to get that $\widetilde{f}_*\circ \psi_{\X_{\eta},*}$ induces an isomorphism in degree $-1$ and by Lemma \ref{obstruction} we have $0\neq o_3\left(\widetilde{\G} \right)$. So altogether
\begin{equation*}
0 \neq \left(\widetilde{f}_*\circ \psi_{\X_{\eta},*}\right)o_3\left(\widetilde{\G}\right)=o_3\left(\Psi\left(G\right)\right)\in \HHH^3\left(\Gamma,\Ext^{-1}_{\PP^{n+1}}\left(\Psi\left(G\right),\Psi\left(G\right)\right)\right).\qedhere
\end{equation*}
\end{proof}
\end{lemma}

Now we can finally finish the proof of Theorem~\ref{Main Theorem}.

\begin{proof}

Assume $\Psi_\eta$ is Fourier-Mukai. Then by Corollary~\ref{Fourier-Mukai compatible with basechange} $\Psi_\eta $ admits  a lift $$\Psi_{\eta,\Gamma}: \D^b\left(\Coh\left(X\right)_\Gamma\right) \to \D^b\left(\Coh\left(\PP^{n+1}\right)_\Gamma\right).$$ This means that $\Psi_{\eta}\left(G \right)\in \D^b\left(\PP^{n+1}\right)_\Gamma$ has a lift to $\D^b\left(\Coh\left(\PP^{n+1} \right)_\Gamma\right)\hookrightarrow \D_\infty\left(\PPP^{n+1}\otimes_\kk \Gamma\right)$. Since we have by Lemma~\ref{obstruction on P^n does not vanish}
$$o_3\left(\Psi_\eta \left(G\right) \right)\neq 0 $$
 such a lift cannot exist. 

So $\Psi_\eta$ cannot be Fourier Mukai. 
\end{proof}

\subsection{Application: odd dimensional Quadrics}

We will show that the tilting bundle $G$ for an odd dimensional quadric hypersurface and its endomorphism algebra $\Gamma$ satisfy the assumptions of Theorem~\ref{Main Theorem}. For this we start by recalling that quadrics admit an exceptional sequence, which gives rise to a tilting bundle.

\begin{theorem}[{\cite[Corollary~3.2.8]{Boehning}}]\label{Theorem quadric has an exceptional sequence}
Let $Q\hookrightarrow \PP^{2k}$ be the embedding of a smooth quadric. Then $Q$ admits an exceptional sequence:
$$\left(S\left(-2k+1\right),\OO_Q\left(-2k+2\right),...,\OO_Q\left(-1\right),\OO_Q\right),$$
where $S$ denotes the spinor bundle.
\end{theorem}

In particular we may consider for the embedding of a smooth quadric $f: Q \hookrightarrow \PP^{2k}$ the tilting bundle:

$$G:= S\left(-2k+1\right)\oplus \bigoplus_{l=0}^{-2k+2}\OO_Q\left(-l\right) \text{ and } \Gamma:=\End\left(G\right).$$

Now we need to verify the assumptions on the concentration of $\Ext^*_Q\left(G(-p),G)\right)$ and $\Ext^*_Q\left(G,G\left(p+d\right)\right)$. We will use $p=-2k-2$ and 
$$0\neq \eta \in \ker\left(f_*:\HHH^{n+3}\left(X,\OO_Q\left(-2k-2\right)\right)\to \HHH^{n+3}\left(\PP^{2k},f_*\OO_Q\left(-2k-2\right)\right)\right)$$
 as we know by Proposition~\ref{proposition guaranteed non-trivial kernel} that
$$f_*:\HHH^{n+3}\left(Q,\OO_Q\left(-2k-2\right)\right)\to \HHH^{n+3}\left(\PP^{n+1},f_*\OO_Q\left(-2k-2\right)\right)$$
has one-dimensional kernel.

For the $\Ext$-calculations we will need the following statement which also holds for even quadrics. However, as in the even case we would need to track the different spinor bundles depending on the equivalence class of the dimension modulo four, we will restrict to the odd case for legibility.

\begin{lemma}\label{lemma selfExt spinor}
Let $Q\hookrightarrow \PP^{2k}$ be a smooth odd dimensional quadric and let $S$ be the spinor bundle. Then the following hold:
\begin{enumerate}
\item\label{lemma selfExt spinor case 1} We have for $i\notin \left\{0,1,n\right\}$
$$\Ext_Q^i\left(S,S\left(m\right)\right)\cong \Ext_Q^{i-1}\left(S,S\left(m+1\right)\right).$$
\item\label{lemma selfExt spinor case 2} If $m\le -1$ we have additionally

\begin{align*}
\Ext_Q^i\left(S,S\left(m\right)\right)&\cong \Ext_Q^{i-1}\left(S,S\left(m+1\right)\right)\\
\Ext_Q^i\left(S,S\left(m\right)\right)&\cong \Ext_Q^{i-1}\left(S,S\left(m+1\right)\right).
\end{align*}

\end{enumerate}
\begin{proof}
Consider the short exact sequence \cite[Theorem~2.8]{Ottaviani}
$$0\mapsto S \mapsto \OO_Q^{2^{k+1}}\to S\left(1\right)\to 0$$
which gives after applying $\Ext^i_Q\left(\_,S\left(m+1\right)\right)$ the long exact sequence 
$$
\tikz[heighttwo,xscale=5.3,yscale=3.5,baseline]{
\node (dotsin) at (0,1) {$\cdots$};
\node (Ei-1OS) at (1,1) {$\Ext^{i-1}_Q\left(\OO_Q^{2^{k+1}},S\left(m+1\right)\right)$};
\node (Ei-1S) at (2,1) {$\Ext^{i-1}_Q\left(S,S\left(m+1\right)\right)$};
\node (EiS) at (0,0) {$\Ext^{i}_Q\left(S\left(1\right),S\left(m+1\right)\right)$};
\node (EiOS) at (1,0) {$\Ext_Q^i\left(\OO_Q^{2^{k+1}},S\left(m+1\right)\right)$};
\node (dotsout) at (2,0) {$\cdots$.};

\draw[->]
(dotsin) edge (Ei-1OS)
(Ei-1OS) edge (Ei-1S)
(Ei-1S)edge[out=355,in=175,overlay] (EiS)
(EiS) edge (EiOS)
(EiOS) edge (dotsout);
}$$
In particular we have $$\Ext^i_Q\left(S,S\left(m\right)\right)\cong \Ext^i_Q\left(S\left(1\right),S\left(m+1\right)\right)\cong \Ext^{i+1}_Q\left(S,S\left(m-1\right)\right)$$ if we have for $j\in \left\{i,i-1\right\}$ 
$$\Ext^j\left(\OO_Q^{2^{k+1}},S\left(m+1\right)\right)\cong \bigoplus_{l=0}^{2^{k+1}}\Ext^j\left(\OO_Q,S\left(m+1\right)\right)\cong \bigoplus_{l=0}^{2^{k+1}}\HH^j\left(X,S\left(m+1\right)\right)=0.$$
By \cite[Theorem~2.3]{Ottaviani} we have $\HH^j\left(X, S\left(m+1\right)\right)=0$ for $j \notin\left\{0,n\right\}$ which implies \ref{lemma selfExt spinor case 1}. 

If $m\le-1$ we have $m+1\le 0$ and so we get by \cite[Theorem~2.3]{Ottaviani} $H^0\left(X,S\left(m+1\right)\right)=0$, which gives \ref{lemma selfExt spinor case 2}.
 \end{proof}
\end{lemma}

\begin{proposition}\label{proposition ext concentrated in one degree}
Let $i \neq 2k-1$. Then we have $$\Ext_Q^*\left(G\left(2k+2\right),G\right)=0$$
\begin{proof}
Since $G$ is a sheaf we may assume $0\le i \le 2k-2$ for dimension reasons.
By definition of $G$ and additivity of $\Ext$ we have
\begin{align*}
\Ext^i_Q\left(G(2k+2),G\right)&=\Ext^i_Q \left(\left(S\left(-2k+1\right)\oplus \bigoplus_{l=0}^{2k-2}\OO_Q\left(-l\right)\left(2k+2\right)\right),S\left(-2k+1\right)\oplus \bigoplus_{l=0}^{2k-2}\OO_Q\left(-l\right)\right)\\
&\cong \bigoplus_{h,l=0}^{2k-2}\Ext_Q^i\left(\OO_Q\left(2k+2 -l\right),\OO_Q\left(-h\right)\right)\\
&\;\oplus \bigoplus_{l=0}^{2k-2}\Ext_Q^i \left(\OO_Q\left(2k+2-l\right),S\left(-2k+1\right)\right)\\
&\;\oplus \bigoplus_{l=0}^{2k-2}\Ext_Q^i \left(S\left(2k+2-2k+1\right),\OO_Q\left(-l\right)\right)\\
&\;\oplus \Ext^i_Q\left(S\left(-2k+1+2k+2\right),S\left(-2k+1\right)\right)
\end{align*}
In particular we can compute these $\Ext$-groups one by one.

We start with $\Ext_Q^i\left(\OO_Q\left(2k+2-l\right),\OO_Q\left(-h\right)\right)$ for which we get
\begin{align*}
\Ext_Q^i\left(\OO_Q\left(2k+2 -l\right),\OO_Q\left(-h\right)\right)&\cong \Ext_Q^i\left(\OO_Q,\OO_Q\left(l-h-2k-2\right)\right)&\text{twisting on both sides}\\
&\cong \HH^i\left(Q,\OO_Q\left(l-h-2k-2\right)\right) &\mathllap{\Ext^i_Q\left(\OO_Q,\_\right)\cong \HH^i\left(Q,\_\right)}\\
&\cong 0. &l-h-2k-2<0
\end{align*}

Since we have $l-4k-1 \le 0$ we get
\begin{align*}
\Ext_Q^i \left(\OO_Q\left(2k+2-l\right),S\left(-2k+1\right)\right)&\cong \Ext_Q^i \left(\OO_Q,S\left(l-4k-1\right)\right)&\text{twisting on both sides}\\
&\cong\HH^i \left(Q,S\left(l-4k-1\right)\right) & \mathllap{\Ext^i_Q\left(\OO_Q,\_\right)\cong \HH^i\left(Q,\_\right)}\\
&\cong 0. &\mathllap{\text{\cite[Theorem~2.3]{Ottaviani}}}
\end{align*}

By \cite[Theorem~2.8]{Ottaviani} we have $S^\vee \cong S\left(1\right)$ which we may use to compute 
\begin{align*}
\Ext_Q^i \left(S\left(3\right),\OO_Q\left(-l\right)\right) &\cong \Ext_Q^i\left(S,\OO_Q\left(-3-l\right)\right)&\text{twisting on both sides}\\
&\cong \Ext_Q^i\left(\OO_Q,S^\vee\left(-3-l\right)\right)&\mathllap{\text{dualizing}}\\
&\cong \Ext_Q^i\left(\OO_Q,S\left(-2-l\right)\right)& \mathllap{S^\vee \cong S\left(1\right)}\\
&\cong \HH^i\left(Q,S\left(-2-l\right)\right)&\mathllap{\Ext^i_Q\left(\OO_Q,\_\right)\cong \HH^i\left(Q,\_\right)}\\
&\cong 0. &\mathllap{\text{ \cite[Theorem~2.3]{Ottaviani}}}
\end{align*}

We may use $i\le 2k-1$ to get
\begin{align*}
\Ext^i_Q\left(S\left(3\right),S\left(-2k+1\right)\right)&\cong \Ext^i_Q\left(S,S\left(-2k-2\right)\right) &\text{twisting on both sides}\\
&\cong \Ext^1_Q\left(S,S\left(-2k-3+i\right)\right)&\mathllap{\text{Lemma~\ref{lemma selfExt spinor}}}\\
&\cong \Ext^{-1}_Q\left(S,S\left(-2k-1+i\right)\right)&\mathllap{\text{Lemma~\ref{lemma selfExt spinor}}}\\
&\cong 0. &\mathllap{S\text{ is a sheaf}}
\end{align*}

So every direct summand vanishes, and in particular 
$$\Ext^i_Q\left(G(2k+2),G\right)=0 \text{ for } i \neq n$$
as desired.
\end{proof}
\end{proposition}

\begin{proposition}\label{proposition Ext concentrated in 2 degrees}
Let $i \notin \left\{0,2k-1\right\}$. Then we have
$$\Ext_Q^i\left(G,G\left(-2k\right)\right)\cong 0.$$
\begin{proof}
Since $Q$ has dimension $2k-1$ and $G$ is a sheaf we may assume $0<i<2k-1$.

By definition of $G$ and additivity of $\Ext^i_Q\left(\_,\_\right)$ we have:
\begin{align*}
\Ext_Q^i\left(G,G\left(-2k\right)\right)&=\Ext_Q^i \left(S\left(-2k+1\right) \oplus\bigoplus_{l=0}^{2k-2}\OO_Q\left(-l\right),S\left(-4k+1\right)\oplus\bigoplus_{h=0}^{2k-2}\OO_Q\left(-2k-h\right)\right)\\
&\cong \bigoplus_{l,h=0}^{2k-2}\Ext^i_Q\left(\OO_Q\left(-l\right),\OO_Q\left(\OO_Q\left(-2k-h\right)\right)\right)\\
&\; \oplus \bigoplus_{l=0}^{2k-2}\Ext^i_Q\left(\OO_Q\left(-l\right),S\left(-4k+1\right)\right)\\
&\; \oplus \bigoplus_{l=0}^{2k-2} \Ext^i_Q\left(S\left(-2k+1\right),\OO_Q\left(-2k-l\right)\right)\\
&\; \oplus \Ext^i_Q\left(S\left(-2k+1\right),S\left(-4k+1\right)\right).
\end{align*}
As above  we can compute the cases separately.

We start with $\Ext^i_Q\left(\OO_Q\left(-l\right),\OO_Q\left(-2k-h\right)\right)$. For this we get:
\begin{align*}
\Ext^i_Q\left(\OO_Q\left(-l\right),\OO_Q\left(-2k-h\right)\right)&\cong \Ext^i_Q\left(\OO_Q,\OO_Q\left(l-2k-h\right)\right)&\text{twisting on both sides}\\
&\cong \HH^i\left(Q,\OO_Q\left(1-2k-h\right)\right) &\mathllap{\Ext^i_Q\left(\OO_Q,\_\right)\cong \HH^i\left(Q,\_\right)}\\
&\cong 0. &\mathllap{i \notin \left\{0,2k-1\right\}}
\end{align*}
For $\Ext^i_Q\left(\OO_Q\left(-l\right),S\left(-4k+1\right)\right)$ we get
\begin{align*}
\Ext^i_Q\left(\OO_Q\left(-l\right),S\left(-4k+1\right)\right) &\cong \Ext^i\left(\OO_Q,S\left(l-4k+1\right)\right)&\text{twisting on both sides}\\
&\cong \HH^i\left(Q,S\left(l-4k+1\right)\right) &\mathllap{\Ext^i_Q\left(\OO_Q,\_\right)\cong \HH^i\left(Q,\_\right)}\\
&\cong 0. & \mathllap{ i\notin\left\{0,2k-1\right\} \; \text{\cite[Theorem~2.3]{Ottaviani}} }
\end{align*}
While for  $\Ext^i_Q\left(S\left(-2k+1\right),\OO_Q\left(-2k-l\right)\right)$ one can compute:
\begin{align*}
 \Ext^i_Q\left(S\left(-2k+1\right),\OO_Q\left(-2k-l\right)\right) &\cong  \Ext^i_Q\left(S,\OO_Q\left(-1-l\right)\right)&\text{twisting on both sides}\\
 &\cong \Ext_Q^i\left(\OO_Q, S^\vee \left(-1-l\right)\right)& \text{dualizing}\\
 &\cong \Ext_Q^i\left(\OO_Q, S \left(-l\right)\right)&\text{\cite[Theorem~2.8]{Ottaviani}}\\
 &\cong \HH^i\left(Q,S\left(-l\right)\right)&\Ext^i_Q\left(\OO_Q,\_\right)\cong \HH^i\left(Q,\_\right)\\
 &=0. &\mathllap{i \notin\left\{0,2k-1\right\}\;\text{ \cite[Theorem~2.3]{Ottaviani}}}
 \end{align*}
 Finally for $\Ext^i_Q\left(S\left(-2k+1\right),S\left(-4k+1\right)\right)$ we get
 \begin{align*}
 \Ext^i_Q\left(S\left(-2k+1\right),S\left(-4k+1\right)\right)
 &\cong \Ext^i_Q\left(S,S\left(-2k\right)\right)&\text{twisting on both sides}\\
 &\mathrlap{\cong \Ext_Q^1\left(S,S\left(-2k+i-1\right)\right)}&\mathllap{\text{Lemma~\ref{lemma selfExt spinor}}\eqref{lemma selfExt spinor case 1}}\\
 &\cong \Ext_Q^0\left(S,S\left(-2k-i\right)\right)&\mathllap{\text{Lemma~\ref{lemma selfExt spinor}}\eqref{lemma selfExt spinor case 2}}\\
 &\mathrlap{\cong \Ext_Q^{-1}\left(S,S\left(-2k+1-i\right)\right)}&\mathllap{\text{Lemma~\ref{lemma selfExt spinor}}\eqref{lemma selfExt spinor case 2}}\\
 &=0, &\mathllap{S\text{ is a sheaf}}
 \end{align*}
 where we used $i<2k-1$ and so $-2k+i\le -1$, respectively $-2k+1+i\le -1$ for the last two lines.
 
 So all the direct summands of $\Ext^i_Q\left(G,G\left(-2k\right)\right)$ vanish for $i\notin\left\{0,2k-1\right\}$ as claimed.
\end{proof}
\end{proposition}

So altogether we can now phrase the following Theorem~\ref{theorem application quadrics} which also recovers the result from \cite{Rizzardo2019} when specialized to the case $k=2$.

\begin{theorem}\label{theorem application quadrics}
Let $Q \hookrightarrow \PP^{2k}$ be the embedding of a smooth odd dimensional quadric for $k\geq2$. Then we have an exact functor:
$$\Psi_\eta :\D^b\left(Q\right)\to \D^b\left(\PP^n\right)$$ that cannot be Fourier-Mukai.
\begin{proof}
We want to apply Theorem~\ref{Main Theorem}. 

First of all we have by Proposition~\ref{proposition guaranteed non-trivial kernel} for $k>2$ an $$0\neq \eta \in \HHH^{2k+2}\left(Q,\OO_Q\left(-2k-2\right)\right)$$ that is in the kernel of $f_*:\HHH^{n+3}\left(Q,\OO\left(-2k-2\right)\right)\to \HHH^{n+3\left(\PP^{2k},f_*\OO\left(-2k-2\right)\right)}$.

For $k=2$ we get that the top Hochschild cohomology is $\HHH^{n+3}\left(Q,\OO\left(-2k-2\right)\right)$, and so by Theorem~\ref{kernels are precisely the middle part} we have $$\dim\ker\left(f_*:\HHH^{n+3}\left(Q,\OO\left(-6\right)\right)\to\HHH^{n+3}\left(\PP^4,f_*\OO\left(-6\right)\right)\right)=\h^{2,1}_1\left(Q\right).$$
Using the formula \eqref{formula middle part} we compute
\begin{align*}
\h^{2,1}_1\left(Q\right)&=\sum_{\mu=0}^5\left(-1\right)^\mu \binom{6}{\mu}\binom{-1+4-\left(\mu-1\right)\left(2-1\right)}{4}\\
&=\sum_{\mu=0}^5 \left(-1\right)^\mu \binom{6}{\mu}\binom{-1+4-\left(\mu-1\right)}{4}\\
&=\sum_{\mu=0}^5 \left(-1\right)^\mu \binom{6}{\mu}\binom{4-\mu}{4}\\
&=\left(-1\right)^0\binom{6}{0}\binom{4}{4}\\
&=1,
\end{align*}
where we used that $\binom{4-\mu}{4}$ only can be non-zero if $\mu=0$. In particular we get a one-dimensional kernel from which we may pick an $\eta \neq 0$.

We now collect the other assumptions which we verified above.

By Theorem~\ref{Theorem quadric has an exceptional sequence} $Q$ admits a tilting bundle $G$ and by Lemma~\ref{Tilting gives HT equivalence} we know that for $\Gamma :=\End\left(G\right)$ the functor $C^Q_{G,\Gamma}$ is an equivalence. In particular we get by  Proposition~\ref{cT=CT} $c_{G,\Gamma}\left(\eta\right)\neq 0$, which is assumption \eqref{assumption characteristic morphism}. Now finally we need to verify that the corresponding $\Ext$-groups are suitably concentrated, which is verified in Proposition~\ref{proposition ext concentrated in one degree} for assumption \eqref{assumption Ext concentrated in one degree} and Proposition~\ref{proposition Ext concentrated in 2 degrees} for assumption \eqref{assumption Ext vanishes in n-1 n-2}.

So we may apply Theorem~\ref{Main Theorem} to get a non-Fourier-Mukai functor $\Psi_\eta$.
\end{proof}
\end{theorem}

\appendix
\section{Modules over $\kk$-linear categories and $\Ainfty$-modules over $\Ainfty$-categories}

In this Appendix we will recall a few basic facts about modules over $\kk$-linear categories and $\Ainfty$-modules over $\Ainfty$-categories, for $\kk$ a field.

\subsection{Modules over $\kk$-linear Categories}

We start by recalling the definition of modules over a $\kk$-linear category and the relationship between those and the classical notion of modules.

The idea of generalizing the notion of modules over rings to categories first was introduced by B. Mitchell \cite{Mitchell1972}. All in all the idea is that one can interpret a $\kk$-algebra as a $\kk$-linear category with one object and under that interpretation a module corresponds to a functor from the $\kk$-linear category to the category of $\kk$-vector-spaces. 

\begin{remark}
Recall that a $\kk$-linear category $\C$ is a category such that every morphism space $\C\left(M,N\right)$ is a $\kk$-vectorspace and composition defines a $\kk$-linear map $\_\circ\_:\C\left(M',M\right)\otimes \C\left(M'',M'\right)\to \C\left( M'',M\right)$.
\end{remark}

\begin{definition}[{\cite{Mitchell1972}}]\label{definition modules over kk-linear category}
Let $\X$ be a small $\kk$-linear category. A $\X$-module is a $\kk$-linear functor
$$\M: \X \to \Vect\left(\kk\right).$$
A morphism of $\X$-modules is a natural transformation between two $\X$-modules $\cN$ and $\M$:
$$f: \cN \to \M .$$
We refer to the category of $\X$-modules $\X\modules$.
\end{definition}

\begin{lemma}
Let $\X$ be a $\kk$-linear category. Then we have that the category $\X\modules$ is a $\kk$-linear abelian category.
\begin{proof}
By Definition~\ref{definition modules over kk-linear category} we have $\X\modules = \mathrm{Fun}_\kk\left(\X,\Vect\left(\kk\right)\right)$. In particular we have immediately a canonical $\kk$-action on the morphism spaces. As kernels and cokernels can be computed objectwise in the target category \cite[A.4.3.]{Weibel1994} we have that $\mathrm{Fun}_\kk\left(\X,\Vect\left(\kk\right)\right)$ is also abelian. In particular we ger that $\X\modules$ is abelian $\kk$-linear.
\end{proof}
\end{lemma}

\begin{remark}
Let $\Gamma$ be a $\kk$-algebra. Then we have that a classically defined $\Gamma$-module $M$ consists of a $\kk$-vector space $V$ together with a $\kk$-algebra morphism $\gamma: \Gamma \to \End\left(V\right)$.

On the other hand, if we consider $\Gamma$ to be a $\kk$-linear category with one object $*$, then $\M$ consists by Definition~\ref{definition modules over kk-linear category} also of a vectorspace $V= \M\left(*\right)$ together with a morphism of $\kk$-algebras (a map of morphism spaces) $$\Gamma\to \End\left(V\right)=\Vect\left(\kk\right)\left(\M\left(*\right),\M\left(*\right)\right).$$ In particular in this case the two notions of modules over $\Gamma$ coincide.

Similarly the notion of natural transformation captures in this case precisely the commuting with the $\Gamma$ action.
\end{remark}

\begin{definition}[{\cite{Mitchell1972}}]
Let $\X $ be a $\kk$-linear category. We define the derived category of $\X$-modules (respectively bounded, bounded below or bounded above) derived category, to be the derived category (respectively bounded, bounded below or bounded above derived category) of the abelian category $\X\modules$.
$$\D^\natural\left(\X \right):=\D^\natural\left(\X\modules\right),$$
for $\natural \in \left\{\_,b,-,+\right\}$.
\end{definition}

As we will later define a $\kk$-linear category corresponding to a scheme and then model morphisms of schemes also as functors between $\kk$-linear categories we will denote the restriction of scalar functors in the following way:

\begin{definition}
Let $f:\X\to \Y$ be a $\kk$-linear functor and let $\M$ be a $\Y$-module. Then we define the module $f_*\M$ to be the $\X$-module defined by 
$$f_*\M:=\M \circ f.$$
\end{definition}

\begin{remark}
We choose the notation $f_*$ over $f^*$ as we will later model the category of sheaves on a projective scheme by modules over a $\kk$-linear category, and under this construction the functor $f_*$ corresponds to the direct image and so the notation turns out to be more consistent and less confusing throughout this work.
\end{remark} 

\begin{lemma}
Let $f:\X \to \Y$ be a $\kk$-linear functor. Then the assignment $\M\mapsto f_* \M$ defines a left exact functor $f_*:\Y\modules \mapsto \X\modules$.
\begin{proof}
As kernels and images are computed on the target category we do not need to worry about left-exactness. It also defines a functor as it is just precomposition with a functor and so it has to be functorial.
\end{proof}
\end{lemma}

\subsection{$\Ainfty$-Structures}

Throughout this section we follow \cite{kellerAinfty} and \cite{Seidel2008}, in particular we will use the sign conventions from \cite{kellerAinfty}. Although B. Keller only talks about $\Ainfty$-algebras, the sign conventions can also be applied to $\Ainfty$-categories and are equivalent to the sign conventions in the book by P. Seidel which is considering $\Ainfty$-categories throughout. Furthermore K. Lef\`evre-Hasegawa \cite{Lefevre} covers the case of $\Ainfty$-categories using the same signs as Keller, however we primarily refer to \cite{Seidel2008} for the category case, as \cite{Lefevre} is in French.

\subsection{$\Ainfty$-Categories and their Functors}

Since we will repeatedly use dg-categories as examples for $\Ainfty$-categories we recall the definition of a dg-category

\begin{definition}
A dg category $\C$ is a category such that we have for all $M,N \in \C$ a chain complex $\C^* \left( M,N\right)$, such that the Leibnitz rule holds $$ d \left( x \circ y \right) =d x \circ y + x \circ d y .$$

\end{definition}

\begin{definition}[{\cite[3.1.]{kellerAinfty}}]\label{Definition An category}
Let $n\in \N\cup \left\{\infty\right\}$. An $\An$-category $\X$ over a field $\kk$ consists of a class of objects $\obj \left(\X\right)$ and $\Z$-graded $\kk$-vector-spaces as morphism spaces 
$$\X\left(a,b\right), $$
for $a,b \in \obj\left(\X\right)$, together with compositions 
\begin{equation*}
\m_i: \underbrace{\X\left(a_i,a_{i-1}\right)\otimes_\kk \X \left(a_{i-1},a_{i-2}\right)\otimes ... \otimes \X \left(a_1,a_0\right)}_i \to \X\left(a_i,a_0\right)
\end{equation*}
of degree $2-i$ for $1\le i\le n$ and $a_0,...,a_i \in \obj\left(\X\right)$ such that 
\begin{equation}\label{equation A-infinity category}
\sum_{r+s+t=k} \left(-1\right)^{r+st}\m_u\circ  \left(\id^{\otimes r} \otimes \m_s \otimes \id^{\otimes t}\right)=0 \tag{$*_k$}
\end{equation}
holds for all $k\le n$, where $u=r+1+t$.

We will sometimes denote $a\in \obj \left(\X\right)$ by $a \in \X$ to avoid clumsy notation.
\end{definition}

\begin{definition}[{\cite[(2a)]{Seidel2008}}]
An $\A_n$-category $\X$ is called unital if every object $a\in \obj\left( \X \right)$ admits a unit $\id\in \X\left(a,a\right)^{0}$ such that 
\begin{align*}
&\m_1\left(\id\right)=0\\
&\m_2\left(x,\id\right)=x=\m_2\left(\id,x\right )\\
&\m_i \left(x_i,...,\id,....,x_1\right)=0 & i \neq 2.
\end{align*}
\end{definition}

\begin{remark}\label{remark equations A-infinity category}
Observe that the first few incarnations of  $\eqref{equation A-infinity category}$  give:
\begin{itemize}
\item[$k=1$:] In this case $\left(*_1\right)$ gives
$$\m_1 \circ \m_1=0.$$
This means that $\m_1$ defines a differential on $\X\left(a,b\right)$.
\item[$k=2$:] Here $\left(*_2\right) $ boils down to 
$$\m_1 \circ \m_2 = \m_2\left(\m_1 \circ \id +\id \circ \m_1\right),$$
which is the Leibnitz rule $d \left(x\circ y\right)=d x \circ y + x\circ d y$.
\item[$k=3$:] And $\left(*_3 \right)$ gives
\begin{align*}
&\m_2\circ \left( \id \otimes \m_2 - \m_2 \otimes \id \right)= \\
&=\m_1 \circ \m_3 + \m_3 \otimes \left( \m_1\otimes \id \otimes \id + \id \otimes \m_1 \otimes \id + \id \otimes \id \otimes \m_1 \right),
\end{align*}
which means that $\m_2$ is associative up to a homotopy given by $\m_3$. More generally one can think of an $\An$-category as a category that is homotopy-associative up to degree $n$. 
\end{itemize}
By Definition~\ref{Definition An category} every $\An$ category defines an $\Am$-category for all $m\le n$ just by forgetting the higher actions.
\end{remark}

\begin{definition}[{\cite[(1a)]{Seidel2008}}]
Let $\X$ be an $\A_n$-category for $n\geq 3$. Then the category $\HH^*\left(\X\right) $ is the graded $\kk$-linear category consisting of the same objects as $\X$ and morphism spaces
$$\HH^*\left(\X\right)\left(a,b\right):=\HH^*\left(\X\left(a,b\right)\right).$$
Where we use Remark~\ref{remark equations A-infinity category} to consider $\X\left(a,b\right)$ as a chain complex with differential $\m_1$.

The $\kk$-linear category $\HH^0\left(\X\right)$ is the category with the same objects as $\X$ and morphism spaces 
$$\HH^0\left(\X\right)\left(a,b\right):=\HH^0\left(\X\left(a,b\right)\right).$$

We have by Remark~\ref{remark equations A-infinity category} that $\HH^*\left(\X\right)$ defines a graded $\kk$-linear category and $\HH^0\left(\X\right)$ defines an ordinary $\kk$-linear category.
\end{definition}

\begin{definition}[{\cite[(2a)]{Seidel2008}}]
An $\A_n$-category is called homologically unital if $\HH^0\left(\X\right)$ admits a unit morphism $\id\in \HH^0\left(\X\right)\left(a,a\right)$ for all $a \in \obj \left(\X\right)$.
\end{definition}

\begin{definition}
An $\An$-category is called small if its objects form a set. It is called essentially small if the isomorphism classes of objects form a set.
\end{definition}

\begin{definition}[{\cite[3.1.]{kellerAinfty}}]
An $\An$-category is an $\An$-algebra if $\obj\left(\X\right)$ consists of only one object for $n \in \N\cup \left\{\infty \right\}$.
\end{definition}

\begin{example}
There are a few obvious examples of $\Ainfty$-categories:
\begin{itemize}
\item Let $\X$ be a $\kk$-linear category, then it is an $\Ainfty$ category via
$$\m_i =\begin{cases} \left(\_\right)\circ \left(\_\right)&i=2\\
0 & i\neq 2.\end{cases}$$
\item More generally, let $\X$ be a dg-category, then $\X$ is an $\Ainfty$-category with
$$\m_i=\begin{cases} 
 d  &i=1 \\
\left(\_\right)\circ \left(\_\right) & i=2\\
0& i\notin \left\{1,2\right\}. \end{cases} $$
\end{itemize}
\end{example}

\begin{definition}[{\cite[3.4.]{kellerAinfty}}]
An $\An$-functor between two $\An$-categories $f:\X \to \Y$ is given by a map on objects
$$f: \obj\left(\X\right) \to \obj\left(\Y\right)$$
and a set of morphisms
$$\left\{f_i: \X \left(a_i,a_{i-1}\right)\otimes \X \left(a_{i-1},a_{i-2}\right)\otimes...\otimes \X\left(a_1,a_2\right)\to \Y\left(f\left(a_i\right),f\left(a_0\right)\right)\right\}$$
of degree $1-i$ for every $i\le n$ and $a_i,...,a_0 \in \obj\left(\X\right)$ such that
\begin{equation}\label{equation A-infinity morphisms}
\sum_{\mathclap{r+s+t=k}}\left(-1\right)^{r+st}f_u\left(\id^{\otimes r}\otimes \m_s \otimes \id^{\otimes t}\right)=\sum_{\mathclap{\stackrel{1\le l\le n}{ k=i_1+...+i_l}}}\left(-1\right)^m \m_r \left(f_{i_1}\otimes f_{i_2}\otimes...\otimes f_{i_l}\right)\tag{$**_k$}
\end{equation}
holds, where $u=r+1+t$ and 
$$m=\left(l-1\right)\left(i_1-1\right)+\left(l-2\right)\left(i_2-1\right)+...+2\left(i_{l-2}-1\right)+\left(i_{l-1}-1\right).$$
\end{definition}

\begin{remark}
Again we compute the first few incarnations of $\eqref{equation A-infinity morphisms}$:
\begin{itemize}
\item[$k=1$:] In this case we have $$f_1 \circ \m_1=m_1\circ f_1,$$ in particular $f_1$ defines a morphism of chain complexes.
\item[$k=2$:] Here we get 
$$f_1 \circ \m_2 = \m_2 \circ \left(f_1 \otimes f_1\right)+\m_1 \circ f_2 + f_2 \left(\m_1\otimes \id +\id \otimes \m_1\right),$$
so $f_1$ commutes with $\m_2$ up to a homotopy given by $f_2$. 
\end{itemize}
More generally one can think of an $\An$-morphism $f$ as commuting with the $\An$-structure up to higher homotopies, whose information $f$ includes in form of the higher $f_i$.
\end{remark}

\begin{definition}[{\cite[3.1.]{kellerAinfty}}]
An $\An$-functor between two unital $\An$-algebras is called an $\An$-morphism for $n \in \N\cup \left\{\infty \right\}$.
\end{definition}

\begin{definition}[{\cite[3.1.]{kellerAinfty}}] 
An $\Ainfty$-functor $f:\A\to \B$ is a quasi-equivalence if 
$$f:\obj\left(\A\right)/_{\mathord{\cong}} \to \obj\left(\B\right)/_{\mathord{\cong}}$$
is surjective and all $f_1$ induce isomorphisms on cohomology
$$\HH^*\left(f_1\right):\HH^*\left(\A\left(a,a'\right)\right)\xrightarrow{\sim} \HH^*\left(\B\left(f a, fa'\right)\right).$$
\end{definition}

\begin{proposition}[{\cite[Proposition~3.2.1]{Lefevre}}]
Every homologically unital $\Ainfty$-category is quasi-equivalent to an unital one.
\end{proposition}

\begin{definition}[{\cite[3.4.]{kellerAinfty}}]
A quasi-equivalence between two $\Ainfty$-algebras is called a quasi-isomorphism.
\end{definition}

\begin{theorem}[{\cite{Kadeishvili1980}}]\label{theorem Kadeishvili}
Let $\X$ be an $\Ainfty$-category. Then the cohomology $\HH^*\left(\X\right)$ has an $\Ainfty$-category structure such that 
\begin{itemize}
\item $\m_1=0$
\item there is a quasi-equivalence $\HH^*\X\xrightarrow{\sim} \X$ lifting the identity on $\HH^*\X$.
\end{itemize}
Moreover, this structure is unique up to (non-unique) isomorphism of $\Ainfty$-categories.

\end{theorem}

\begin{remark}
From now on we will assume that the cohomology $\HH^*\left(\X\right)$ of an $\Ainfty$-category is equipped with the $\Ainfty$-structure arising by Theorem~\ref{theorem Kadeishvili} instead of just regarding it as a graded category interpreted as an $\Ainfty$-category. The $\Ainfty$-category constructed in Theorem~\ref{theorem Kadeishvili} is also referred to as the minimal $\Ainfty$-model of $\X$.
\end{remark}

\subsection{$\Ainfty$-Modules and their Functors}

\begin{definition}[{\cite[4.2.]{kellerAinfty}}]
Let $\X$ be a small $\An$-category for $n\in \N\cup \left\{\infty\right\}$. An $\A_n$-module over $\X$ consists of a $\Z$-graded space
$$\M\left(a,b\right)$$
for every pair of objects $a,b \in \obj \X$ and higher composition morphisms
$$\m_i: \underbrace{\M\left(a_i,a_{i-1}\right)\otimes \X\left(a_{i-1},a_{i-2}\right) \otimes ... \otimes \X\left(a_1,a_0\right)}_{i}\to \M\left(a_i,a_0\right)$$
of degree $2-i$ such that the following equation holds
\begin{equation}\label{equation A-infinity Module}
\sum_{r+s+t=k} \left(-1\right)^{r+st}\m_u\circ  \left(\id^{\otimes r} \otimes \m_s \otimes \id^{\otimes t}\right)=0, \tag{$**_k$}
\end{equation}
where depending on the input $\m_i$ needs to be considered as the $i$th higher composition morphism of $\X$ or $\M$.
\end{definition}

\begin{remark}
We again compute a few incarnations of $\eqref{equation A-infinity Module}$ to give some intuition on the modelled structure.
\begin{itemize}
\item[$k=1$:] In this case we get $$\m_1^\M\circ \m_1^\M=0.$$ So $\m_1$ defines a differential.
\item[$k=2$:] Here we get $$\m_1^\M \circ \m_2^\M = \m^\M_2 \circ \left(\m^\M_1 \otimes \id_\M+\id_\M\otimes \m_1^\A\right),$$
which means that $\m_2$ suffices the Leibnitz rule.
\item[$k=3$:] For this we get similar to the $\Ainfty$-algebra case that the action of $\M$ induced by $\m_2$ is associative up to a homotopy, which is given by $\m_3$. 
\end{itemize}
So one can think about an $\Ainfty$-module as a homotopy coherent module over $\X$.
\end{remark}
 
 \begin{example}
 We collect once more the standard examples.
 \begin{itemize}
 \item Let $\M$ be a graded module over a $\kk$-linear category $\X$, then it is an $\Ainfty$-module over $\X$ via 
 $$\m_i=\begin{cases}\left(\_\right)\circ \left(\_\right) & i=2\\
 0 & i\neq 2. \end{cases}$$
 \item Let $\M$ be a dg-module over a dg-algebra $\X$. Then it defines an $\Ainfty	$-module over $\X$ via
 $$\m_i= \begin{cases}d_\M &i=1\\
 \left(\_\right)\circ \left(\_\right)& i=2\\
 0 & i \notin \left\{1,2\right\}.\end{cases}$$
 \end{itemize} 
 \end{example}
 
 \begin{definition}[{\cite[4.2.]{kellerAinfty}}]
 Let $\M, \cN$ be $\An$-modules over an $\An$-category $\X$ for $n\in \N \cup \left\{\infty\right\}$. A morphism of $\An$-modules consists of a set of morphisms:
 $$f_i: \underbrace{\M\left(a_i,a_{i-1}\right)\otimes \X\left(a_{i-1},a_{i-2}\right)\otimes ...\otimes \X\left(a_1,a_0\right)}_i  \to \cN\left({a_{i},a_0}\right)$$
 of degree $1-i$ for $i\le n$, such that we have for every $k<n$
 \begin{equation}\label{formula A-infinity module homomorphism}
 \sum_{\mathclap{r+s+t}}\left(-1\right)^{r+st}f_u \circ \left(\id^{\otimes r}\otimes \m_s \otimes \id^{\otimes t}\right)=\sum_{\mathclap{n=r+s}}\left(-1\right)^{\left(r-1\right)s}\m_{u'}\left(f_r \otimes \id^s\right)\tag{$**_k$}
 \end{equation}
 where $u=r+s+t$ and $u'=1+s$.
 \end{definition}
 
 \begin{example}
 We compute again $\eqref{formula A-infinity module homomorphism}$ for small $k$:
 \begin{itemize}
 \item[$k=1$:] Similar to the cases above $\left(**_1\right) $ boils down to $$f_1 \circ \m_1 =\m_1 \circ f_1,$$
 which means that $f_1$ defines a morphism of chain complexes.
 \item[$k=2$:] Here we get 
 $$f_1\circ \m_2 - f_2 \circ \left(\m_1 \otimes \id + \id \otimes \m_1\right)= \m_2 \circ \left(f_1 \otimes  \id_\X \right) + \m_1 \circ f_2.$$
This means  that similarly to the case of an $\A_n$-functor between $\A_n$-categories the equation $\left(**_2\right)$ encodes that $f_1$ is compatible with the action induced by $\m_2$ up to a homotopy given by $f_2$.
\end{itemize} 
 These examples are another reason one can think about $\Ainfty$-structure as a notion for inductive homotopy coherent algebraic structures. 
 \end{example}
 
 \begin{definition}[{\cite[4.2.]{kellerAinfty}}]
 An $\An$-morphism $f: \M \to \cN$ is a quasi-isomorphism if it induces an isomorphism on cohomology
 $$\HH^*\left(f\right): \HH^*\M \xrightarrow{\sim} \HH^*\cN.$$
 \end{definition}

 \begin{definition}[{\cite[4.2.]{kellerAinfty}}]
 Let $f: \M \to \M'$ and $g: \M' \to \M''$ be morphisms of $\Ainfty$-modules over a homologically unital $\Ainfty$-algebra $\X$. Then the composition $f\circ g: \M \to \M''$ is given by
 $$\left(f\circ g\right)_n= \sum_{\mathclap{n=r+s}} \left(-1\right)^{\left(r-1\right)s} f_u\left(g_r\otimes \id^{\otimes s}\right),$$
 where we put $u=1-s$.
 \end{definition}

\begin{definition}[{\cite[4.2.]{kellerAinfty}}]
Let $\X$ be a homologically unital $\Ainfty$-algebra, then we define the category of $\Ainfty$-modules $\C_\infty\left(\X\right) $ to be the category consisting of $\Ainfty$-modules and morphisms given by $\Ainfty$-morphisms.
\end{definition}

\begin{remark}
The identity of an object in $\C_\infty\left(\X\right)$ is given by $$\id=\left(\id,0,...\right).$$
\end{remark}
 
\begin{definition}[{\cite[1k]{Seidel2008}}]
Let $f:\X \to \Y$ be an $\A_i$-functor. Then the functor 
$$f_*:\C_\infty \left(\Y\right)\to \C_\infty\left(\X\right)$$ is given on modules by
$$f_*\M\left(a\right):= \M\left(f\left(a\right)\right)$$
for objects $a \in \obj\left(\X\right)$. Higher compositions are given by 
$$\m_k\left(m,x_{k-1},...,x_1\right)=\sum_{l<k} \sum_{\mathclap{s_1,...,s_l}} \m_l\left(m,f_{s_l}\left(x_{k-1},...,x_{k-s_l}\right),...,f_{s_1}\left(a_{s_1},...,a_1\right)\right) .$$
On morphisms $f^*$ is given by 
$$f_* \varphi_k\left(m,x_{k-1},...,x_1\right)=\sum_{l<k} \sum_{s_1,...,s_l} \varphi_l\left(m,f_{s_l}\left(x_{k-1},...,x_{k-s_l}\right),...,f_{s_1}\left(a_{s_1},...,a_1\right)\right) .$$
\end{definition}

\begin{remark}
We again choose the notation $f_*$ over $f^*$ as we will later model the category of sheaves on a projective scheme by modules over a $\kk$-linear category, and under this construction the functor $f_*$ corresponds to the direct image and so the notation turns out to be more consistent and less confusing throughout this work.
\end{remark}

\begin{definition}[{\cite[4.2.]{kellerAinfty}}]
Let $\X$ be a homologically unital small $\Ainfty$-category. Then we define the category
$$\D_\infty \left(\X \right):=\C_\infty\left(\X\right)[\left\{\Ainfty-\text{quasi-isomorphism}\right\}^{-1}].$$
\end{definition} 
 
\begin{remark}[{\cite[4.2.]{kellerAinfty}}]
More generally one could consider $\Ainfty$-categories over commutative rings instead of a field $\kk$. In this case we would have to distinguish between the derived category of $\Ainfty$-modules, as we defined it, and the category of $\Ainfty$-modules up to homotopy. However, over a field one can prove that actually every quasi-isomorphism of $\Ainfty$-modules is a homotopy equivalence and vice versa. In particular in this case the naively derived category arising by formally inverting quasi-isomorphisms and the category of $\Ainfty$-modules up to homotopy coincide.

The interpretation of $\D_\infty\left(\X\right)$ as arising via $\Ainfty$-modules up to homotopy immediately gives that $\D_\infty\left(\X\right)$ is well-defined and there are no set-theoretic issues arising.
\end{remark} 

 \bibliographystyle{hep}
\bibliography{Alexandria}  

\end{document}